\documentclass
[
    a4paper,
    DIV=10,
    abstracton
]
{scrartcl}

\usepackage{amsfonts,amssymb,amsmath,amsthm}
\usepackage{fullpage}
\usepackage{bm}
\usepackage{bbm}
\usepackage{dsfont}
\usepackage[pdftex]{graphicx}
\usepackage{xcolor}
\usepackage{epstopdf}
\usepackage{enumerate}  %for custom lists
\usepackage{tikz} %for nice pictures
\usetikzlibrary{decorations.markings} %for halfway arrows
\usepackage[colorlinks,
			pdffitwindow=false,
      plainpages=false,
      pdfpagelabels=true,
     	pdfpagemode=UseOutlines,
      pdfpagelayout=SinglePage,
			colorlinks=true,
      hyperfootnotes=false,
			linkcolor=blue,
			citecolor=green!50!black]{hyperref}

\usepackage{nccmath}   %for mfrac; medium-sized fractions
\usepackage{algorithm}
\usepackage[noend]{algpseudocode}

\graphicspath{{./graphics/}}

\newtheoremstyle{custom}{3pt}{3pt}{}{}{\bfseries}{:}{.5em}{}
\theoremstyle{custom}
\newtheorem{example}    {Example}%[section]

\newtheorem{remark}[example]{Remark}
\numberwithin{example}{section}
\newcommand{\A}{\mathcal{A}}
\newcommand{\e}{\mathrm{e}}
\newcommand{\ep}{\varepsilon}
\newcommand{\fw}{\mathrm{fw}}
\newcommand{\bw}{\mathrm{bw}}
\newcommand{\cross}{\mathrm{cross}}
\newcommand{\meet}{\mathrm{meet}}
\newcommand{\R}{\mathbb{R}}

\newcommand{\J}{\mathcal{J}}
\newcommand{\N}{\mathbb{N}}

\newcommand{\suber}[1]{_{\scriptscriptstyle{(#1)}}}
\newcommand{\super}[1]{^{\scriptscriptstyle{(#1)}}}
\newcommand{\smallto}{{\scriptstyle \rightarrow}}
\newcommand{\rv}[1]{\bm{#1}}		% random variables are bold
\def \prob {\mathbb{P}} % probability sign (measure)

\newcommand{\dist}{\text{dist}} % distance function in algorithm
\newcommand{\spath}{\text{path}} % path function in algorithm

\DeclareMathOperator\supp{\mathrm{supp}}
\DeclareMathOperator*\argmax{\arg\,\max}

\title{From large deviations to semidistances of transport and mixing: coherence analysis for finite Lagrangian data}
\author{P\'eter Koltai\thanks{Institute of Mathematics, Freie Universit\"at Berlin, Arnimallee 6, 14195 Berlin, Germany. E-mail: peter.koltai{@}fu-berlin.de}
          \and 
        D.R. Michiel Renger\thanks{Weierstra\ss-Institut, Mohrenstra{\ss}e 39, 10117 Berlin, Germany. E-Mail: d.r.michiel.renger@wias-berlin.de}}
\date{}

\begin{document}
\maketitle

\begin{abstract}

One way to analyze complicated non-autonomous flows is through trying to understand their transport behavior. In a quantitative, set-oriented approach to transport and mixing, finite time coherent sets play an important role. These are time-parametrized families of sets with unlikely transport to and from their surroundings under small or vanishing random perturbations of the dynamics. Here we propose, as a measure of transport and mixing for purely advective (i.e., deterministic) flows, (semi)distances that arise under vanishing perturbations in the sense of large deviations.
Analogously, for given finite Lagrangian trajectory data we derive a discrete-time and space semidistance that comes from the ``best'' approximation of the randomly perturbed process conditioned on this limited information of the deterministic flow. It can be computed as shortest path in a graph with time-dependent weights.
Furthermore, we argue that coherent sets are regions of maximal farness in terms of transport and mixing, hence they occur as extremal regions on a spanning structure of the state space under this semidistance---in fact, under any distance measure arising from the physical notion of transport. Based on this notion we develop a tool to analyze the state space (or the finite trajectory data at hand) and identify coherent regions. We validate our approach on idealized prototypical examples and well-studied standard cases.

\end{abstract}

%\tableofcontents

%%%

\section{Introduction}

\paragraph{Transport in dynamical systems.}

Instrumental to understanding the essential behavior of complicated non-autonomous flows is to grasp how transport is happening in them. This leads on a qualitative level to objects that prohibit transport, commonly named \emph{transport barriers}; often originating from the geometric picture for autonomous systems and that trajectories are unable to cross co-dimension 1 invariant manifolds~\cite{Hal00,Hal01,HaBV12,HaBV13}. For periodically-forced systems, invariant manifolds enclose regions called ``lobes'' that get transported across these periodically varying manifolds~\cite{MMP84,RKW90}.

On a quantitative level, one searches for surfaces of small flux~\cite{BaFrSa14,Kar16,FrKo17}, so-called partial barriers~\cite{Wig37,Mei92}. Instead of characterizing regions that do not mix with one another via enclosing them by boundaries of low flux, there are approaches that aim to describe these sets directly. Such \emph{set-oriented} concepts are strongly interwoven with the theory of transfer operators (Perron--Frobenius and Koopman operators), and comprise \emph{almost-invariant sets}~\cite{DeJu99}, \emph{ergodic partitions}~\cite{MeWi99} in autonomous, and \emph{coherent sets}~\cite{FrSaMo10,Froyland2013} in the non-autonomous cases.

Distinctive attention has been given to coherent sets, which are a (possibly time-dependent) family of sets having little or no exchange with their surrounding in terms of transport, and are robust to small diffusion over a finite time of consideration~\cite{Froyland2013,Froyland2015}. Natural examples include moving vortices in atmospheric~\cite{RypEtAl07,FrSaMo10}, oceanographic~\cite{TBBM03,DFHPG09,FrEtAl15}, and plasma flows~\cite{PHJJ07}. In such applications, one would like to be able to find coherent sets even in the cases when a dynamical model that can be evaluated arbitrarily often is not available, only a finite set of \emph{Lagrangian trajectory data} (passive tracers moving with flow with positions sampled at discrete time instances). This problem has received lot of attention in recent years, and a diverse collection of tools has been developed to tackle it~\cite{BuMe12,AlTh12,SGetal15,FrPa15,AlPe15, WRR15,HaEtAl16,BaKo17,SchDa17,PGSc17, RyPr17,FaMePo16,FrJu17}.

While other current methods aim at collecting trajectories into coherent sets, in~\cite{BaKo17} it has been proposed to go one step further and analyze the connectivity structure of  the state space under transport and mixing with ``transport coordinates'' and the ``skeleton of transport''. Very similar observations have been made earlier in~\cite{BuMe12} in the infinite-time limit for periodically-forced systems. While coherent sets (and transport barriers) aim at \emph{partitioning} the state space, the skeleton is aiming at \emph{``spanning''} the state space with respect to transport. In this respect, coherent sets can be associated with distinct ``extremal regions'' of the skeleton. Here we will only use this idea of extremality, more precisely that coherent sets are ``maximally far'' from one another, as measured by transport. To this end we will need to measure ``farness'' of dynamical trajectories.

Several \emph{dynamical distance measures} have been put forward already to measure the ``distance'' or ``dissimilarity'' of trajectories or initial states in dynamical systems~\cite{MeBa04,BuMe12,FrPa15,HaEtAl16,FaMePo16, KarraschKeller2016}. The majority of them are shown to serve their purpose well in revealing coherent structures efficiently and reliably. However, they are either heuristic in the sense that they are not derived from the physical notion of transport and mixing, or no discretizations to finite scattered trajectory data have been developed.

The purpose of this paper is thus twofold. On the one hand, we develop a distance measure (a semidistance) between trajectories that is derived from the physical notion of transport and mixing subject to diffusion of vanishing strength, and we also derive a discretized distance measure for finite (also possibly sparse and incomplete) Lagrangian data that is consistent with its continuous counterpart in the limit of infinite data.
On the other hand, we construct a tool to analyze with such distances the structure of the state space under transport, especially to find coherent sets. This tool makes use of the idea that coherent sets are some sort of extremal regions on a spanning structure with respect to transport, although in this work we will not investigate this ``skeleton'' in its entirety.

\paragraph{Finite time coherent sets.}

Let us consider the ordinary differential equation (ODE)
\begin{equation}
  \dot x_t = v(t,x_t)
\label{eq:ODE}
\end{equation}
on some bounded~$X\subset\R^d$ and on a finite time interval~$[0,T]$ for some~$T>0$. Throughout the paper we will assume that~$v:[0,T]\times X\to\R^d$ is a continuous velocity field tangential at the boundary, such that the flow of~\eqref{eq:ODE}, denoted by~$\phi_{s,t}[\cdot]$, $0\le s, t\le T$, is a diffeomorphism on appropriate subsets of~$X$. For~$t<s$ we flow backward in time:~$\phi_{s,t}= \phi_{t,s}^{-1}$.

Many different notions to characterise coherent sets have been proposed in the literature. Central to all of these notions is the idea that coherent sets should be \emph{robust under noise}; without such a requirement any non-intersecting characteristic of a singleton could be considered a coherent set. To this end one typically perturbs the ODE~\eqref{eq:ODE} by a random noise~\cite{DeMaJu16,FrKo17,KarraschKeller2016}, leading to the It{\^o} stochastic differential equation (SDE)\footnote{We denote random variables by boldface symbols.}
\begin{equation}
  d\rv{x}_t\super{\ep} = v(t,\rv{x}_t\super{\ep})dt + \sqrt\ep d\rv{w}_t\,,
\label{eq:SDE}
\end{equation}
where~$\{\rv{w}_t\}_{t\in [0,T]}$ is a Wiener process (Brownian motion) with generator~$\phi\mapsto\frac12\Delta\phi$, reflecting boundaries, and starting from $\rv{w}_0=0$ (deterministically) and $\ep>0$ is, at least for now, a given small constant. In fact, the rigorous mathematical formulation of an SDE with reflecting boundaries can be quite subtle, see \cite{Andres2009}. We ignore this issue as it does not affect our analysis.

According to the definition of \emph{finite time coherent pairs}~\cite{FrSaMo10,Froyland2013,KoltaiCiccottiSchuette2016}, two sets~$A,B\subset X$ are coherent for times~$0$ and~$T$ if most mass from set~$A$ is likely to end up in set~$B$, and most mass ending up in set~$B$ is likely to originate from set~$A$, that is,
\begin{align}
  \prob\big[\rv{x}\super{\ep}_T\in B \mid \rv{x}\super{\ep}_0\in A\big] \approx 1, &&\text{ and }&& \prob\big[\rv{x}\super{\ep}_0\in A \mid \rv{x}\super{\ep}_T\in B\big]\approx 1.
\label{eq:forward backward conditions}
\end{align}
Naturally, for practical purposes one would need to choose how small~$\ep$ and how large these probabilities should be. As the systems we are dealing with are often deterministic by nature, and there is no ``physically straightforward'' choice of the diffusion strength~$\ep$, our first aim is to remove some of this indeterminacy by quantifying what it means for probabilities to be close to~$1$ for small~$\ep$, in terms of \emph{large deviations} as we explain below\footnote{A different way of factoring out diffusion to obtain coherent sets for deterministic flows appeared in the set-oriented transfer-operator based characterization in~\cite{Froyland2015,FrKw16}, leading to the notion of the \emph{dynamic Laplacian}. See also our concluding remarks in Section~\ref{ssec:dynLap}.}. However, it turns out that the forward and backward conditions~\eqref{eq:forward backward conditions} are essentially equivalent in the large-deviation regime, 
%see Appendix~\ref{app:two-time LDP}. Therefore they do not need to be treated separately, and we will consider ``forward probabilities'' only. Moreover,
and even worse, the large-deviation limits of~\eqref{eq:forward backward conditions} hardly give any quantitative information about how coherent two sets might be, as discussed in Appendix~\ref{app:two-time LDP}. To conclude, the large deviations of conditions~\eqref{eq:forward backward conditions} do not yield sensible conditions for coherence.

\paragraph{Large-deviation based semidistances.}

In the current paper we take a different approach. We study semidistances that quantify how unlikely it is for mass to flow from one point to another. These are semidistances in the sense that they satisfy all properties of a metric except for the triangle inequality. In the first part of this paper, sections~\ref{sec:cont setting} and \ref{sec:discr setting}, we show how such semidistances can arise naturally from probabilistic arguments via large-deviation principles, as we explain below. In the second part, sections~\ref{sec:method} and \ref{sec:numerics}, we discuss how (general) semidistances can be used to analyse coherent sets, and we apply the concepts of this paper to a number of examples.

In Section~\ref{sec:cont setting} we derive two different semidistances from the \emph{large deviations} of two probabilities. The first one is related to the probability that the endpoint~$\rv{x}_T\super{\ep}$ of the random path is~$\phi_{0,T}\lbrack y\rbrack$, given that it starts in $\rv{x}_0\super{\ep}=x$, for any two initial positions $x,y\in X$. As $\ep\to0$, the process can no longer deviate from the deterministic flow of~\eqref{eq:ODE}, and hence this probability will converge to $0$ whenever $x\neq y$. In fact, it converges exponentially fast~\cite{FreidlinWentzell1998}, i.e.,
\begin{equation}
  \prob\big[ \rv{x}_T\super{\ep}\asymp \phi_{0,T}\lbrack y\rbrack \mid \rv{x}_0\super{\ep}=x\big] \sim \e^{-\frac1\ep \mu_T(x\smallto y)}
\label{eq:cont ldp one traj}
\end{equation}
for some function $\mu_T(x\smallto y)\geq0$, where this statement and notation are made precise in Section~\ref{subsec:cont traj ldp}. Such exponential convergence results are called \emph{large-deviation principles}, and $\mu_T(x\to y)$ is the \emph{large-deviation rate}. The less probable it is to reach one point from another, the larger the rate between them. The first semidistance is then obtained via symmetrisation:
\begin{equation}
  \mu_T^\cross(x,y):=\mu_T(x\smallto y) + \mu_T(y\smallto x).
\label{eq:cont fwbw distance}
\end{equation}
We call this the \emph{cross} semidistance, since it arises from mass flowing from $x$ to $\phi_{0,T}\lbrack y\rbrack$ and mass flowing from $y$ to $\phi_{0,T}\lbrack x\rbrack$ simultaneously and independently, see Figure~\ref{fig:cont symm dist}.

The second semidistance arises as the large deviations of the probability for two independent random trajectories $\rv{x}\super\ep,\rv{y}\super\ep$ starting at $x$ and $y$, respectively, to meet at or before time~$T$ (Figure~\ref{fig:cont meet dist}): 
\begin{equation}
  \prob\big[ \rv{x}_T\super{\ep}\asymp \rv{y}_T\super{\ep} \mid \rv{x}_0\super{\ep}=x, \rv{y}_0\super{\ep}=y\big] \sim \e^{-\frac1\ep \mu_T^\meet(x,y)},
\label{eq:cont ldp meet}
\end{equation}
where the \emph{meeting} semidistance is given by
\begin{equation*}
  \mu_T^\meet(x,y):=\inf_{z\in X} \mu_T(x\smallto z) + \mu_T(y\smallto z).
\end{equation*}

\begin{figure}[h!]
\begin{minipage}{0.45\textwidth}
\centering
\begin{tikzpicture}[scale=1.1]
\begin{scope}[decoration={markings,mark=at position 0.5 with {\arrow{>}}}] 
  \tikzstyle{every node}=[font=\small]
  \draw[->](0,0)--(4.5,0) node[anchor=west]{$t$};
  \draw(4,0.05)--(4,-0.05) node[anchor=north]{$T$};
  \draw(0,0)--(0,2.5);
  \filldraw(0,1) circle (0.05) node[anchor=east]{$x$};
  \draw[dashed](0,1).. controls (2,1.1) and (3.5,0.8) .. (4,0.5);
  \filldraw(4,0.5) circle (0.05) node[anchor=west]{$\phi_{0,T}\lbrack x\rbrack$};
  \filldraw(0,2) circle (0.05) node[anchor=east]{$y$} ;
  \draw[dashed](0,2).. controls (1,1.5) and (3,1.3) .. (4,1.5);
  \filldraw(4,1.5) circle (0.05) node[anchor=west]{$\phi_{0,T}\lbrack y\rbrack$};
  \draw(0,1)[postaction={decorate}].. controls (1.5,1.1) and (3.5,1.2) .. (4,1.5);
  \draw(0,2)[postaction={decorate}].. controls (1,1.4) and (3,1.2) .. (4,0.5);
\end{scope}
\end{tikzpicture}
\caption{$\mu_T^\cross(x,y)$ is the cost to move from $x$ to $\phi_{0,T}\lbrack y\rbrack$ and from $y$ to $\phi_{0,T}\lbrack x\rbrack$.}
\label{fig:cont symm dist}
\end{minipage}
\qquad
\begin{minipage}{0.45\textwidth}
\centering
\begin{tikzpicture}[scale=1.1]
\begin{scope}[decoration={markings,mark=at position 0.5 with {\arrow{>}}}]
  \tikzstyle{every node}=[font=\small]
  \draw[->](0,0)--(4.5,0) node[anchor=west]{$t$};
  \draw(4,0.05)--(4,-0.05) node[anchor=north]{$T$};
  \draw(0,0)--(0,2.5);
  \filldraw(0,1) circle (0.05) node[anchor=east]{$x$};
  \draw[dashed](0,1).. controls (2,1.1) and (3.5,0.8) .. (4,0.5);
  \filldraw(4,0.5) circle (0.05) node[anchor=west]{$\phi_{0,T}\lbrack x\rbrack$};
  \filldraw(0,2) circle (0.05) node[anchor=east]{$y$} ;
  \draw[dashed](0,2).. controls (1,1.5) and (3,1.3) .. (4,1.5);
  \filldraw(4,1.5) circle (0.05) node[anchor=west]{$\phi_{0,T}\lbrack y\rbrack$};
  \draw(0,1)[postaction={decorate}].. controls (1.5,1.3) and (3.5,1) .. (4,1);
  \draw(0,2)[postaction={decorate}].. controls (1,1.4) and (3,1.2) .. (4,1);
\end{scope}
\end{tikzpicture}
\caption{$\mu_T^\meet(x,y)$ is the cost for two trajectories to meet.}
\label{fig:cont meet dist}
\end{minipage}
\end{figure}

By this procedure we find two semidistances~$\mu_T^\cross$ and~$\mu_T^\meet$ that can be used as a measure of ``farness'' of points~$x,y$, which will be low for points in the same coherent set, and high otherwise. Since both arise from large-deviation principles, they have a nice additional interpretation as a probabilistic cost or free energy that needs to be paid in order to deviate from the expected flows; such interpretation is common in statistical physics, see for example~\cite{Onsager1953I}.

Nevertheless, we will see that in order to calculate these costs explicitly,
%one would need to differentiate in path space with respect to the exact flow maps~$\phi_{0,T}$, which can be hard to calculate in practice. This becomes even more so problematic as the exact flow map $\phi_{0,T}$ nor the velocity field $v$ may be known.
the velocity field~$v$ needs to be known. As discussed above, this is in practice seldom the case; mostly one can only assume to have discrete-time snapshots of the positions of a limited number of floaters.
%Indeed, in practice one often uses a large number of floaters or trackers to measure the flow; the resulting data set consists of discrete-time snap-shots of the positions of all floaters \red{Citations needed!}.
With this in mind, we derive similar cost functions as above, that are based on such a finite data set \emph{only}. This will be the content of Section~\ref{sec:discr setting}. First, the \emph{dynamics} is discretized in time and space by conditioning a usual time-stepping method for the SDE~\eqref{eq:SDE} on the event that the random continuous trajectories are to be found in the set of known floater positions at the~$K\in\N$ given time instances. As above, we then derive two large-deviation semidistances~$\nu_K^\cross(x,y)$ and~$\nu_K^\meet(x,y)$ that have a clear probabilistic interpretation, that can be used to characterise coherent sets, and that are based on the finite data set rather than on the explicit velocity field. In fact, we will show that these discrete-space-time semidistances are really specific discretisations of the continuous-space-time semidistances~$\mu_T^\cross$ and~$\mu_T^\meet$. As shown in Section~\ref{subsec:discr end ldp}, they can be computed as shortest path lengths in a time-dependent weighted graph. We give an algorithm to compute these shortest paths in Appendix~\ref{app:algo}.

Let us stress that these semidistances are defined for deterministic dynamical systems. The random perturbation that is factored out by the large-deviation principle is merely acting as a catalyst to help quantify how strongly distinct trajectories mix---or, we should rather say how \emph{poorly}, as the transport from one trajectory to another is inversely proportional to their semidistance.

\paragraph{Coherence analysis with semidistances.}

In Section~\ref{sec:method} we describe how in general a semidistance on finite Lagrangian data can be used to analyse coherence. Key to our method is the notion of \emph{cornerstone}: a point that is furthest away from all other points. Cornerstones are though of as ``endpoints'' of a spanning structure, and ideally each cornerstone is in some sense the center of a coherent set. As a next step, trajectories can be clustered around cornerstones to yield coherent sets. Of course, this approach is very close to the \emph{k-means-} and \emph{fuzzy c-means clustering} of trajectories with respect to dynamical distances in~\cite{HaEtAl16,FrPa15}, with the important difference that the centers are not chosen by the heuristics of these clustering approaches, but with regard to the properties of coherent sets in the light of transport and mixing.

To exemplify the usefulness of the theory put forth in this paper, in sections~\ref{sec:method} and~\ref{sec:numerics} we test our approach on a number of standard test cases. Finally, Section~\ref{sec:outlook} discusses possible combinations of this work with other concepts.

\section{Large-deviation semidistances in continuous time and space}
\label{sec:cont setting}

In this section we study large deviations of the forms~\eqref{eq:cont ldp one traj} and \eqref{eq:cont ldp meet}. In large-deviation theory it is ofter easier to first study large deviations in a larger space. In our setting, we first study the large deviations of paths in Section~\ref{subsec:cont traj ldp} before transforming to the large deviations of the endpoints in Section~\ref{subsec:cont end ldp}. We end with a discussion of the resulting semidistances $\mu_T^\cross,\mu_T^\meet$ in Section~\ref{subsec:cont semidistances}.

\subsection{Large deviations of paths}
\label{subsec:cont traj ldp}

We denote paths by $w\suber{\cdot}$ to distinguish them from points $w$. Let $\prob$ be the Wiener measure, i.e., the probability that a Brownian path lies in a set $U\subset C(0,T;\R^d)$ is $\prob[\rv{w}_{(\cdot)}\in U]$. Recall that there does not exist a canonical probability measure on the space of paths, and so the Wiener measure can not be identified with a meaningful density. This means that one always needs to consider sets rather than particular realisations of the Brownian path. Nevertheless, large-deviation rates are always local, in the sense that they depend on one realisation only (the most likely one in the set $U$ under consideration). This motivates writing~$\rv{w}\suber{\cdot}\asymp f\suber{\cdot}$ if $\rv{w}\suber{\cdot}$ lies in an infinitesimal neighborhood~$U$ of the path~$f\suber{\cdot}$. %Asymptotics in~$\ep\to 0$ is expressed by the symbol~$\sim$ 
We will make this more precise below.

The large deviations for the SDE~\eqref{eq:SDE} are a standard result by Freidlin--Wentzell~\cite{FreidlinWentzell1998}. This result can be derived via a combination of \emph{Schilder's Theorem} and a \emph{Contraction Principle} as we now explain.

We first consider the noise part $\sqrt\ep \rv{w}_t$, which clearly converges (almost surely uniformly) to the constant path~$0$ as~$\ep\to0$. The corresponding large-deviation principle is given by Schilder's Theorem~\cite[Th.~5.2.3]{Dembo1998}:
\begin{equation}
  -\ep\log\prob\left[\sqrt{\ep}\rv{w}_{(\cdot)} \asymp w_{(\cdot)}\right] \xrightarrow[\ep\to0]{}
%  \left\{
%    \begin{array}{ll}
      \frac12\int_0^T\!\lvert \dot w_t\rvert^2 dt, %& w\suber{\cdot}\in H^1(0,T,\R^d).\\
%      \infty, & \text{otherwise.}
%    \end{array}
%  \right.
\label{eq:Schilders}
\end{equation}
for differentiable paths $w\suber{\cdot}$ starting from~$w_0=0$ (otherwise the limit will be~$\infty$).\footnote{More rigorously,~\eqref{eq:Schilders} means
$
\ep\log\prob\left[\sqrt{\ep}\rv{w}_{(\cdot)} \in U\right] \xrightarrow[\ep\to0]{}  -\inf_{w_{(\cdot)} \in U}\frac12\int_0^T\!\lvert \dot w_t\rvert^2 dt,
$
where, for technical reasons, this convergence is realized by a liminf lower bound for open sets~$U$ and limsup upper bound for closed sets~$U$; see~\cite{Dembo1998}.}

Let us assume that the velocity field $v(t,\cdot)$ is Lipschitz, so for each realisation of the Brownian path~$\rv{w}\suber{\cdot}=w\suber{\cdot}$ corresponds a unique solution~$\rv{x}\super{\ep}\suber{\cdot}$ of the SDE, starting from some given~$\rv{x}\super{\ep}_0=x$, see~\cite[Th.~5.2.1]{Oksendal2003}. The Contraction Principle~\cite[Th.~4.2.1]{Dembo1998} then states that the large-deviation rate of a path $x\suber{\cdot}$ is given by the minimum of~\eqref{eq:Schilders} over all realisations of the noise that give rise to that path, i.e.,
\begin{equation}
  -\ep\log\prob\left[\rv{x}\super{\ep}_{(\cdot)} \asymp x_{(\cdot)}\right]
  \xrightarrow[\ep\to0]{} \inf_{w\suber{\cdot} \,:\, \dot x_t=v(t,x_t) + \dot w_t} \,  \mfrac12\int_0^T\!\lvert \dot w_t\rvert^2\,dt =
%  & = \begin{cases}
        \mfrac12\int_0^T\!\lvert \dot x_t-v(t,x_t)\rvert^2\,dt, %&x\suber{\cdot}\in H^1(0,T,\R^d),\\
%        \infty,                                                 &\text{otherwise}.
%      \end{cases}
\label{eq:Freidlin-Wentzell}
\end{equation}
for differentiable paths~$x\suber{\cdot}$ starting from~$x_0=x$.

\subsection{Large deviations of endpoints}
\label{subsec:cont end ldp}

We now derive the large-deviation principle of the type~\eqref{eq:cont ldp one traj} as discussed in the introduction. In a sense, the pathwise large deviations~\eqref{eq:Freidlin-Wentzell} encode more information than is needed if we are only interested in the endpoint~$\rv{x}\super{\ep}_T\asymp \phi_{0,T}(y)$ of the random path. Another application of the Contraction Principle then states that the large-deviation rate for the endpoint is the minimum of~\eqref{eq:Freidlin-Wentzell} over all paths starting from~$x$ and ending in that given endpoint~$\phi_{0,T}\lbrack y\rbrack$, i.e.:
\begin{equation}
  -\ep\log\prob\big[ \rv{x}_T\super{\ep}\asymp \phi_{0,T}\lbrack y\rbrack \mid \rv{x}_0\super{\ep}=x\big] \xrightarrow[\ep\to0]{} \inf_{x\suber{\cdot}\,: x_0=x,x_T=\phi_{0,T}\lbrack y\rbrack}\, \mfrac12\int_0^T\!\lvert \dot x_t-v(t,x_t)\rvert^2\,dt =: \mu_T(x\smallto y).
\label{eq:cont ldp one traj explicit}
\end{equation}
This defines the `one-way' rate that we are after. 

The sum~\eqref{eq:cont fwbw distance} then defines the cross semidistance $\mu_T^\cross(x,y)$, and has a natural interpretation in terms of large deviations: As mentioned in the introduction, it arises from two independent and simultaneous copies $\rv{x}\super{\ep}_t,\rv{y}\super{\ep}_t$. By independence, the probability that $(\rv{x}\super{\ep}_T,\rv{y}\super{\ep}_T)\asymp(\phi_{0,T}\lbrack y\rbrack,\phi_{0,T}\lbrack x\rbrack)$ given $(\rv{x}\super{\ep}_0,\rv{y}\super{\ep}_0)=(x,y)$ is a product of one-way probabilities, yielding the sum of two one-way rates in the large deviations, see Figure~\ref{fig:cont symm dist}.
%For the other interpretation, let $\rv{x}\super{\ep}\suber{\cdot}$ be the random process on $(0,2T)$ that solves the SDE in $(0,T)$ with initial condition $\rv{x}\super{\ep}_0=x$, jumps from the endpoint $\rv{x}\super{\ep}_{T^-}$ to the pulled-back point $\rv{x}\super{\ep}_{T^+}:=\phi_{0,T}^{-1}\lbrack \rv{x}\super{\ep}_{T^-}\rbrack$ at time $T$, and then runs the SDE further on $(T,2T)$ with the new initial condition $\rv{x}\super{\ep}_{T^+}$. One can then study the probability that this process, starting in $x$, ends in $\phi_{0,T}\lbrack y\rbrack$ at time $T^-$ and in $\phi_{0,T}\lbrack x\rbrack$ at time $2T$. The large-deviation rate corresponding to this probability is exactly the semidistance $\mu_T^\cross(x,y)$.

%\label{rem:cont fwbw ldp}
%\end{remark}

A similar argument can be used to derive the meeting large deviations~\eqref{eq:cont ldp meet}. Let~$\rv{x}\super{\ep}_t$ and~$\rv{y}\super{\ep}_t$ be two independent solutions of the SDE~\eqref{eq:SDE}, starting from given~$x$ and~$y$, respectively. We consider the probability that both trajectories end in a given point, say~$\phi_{0,T}\lbrack z\rbrack$ for some~$z\in\R^d$, see Figure~\ref{fig:cont meet dist}. Assuming independence of the two trajectories, we immediately get
\begin{align*}
  &-\ep\log\prob\big[ \rv{x}_T\super{\ep}\asymp \phi_{0,T}\lbrack z\rbrack, \rv{y}_T\super{\ep} \asymp \phi_{0,T}\lbrack z\rbrack \mid \rv{x}_0\super{\ep}=x, \rv{y}_0\super{\ep}=y\big] \\
    &\qquad =
  -\ep\log\prob\big[ \rv{x}_T\super{\ep}\asymp \phi_{0,T}\lbrack z\rbrack \mid \rv{x}_0\super{\ep}=x \big] -\ep\log\prob\big[ \rv{y}_T\super{\ep}\asymp \phi_{0,T}\lbrack z\rbrack \mid \rv{y}_0\super{\ep}=y \big]\\
    &\qquad\!\xrightarrow[\ep\to0]{\eqref{eq:cont ldp one traj explicit}} \mu_T(x\smallto z) + \mu_T(y\smallto z).
\end{align*}
However, we are only interested in the probability that the two trajectories meet, and not in the point where they meet. A final Contraction Principle thus yields:
\begin{equation}
  -\ep\log\prob\big[ \rv{x}_T\super{\ep}\asymp \rv{y}_T\super{\ep} \mid \rv{x}_0\super{\ep}=x, \rv{y}_0\super{\ep}=y\big]
    \xrightarrow[\ep\to0]{}
  \inf_{z\in X} \mu_T\big(x\smallto z \big) + \mu_T\big(y\smallto z\big)
  =:  \mu_T^\meet(x,y).
\label{eq:cont ldp meet explicit}
\end{equation}
Observe that the two paths could also meet earlier and subsequently follow the same trajectory up until time $T$ with zero cost; the time $T$ thus acts as a maximum time at which the paths should meet.

\subsection{The semidistances}
\label{subsec:cont semidistances}

We now discuss some metric properties of the rate functionals. Recall from the introduction that we assumed that the flow is a diffeomorphism. Therefore $\mu_T(x\to y)=0$ if and only if $x=y$. It is then easy to see that, for any $x,y$,
\begin{enumerate}[(i)]
\item $\mu_T^{\cross}(x,y)\geq0$ \hspace{2.3cm}       and\qquad $\mu_T^{\meet}(x,y)\geq0$,
\item $\mu_T^{\cross}(x,y)=0 \iff x=y$ \hspace{0.25cm} and\qquad $\mu_T^{\meet}(x,y)=0 \iff x=y$,
\item $\mu_T^{\cross}(x,y)=\mu_T^{\cross}(y,x)$ \hspace{0.7cm} and\qquad $\mu_T^{\meet}(x,y)=\mu_T^{\meet}(y,x)$.
\end{enumerate}
However, the triangle inequality can fail, and so~$\mu_T^{\cross}$ and~$\mu_T^{\meet}$ are semidistances only.

We point out the following useful relation between the two. Observe that by the definition, $\mu_T^{\meet}(x,y)\leq \mu_T(x\smallto y) + \mu_T(y,y)=\mu_T(x\smallto y)$, and similarly $\mu_T^{\meet}(x,y)\leq \mu_T(y\smallto x)$. Therefore,
\begin{equation*}
  \mu_T^{\meet}(x,y) \leq \min\big\{ \mu_T(x\smallto y) , \mu_T(y\smallto x)\big\} \leq \max\big\{ \mu_T(x\smallto y) , \mu_T(y\smallto x)\big\} \leq \mu_T^{\cross}(x,y).
\end{equation*}
In order to investigate which semidistance is more suitable to study coherence, one would need to study in which setting the gap $\mu_T^{\cross}(x,y)-\mu_T^{\meet}(x,y)$ becomes large. This is beyond the scope of this paper, but we will show for several examples that both work as they should.

%The inequality does suggest a third candidate: the maximum. Since the maximum of two large-deviation rates does not have a clear large-deviation interpretation we do not study this candidate.

\begin{remark}[Invariance under time-reversal]
Note the following invariance property of $\mu_T$ under time-reversal:
\begin{equation*}
  \mu_T(x\smallto y)= \inf_{ \substack{ y\suber{\cdot}\,: y_0=\phi_{0,T}\lbrack y\rbrack\\ y_T=\phi_{0,T}^{-1}\lbrack \phi_{0,T}\lbrack x\rbrack\rbrack } }\, \mfrac12\int_0^T\!\lvert \dot y_t+v(T-t,y_t)\rvert^2\,dt =:{\overleftarrow\mu}\!_T\big( \phi_{0,T}\lbrack y\rbrack \smallto \phi_{0,T}\lbrack x \rbrack \big),
\end{equation*}
where $\overleftarrow\mu\!_T$ is the one-way rate associated to the backward system $\dot y_t = - v(T-t, y_t)$. This time-reversal property is retained for the cross semidistance: $\mu_T^\cross(x,y)=\overleftarrow{\mu}\!_T^\cross\big( \phi_{0,T}\lbrack x\rbrack,\phi_{0,T}\lbrack y\rbrack \big)$. However, the meeting semidistance $\mu_T^\meet(x,y)=\inf_z \overleftarrow\mu\!_T\big(z\smallto\phi_{0,T}\lbrack x\rbrack\big)+\overleftarrow\mu\!_T\big(z\smallto\phi_{0,T}\lbrack y\rbrack\big)$ is the same as the cost for the backward trajectories to \emph{start} in a joint position and end in $\phi_{0,T}\lbrack x\rbrack$ and $\phi_{0,T}\lbrack y\rbrack$.
\label{rem:cont time-reversal}
\end{remark}

\subsection{A simple example}
\label{ssec:simpexam1}

Let us consider a very simple example, where the domain of interest is the interval~$[0,L]$, and there is no dynamics, i.e.,~$v\equiv 0$. Its primary purpose is to form our intuition and expectations about how the semidistances work in more complicated settings. In particular, we shall see how the semidistances scale in time and system size. 

The system is considered on the time interval~$[0,T]$. One can then easily see that~$\dot{x}_t \equiv L/T$ is an optimal path in~\eqref{eq:cont ldp one traj explicit}, thus giving
\[
\mu_T(0\smallto L) = \mu_T(L\smallto 0) = \mfrac12\int_0^T\!\left(\mfrac{L}{T}\right)^2\,dt = \frac{L^2}{2T}\,,
\]
and so $\mu^\cross_T(0,L)= \frac{L^2}{T}$. Thus, also~ $\mu_T(0\smallto L/2) = \mu_T(L\smallto L/2) = \frac{L^2}{8T}$. In general, the one-way cost is proportional to the squared distance and inversely proportional to time. This also gives
\[
\mu_T^{\meet}(0,L) = \frac{L^2}{4T}\,,
\]
so, in this symmetric situation the meeting distance is half of one-way cost and quarter of the cross semidistance.

We will revisit this example in the next section, and will realize that the behavior of the discrete semidistances deviates from the one observed here for continuous space and time.

\section{Large-deviation semidistances in discrete time and space}
\label{sec:discr setting}

As mentioned in the introduction, the cost functions $\mu_T^\cross$ and $\mu_T^\meet$ are difficult to calculate explicitly, and impossible if the velocity or flow field is not explicitly known. In this section we take a more practical approach. We will assume that the only information at hand is the position at finite times of a finite number $I$ of floaters.

To be more specific, let~$\{x_k\super{i}\}_{k=0,\hdots; K, i=1,\hdots, I} \subset \R^d$ be given positions of floaters~$i=1,\hdots,I$ at time~$k\tau$ for~$k=0,\hdots,K$ for some~$\tau>0$. Assuming that the floaters sample from the deterministic flow field $\phi_{s,t}$, we know that for each floater $i$,
\begin{equation}
  x\super{i}_{k+1} = \phi_{k\tau,(k+1)\tau}\lbrack x\super{i}_k\rbrack\,.
\label{eq:Iflow}
\end{equation}
If we would add noise to the system, we would find random particles described by the set of SDEs
\begin{align}
  d\rv{x}_t\super{i,\ep} = v(t,\rv{x}_t\super{i,\ep})dt + \sqrt\ep d\rv{w}\super{i}_t, && \rv{x}_0\super{i,\ep}=x_0\super{i}, &&\text{for } i=1,\hdots, I,
\label{eq:SDEs}
\end{align}
where~$\rv{w}\super{i}$ are now independent standard Brownian motions.

Our strategy is to study the probability that random particles described by the SDEs~\eqref{eq:SDEs} deviate from the given floater trajectories~\eqref{eq:Iflow}, conditional to the fact that all our knowledge about the otherwise unknown flow field~$\phi_{s,t}$ comes from the time- and space-discrete set of trajectory data~\eqref{eq:Iflow}. We first approximate the SDEs~\eqref{eq:SDEs} by discrete-time, continuous-space Markov processes in Section~\ref{subsec:discrete time}, as it is done in standard time-stepping methods for SDEs~\cite{KlPl10}. Next, in Section~\ref{subsec:discrete time and space} we condition these discrete-time processes on the given floater positions. Then we calculate the large-deviation rate for trajectories in Section~\ref{subsec:discr traj ldp}, and for endpoints in Section~\ref{subsec:discr end ldp}. Finally, we end the section with a discussion of the metric properties of the resulting large-deviation rates in Section~\ref{subsec:discr semidistances}.

\subsection{Discrete-time approximation}
\label{subsec:discrete time}

We first focus our attention to one time step $k\tau\to(k+1)\tau$ of one trajectory $i$, and temporarily drop the superindex for brevity. Since the noise process is a standard Brownian motion, we know its density,
\begin{equation}
  \frac{d\prob[\sqrt{\ep} \rv{w}_{(k+1)\tau}\in dy\,\big\vert\, \sqrt{\ep} \rv{w}_{k\tau}=x]}{dy} = \left(2\pi\ep\tau\right)^{-d/2} \exp\left(-\frac{|x-y|^2}{2\ep\tau}\right)\,.
\label{eq:Inoise}
\end{equation}
Hence, we have exact information on the purely deterministic part of the SDE by~\eqref{eq:Iflow}, and on the purely noise part by~\eqref{eq:Inoise}. We combine this information by using the following time-stepping approximation for the SDE~\eqref{eq:SDE}.

Fix an $\alpha\in\lbrack0,1\rbrack$, and let $(\rv{\xi}_k^+,\rv{\xi}_k^-)_{k=0,\hdots,K}$ be independent normally distributed $\R^d$-valued random variables with unit variance. Given the approximated random position $\tilde{\rv{x}}_k$ at time $k\tau$, we iterate
\begin{equation}
\begin{aligned}
  \tilde{\rv{x}}^+_{k} &:= \tilde{\rv{x}}_{k} + \sqrt{\alpha\tau\ep}\, \rv{\xi}^+_k \\
  \tilde{\rv{x}}^-_{k+1} &:= \phi_{k\tau,(k+1)\tau}\lbrack\tilde{\rv{x}}^+_k\rbrack \\
  \tilde{\rv{x}}_{k+1} &:= \tilde{\rv{x}}^-_{(k+1)\tau} + \sqrt{(1-\alpha)\tau\ep}\, \rv{\xi}^-_{k+1}
\end{aligned}
\label{eq:time-st}
\end{equation}
Here,~$\tilde{\rv{x}}^+_k$ and~$\tilde{\rv{x}}^-_{k+1}$ are only auxiliary (intermediate) steps. The method~\eqref{eq:time-st} is a special case of a \emph{splitting method}, since the deterministic evolution and purely noise parts of the SDE~\eqref{eq:SDE} are handled separately in the distinct steps; here it would be ``noise-flow-noise''.

We would like to stress that our choice of discretization is made on the basis that we can use the available information on the flow given by~\eqref{eq:Iflow}. In the realm of one-step methods for SDEs we are bound to choices of the form~\eqref{eq:time-st}, because there is no information on the drift other than the flow generated by it on prescribed time intervals~$[k\tau, (k+1)\tau)$. Given the form of time-stepping \eqref{eq:time-st}, the optimal (in the sense of highest weak consistency order~\cite{KlPl10}) approximation of the SDE is obtained by choosing~$\alpha=1/2$. That is the so-called \emph{Strang-splitting}~\cite{Str68}, and has weak order two, while for~$\alpha\neq 1/2$ we only get order one.\footnote{Formally, this can be seen by denoting the generators of the noise process and advection by~$A$ and~$B$, respectively, and estimating the difference of the Markov propagators associated with~\eqref{eq:SDEs} and~\eqref{eq:time-st} by performing formal Taylor expansions in~$\tau$ with the non-commuting operators~$A$ and~$B$, to obtain
\[
e^{\tau (A+B)} - e^{(1-\alpha)\tau A} e^{\tau B} e^{\alpha\tau A} = 
\left\{ \begin{array}{ll}
\mathcal{O}(\tau^3), & \alpha = 1/2,\\
\mathcal{O}(\tau^2), & \alpha \neq 1/2.
\end{array}
\right.
\]
}
%In this sense, with~$\alpha=1/2$ the time-stepping~\eqref{eq:time-st} (and its restriction to the finite floater positions,~\eqref{eq:onejump}) can be seen as an optimal approximation to the SDEs~\eqref{eq:SDEs}.

Having performed discretization in time, in the next section we derive a discrete-time, discrete-space, $\alpha$-dependent Markov chain that we will use to derive discrete semidistances.

\subsection{Conditioning on finite data}
\label{subsec:discrete time and space}

Recall that we considered one discrete-time process $(\tilde{\rv{x}}_k)_{k=1,\hdots,K}$ with initial condition $\tilde{\rv{x}}_0=x_0^{(i)}$, and that we suppressed the dependency on $i$. For each $k=0,\hdots,K$, we introduce the set
\begin{equation*}
  \A_k:=\{ x\super{j}_k\}_{j=1,\hdots I}.
\end{equation*}
of available points at time $t=k\tau$. We now condition the random process on the event that for each realisation $\tilde{\rv{x}}_k\in\A_k$ \emph{and} for each intermediate point $\tilde{\rv{x}}_k^+\in\A_k$. This automatically implies conditioning of the other intermediate points $\tilde{\rv{x}}_{k+1}^-\in\A_{k+1}$ due to~\eqref{eq:Iflow}. We choose to condition on the intermediate points for practical reasons; otherwise we would not be able to perform the second step in~\eqref{eq:time-st}, since the discrete trajectories are our only information about the flow, cf.~Remark~\ref{rem:partial conditioning} below.

\begin{figure}[h]
\centering
\begin{tikzpicture}[font=\footnotesize, scale=1.4]
\begin{scope}[decoration={markings,mark=at position 0.5 with {\arrow{>}}}]
  
  \draw[gray!30, fill = gray!30] (1.3,-0.2) rectangle (1.7,2.2);
  \node[gray] at (1.5,-0.5) {$\A_k$};
  
  \filldraw(1.5,0)   circle (0.05);
  \filldraw(1.5,0.5) circle (0.05) node[anchor=east]{$\rv{i}_k=j$};
  \filldraw(1.5,1)   circle (0.05);
  \filldraw(1.5,1.5) circle (0.05);
  \filldraw(1.5,2)   circle (0.05);

  \draw[gray!30, fill = gray!30] (2.8,-0.2) rectangle (3.2,2.2);
  \node[gray] at (3,-0.5) {$\A_{k+1}$};
  
  \filldraw(3,0.1)   circle (0.05);
  \filldraw(3,0.5)   circle (0.05);
  \filldraw(3,0.9)   circle (0.05);
  \filldraw(3,1.6)   circle (0.05) node[anchor=west]{$\rv{i}_{k+1}=m$};
  \filldraw(3,1.9)   circle (0.05);

  \draw[dashed](1.5,0)--(3,0.1);
  \draw[dashed](1.5,0.5)--(3,0.5);
  \draw[dashed](1.5,1.5)--(3,1.6);
  \draw[dashed](1.5,2)--(3,1.9);

  \draw[-,postaction={decorate}](1.5,0.5).. controls (1.3,0.6) and (1.3,0.9).. node[near end,anchor=east]{$p_k^+(j,\ell)$}(1.5,1);
  \draw[-,postaction={decorate}](1.5,1)--node[midway,anchor=south]{$\phi_{k\tau,(k+1)\tau}$}(3,0.9);
  \draw[-,postaction={decorate}](3,0.9).. controls (3.2,1.1) and (3.2,1.4) .. node[midway,anchor=west]{$p_{k+1}^-(\ell,m)$}(3,1.6);

\end{scope}
\end{tikzpicture}
\caption{One time step of the discrete-time discrete space Markov chain~$\rv{i}_k$.}
\label{fig:time-stepping}
\end{figure}
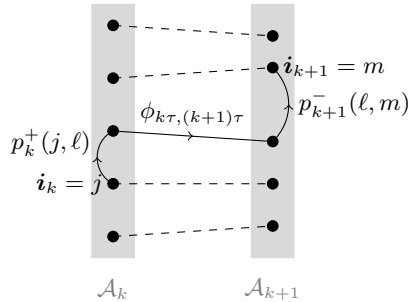

The conditioning on the finite data set results in replacing the discrete-time continuous-space process by a fully discrete-time discrete-space Markov chain that hops between the given trajectories. Therefore, the state of the new Markov chain can be represented by the labels $j=1,\hdots,I$; this is particularly useful since the deterministic flow~\eqref{eq:Iflow} will change the positions but not the labels. Since the resulting process is still Markovian, we can fully characterise its behaviour through its transition probabilities for one time-step $k\to k+1$. We now calculate these transition probabilities, dealing with each step in \eqref{eq:time-st} separately. See Figure~\ref{fig:time-stepping} for a sketch.

For the transition from $\tilde{\rv{x}}_k$ to $\tilde{\rv{x}}_k^+$, where we know the increment distribution~\eqref{eq:Inoise}, note that we are in fact conditioning on a null set, so that the conditional probabilities are sensibly defined as the limits over balls~$B_r(\cdot)$ of small radii~$r\to0$ around these points. We thus obtain, for any~$j,\ell=1,\hdots I$:
\begin{align}
  p_{k}^+(j,\ell) &:= \prob\left[ \tilde{\rv{x}}^+_k = x\super{\ell}_k\,\big\vert\, \tilde{\rv{x}}_k = x\super{j}_k\ \text{and}\ \tilde{\rv{x}}^+_k\in \A_k \right] \notag\\
    & = \lim_{r\to0} \frac{\prob\left[ \tilde{\rv{x}}^+_k\in B_r(x\super{\ell}_k)\, \big\vert\, \tilde{\rv{x}}_k = x\super{j}_k\right]}
                          {\prob\left[\tilde{\rv{x}}^+_k\in B_r(\A_k)  \big\vert\, \tilde{\rv{x}}_k = x\super{j}_k\right]} \notag\\
    & = \frac{\exp\left(-\left|x\super{\ell}_k-x\super{j}_k\right|^2/(2\alpha\ep\tau)\right)}{\sum_{\hat\ell=1}^I \exp\left(-\left|x\super{\hat\ell}_k-x\super{j}_k\right|^2/(2\alpha\ep\tau)\right)}\,,
\label{eq:pplus conditioning}
\end{align}
and similarly for the transition from $\tilde{\rv{x}}_{k+1}^-$ to $\tilde{\rv{x}}_{k+1}$:
\begin{align}
  p_{k+1}^-(\ell,m) &:= \prob\left[ \tilde{\rv{x}}_{k+1} = x\super{m}_{k+1}\,\big\vert\, \tilde{\rv{x}}^-_{k+1} = x\super{\ell}_{k+1}\ \text{and}\ \tilde{\rv{x}}_{k+1}\in \A_{k+1} \right]\,\notag\\
    & = \frac{\exp\left(-\left|x\super{m}_{k+1}-x\super{\ell}_{k+1}\right|^2/\left(2(1-\alpha)\ep\tau\right)\right)}{\sum_{\hat m=1}^I \exp\left(-\left|x\super{\hat m}_{k+1}-x\super{\ell}_{k+1}\right|^2/\left(2(1-\alpha)\ep\tau\right)\right)}\,.
\label{eq:pminus conditioning}
\end{align}
Since the transition from~$\tilde{\rv{x}}^+_k$ to~$\tilde{\rv{x}}^-_{k+1}$ is deterministic (middle equation in~\eqref{eq:time-st}), we have that,
\begin{equation}
  P_k(j,m):=\prob\left[\tilde{\rv{x}}_{k+1} = x\super{m}_{k+1}\,\big\vert\, \tilde{\rv{x}}_k=x\super{j}_k\text{ and } \tilde{\rv{x}}^+_{k}\in \A_k,  \tilde{\rv{x}}^-_{k+1} \in \A_{k+1} \right] = \sum_{\ell=1}^I p_k^+(j,\ell)p_{k+1}^-(\ell,m)\,.
\label{eq:onejump}
\end{equation}
In words, the process performs the following three subsequent steps for \emph{one} time step (see Figure~\ref{fig:time-stepping}):
\begin{enumerate}
\item Start in~$x\super{j}_k$, and perform a jump to some~$x\super{\ell}_k$ with probability~$p\super{k,+}_{j,\ell}$,
\item Perform a deterministic jump from~$x\super{\ell}_k$ to~$x\super{\ell}_{k+1}$,
\item Perform a jump from~$x\super{\ell}_{k+1}$ to~$x\super{m}_{k+1}$ with probability~$p\super{k+1,-}_{\ell,m}$.
\end{enumerate}

The transition probabilities $P_k(j,m)$ define our new, discrete-time Markov chain $(\rv{i}_k)_{k=0,\hdots,K}$ on the discrete space $\{1,\hdots,I\}$. To shorten notation, we will write $i\suber{\cdot}:=(i_k)_{k=0,\hdots,K}$ for a discrete path, analogous to the continuous-time setting. By the Markov property, the probability that the Markov chain realizes such a path is simply
\begin{equation}
  \prob\big\lbrack \rv{i}\suber{\cdot} = i\suber{\cdot} \big\rbrack = \prod_{k=0}^{K-1} P_k(i_k,i_{k+1}),
\label{eq:discrete trajectory prob}
\end{equation}
where we assumed that the chain starts (deterministically) from~$i_0$.

\subsection{Large deviations of discrete trajectories}
\label{subsec:discr traj ldp}

We now study the large deviations of the discrete Markov chain $\rv{i}_k$. Similarly to the continuous setting from Section~\ref{sec:cont setting} we start from the large deviations of paths. First we calculate the large deviations for $p_k^+(j,\ell)$ and $p_{k+1}^-(\ell,m)$. By the Laplace principle~\eqref{eq:Laplace}, 
\begin{align*}
  -\ep\log p_k^+(j,\ell) %&= -\ep\log\frac{\exp\left(-\left\lvert x\super{\ell}_k-x\super{j}_k\right\rvert^2/(2\alpha\ep\tau)\right)}{\sum_{\hat\ell=1}^I \exp\left(-\left\lvert x\super{\hat\ell}_k-x\super{j}_k\right\rvert^2/(2\alpha\ep\tau)\right)} \\
    &\stackrel{\eqref{eq:pplus conditioning}}{=}\ep\log \sum_{\hat\ell=1}^I \exp\left(\mfrac{\left\lvert x\super{\ell}_k-x\super{j}_k\right\rvert^2 - \left\lvert x\super{\hat\ell}_k-x\super{j}_k\right\rvert^2}{2\alpha\ep\tau}\right) \\
    &\xrightarrow[\ep\to0]{} \max_{\hat\ell=1,\hdots,I} \mfrac{\left\lvert x\super{\ell}_k-x\super{j}_k\right\rvert^2 - \left\lvert x\super{\hat\ell}_k-x\super{j}_k\right\rvert^2}{2\alpha\tau} 
    =\mfrac{\lvert x\super{\ell}_k-x\super{j}_k\rvert^2}{2\alpha\tau}.
\end{align*}
We will make this simplification again below. Similarly, we obtain
\begin{equation*}
-\ep\log p_{k+1}^-(\ell,m) 
     \xrightarrow[\ep\to0]{\eqref{eq:pminus conditioning}}\mfrac{\lvert x\super{m}_{k+1}-x\super{\ell}_{k+1}\rvert^2}{2(1-\alpha)\tau}.
\end{equation*}
Using these two exponential approximations, we can again use the Laplace principle~\eqref{eq:Laplace} to find for the jump probability of one time step:
\begin{align*}
    \lim_{\ep\to 0}-\ep\log P_k(j,m) &\stackrel{\eqref{eq:onejump}}{=}   \lim_{\ep\to 0}-\ep\log \sum_{\ell=1}^I p_k^+(j,\ell)p_{k+1}^-(\ell,m) \\
    &= \min_{\ell=1,\ldots, I} \lim_{\ep\to 0} \left(-\ep\log p_k^+(j,\ell) -\ep\log p_{k+1}^-(\ell,m) \right) \\
    &= \min_{\ell=1,\hdots,m}    
      \mfrac{\lvert x\super{\ell}_k-x\super{j}_k\rvert^2}{2\alpha\tau}
      +
      \mfrac{\lvert x\super{m}_{k+1}-x\super{\ell}_{k+1}\rvert^2}{2(1-\alpha)\tau}.
\end{align*}
Finally the large-deviation rate of a discrete path is
\begin{align}
  -\ep\log \prob\left\lbrack \rv{i}\suber{\cdot} = i\suber{\cdot} \right\rbrack &\stackrel{\eqref{eq:discrete trajectory prob}}{=} -\ep\log \prod_{k=0}^{K-1} P_k(i_k,i_{k+1}) \notag\\
    &\xrightarrow[\ep\to0]{} \sum_{k=0}^{K-1} \min_{\ell=1,\hdots,I}
      \mfrac{\lvert x\super{\ell}_k-x\super{i_k}_k\rvert^2}{2\alpha\tau}
      +
      \mfrac{\lvert x\super{i_{k+1}}_{k+1}-x\super{\ell}_{k+1}\rvert^2}{2(1-\alpha)\tau}
    := \J( i\suber{\cdot}),
\label{eq:discrete trajectory ldp}
\end{align}

\begin{remark}
\label{rem:partial conditioning}

Recall that we conditioned on the event that all $\tilde{\rv{x}}_k$ as well as the intermediate points $\tilde{\rv{x}}^+_{k}$ lie in the set $\A_k$ of available points. One might argue that in practice only the points $\tilde{\rv{x}}_{k}$ are measured to lie in $\A_k$, while the other two are mathematical constructs that may lie anywhere. However, if we would relax this conditioning and follow the calculations as above, we would find:
\begin{align*}
  -\ep\log P_k(j,m) \xrightarrow[\ep\to0]{}
   \min_{x\in\R^d}\Big\{ \mfrac{\lvert x-x\super{j}_k\rvert^2}{2\alpha\tau} + \mfrac{\left|x\super{m}_{k+1}-\phi_{t_k,t_{k+1}}\lbrack x\rbrack \right|^2}{2(1-\alpha)\tau} - \min_{\hat m=1,\hdots,I} \mfrac{\left|x\super{\hat m}_{k+1}-\phi_{t_k,t_{k+1}}\lbrack x\rbrack \right|^2}{2(1-\alpha)\tau} \Big\}.
\end{align*}
Since this large-deviation rate still depends on the unknown flow field~$\phi$, it can not be used if only the data of a finite number of floaters is available.
\end{remark}

\begin{remark}[Missing data and non-uniform time-sampling]
\label{rem:missing data}

Note that the construction works exactly as described above even if information about trajectories is partially missing. The conditioning on the set~$\mathcal{A}_k$ works identically, but now these sets might have different cardinalities smaller or equal~$I$. Observe that our only information about the deterministic flow for times in~$[k\tau,(k+1)\tau)$ comes from those trajectories that are available both in~$\mathcal{A}_k$ and~$\mathcal{A}_{k+1}$. If this intersection is empty, we need to skip that time slice completely. This is not a problem, since our choice of sampling time uniformly by the step size~$\tau$ was solely in order to ease presentation. As the reader has probably observed, the extension for varying time steps~$\tau_k$ is straightforward.
\end{remark}

\subsection{Large deviations of endpoints}
\label{subsec:discr end ldp}

Analogously to the continuous setting, we study the large deviations of the one-way probability to hop from $i$ to $j$ in discrete time $K$, and the meeting probability that two independent chains, starting from $i$ and $j$ respectively, meet by discrete time~$K$ or earlier. Since the paths~\eqref{eq:discrete trajectory ldp} encode more information than the endpoints, we can now easily derive the large deviations of the one-way probability by a Contraction Principle. Indeed, for any two indices $i,j=1,\hdots,I$,
\begin{equation}\label{eq:discr rate}
  -\ep\log\prob\lbrack \rv{i}_K=j \mid \rv{i}_0=i\rbrack \xrightarrow[\ep\to0]{} \min_{i\suber{\cdot}:\,i_0=i,i_K=j} \J(i\suber{\cdot}) =: \nu_K(i\smallto j),
\end{equation}
where $\J$ is the discrete-path large-deviation rate~\eqref{eq:discrete trajectory ldp}. Note that~$\J$ is the shortest path length in a graph with time-dependent edge weights
\[
w_k(i,j) = \min_{\ell=1,\hdots,I}
      \mfrac{\lvert x\super{\ell}_k-x\super{i}_k\rvert^2}{2\alpha\tau}
      +
      \mfrac{\lvert x\super{j}_{k+1}-x\super{\ell}_{k+1}\rvert^2}{2(1-\alpha)\tau}\,.
\]

Again, the sum $\nu_K^\cross(i,j):=\nu_K(i\smallto j)+\nu_K(j\smallto i)$ can be given an interpretation in terms of large deviations as in Section~\ref{subsec:cont end ldp}. Moreover, following the same argument as in~\eqref{eq:cont ldp meet explicit}, if we take two independent trajectories $\rv{i}\suber{\cdot}$ and $\rv{j}\suber{\cdot}$, then
\begin{equation*}
  -\ep\log\prob\lbrack \rv{i}_K=\rv{j}_K \mid \rv{i}_0=i, \rv{j}_0=j \rbrack \xrightarrow[\ep\to0]{} \min_{\ell=1,\hdots,I} \nu_K(i\smallto\ell) + \nu_K(j\smallto\ell) =: \nu_K^\meet(i,j).
\end{equation*}

\subsection{The semidistances}
\label{subsec:discr semidistances}

It is easily checked that in the discrete setting the properties of a semidistance are also satisfied:
\begin{enumerate}[(i)]
\item $\nu_K^{\cross}(i,j)\geq0$ \hspace{2.5cm}       and\qquad $\nu_K^{\meet}(i,j)\geq0$,
\item $\nu_K^{\cross}(i,j)=0 \iff i=j$ \hspace{0.25cm} and\qquad $\nu_K^{\meet}(i,j)=0 \iff i=j$,
\item $\nu_K^{\cross}(i,j)=\nu_K^{\cross}(j,i)$ \hspace{0.7cm} and\qquad $\nu_K^{\meet}(i,j)=\nu_K^{\meet}(j,i)$.
\end{enumerate}
Furthermore, the triangle inequality fails, but we again have the following estimate:
\begin{equation*}
  \nu_K^{\meet}(i,j) \leq \min\big\{ \nu_K(i\smallto j) , \nu_K(j\smallto i)\big\} \leq \max\big\{ \nu_K(i\smallto j) , \nu_K(j\smallto i)\big\} \leq \nu_K^{\cross}(i,j).
\end{equation*}

Both semidistances can be computed from shortest-path costs, where the cost of a path is given by~\eqref{eq:discrete trajectory ldp}. We stress that this expression is fairly simple, and depends on the flow field through the \emph{known} positions of the floaters $x_k\super{\ell}$ only. Because of this: 1) these costs can be used in practice if the velocity field is unknown (Section~\ref{sec:method} and Section~\ref{sec:numerics}); 2) these costs can even be applied to cases where there may not be an underlying velocity field, as for example in discrete-time dynamical system (Section~\ref{subsec:two mixing subdomains}).

These semidistances can be computed by first computing the one-way rates~$\nu_K(i\to j)$ using Algorithm~\ref{algo1}, see Appendix~\ref{app:algo}. From these rates one readily obtains the semidistances via
\begin{align*}
  \nu_K^\cross(i,j)=\nu_K(i\to j) + \nu_K(j\to i) &&\text{and}&& \nu_K^\meet(i,j)=\min_{\ell=1,\hdots,I} \nu_K(i\to \ell) + \nu_K(j\to \ell).
\end{align*}

\begin{remark}[Time-reversal for discrete semidistances] Similarly to Remark~\ref{rem:cont time-reversal}, the one-way cost satisfies the time-reversal property $\nu_K(i\to j) = {\overleftarrow \nu}\!\!_K(j\to i)$, provided $\alpha=1/2$, where $\overleftarrow\nu\!\!_K$ is the cost associated to the backward dynamics. Moreover, this time-reversal property also holds for the cross semidistance, whereas for the meeting semidistance $\nu_K^\meet(i,j)=\min_{\ell=1,\hdots,I} \overleftarrow\nu\!\!_K(\ell\to i) + \overleftarrow\nu\!\!_K(\ell \to j)$. Apart from the superior consistency order discussed in Section~\ref{subsec:discrete time}, the invariance of semidistances under time reversal is another reason for choosing~$\alpha=1/2$.
\label{rem:discr time-reversal}
\end{remark}

\begin{remark} Other large-deviation-based semidistances are also possible. If one considers the ``noise-flow'' (i.e.,~$\alpha=1$) time-stepping scheme for the SDE rather than ``noise-flow-noise'', expression~\eqref{eq:discrete trajectory ldp} simplifies a bit. As another example of a large-deviation-based semidistance between two given discrete paths $\{x_k\super{i},x_k\super{j}\}_{k=0,\dots K}$, one could consider the probability to hop back and forth between the two trajectories, see Figure~\ref{fig:L2 distance}. In that case we find in the large-deviation scaling for~$\alpha=1$:
\begin{equation}
  -\ep\log\prob\lbrack \rv{i}_1=j_1,\rv{i}_2=i_2,\dots \mid \rv{i}_0=i_0\rbrack \xrightarrow[\ep\to0]{} \sum_{k=0}^{K-1} \frac{\lvert x_k\super{i}-x_k\super{j}\rvert^2}{2\tau},
\label{eq:ldp L2}
\end{equation}
for~$\alpha=0$ the sum would go from~$k=1$ to~$K$. Naturally, this is simply the $L^2$-distance between two trajectories, as considered earlier in~\cite{FrPa15}. Although this construction is very easy to calculate and its square root is a genuine metric, it is less interpretable as a cost for transport and mixing.

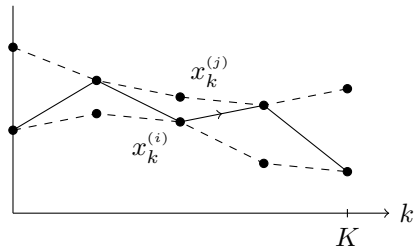
\begin{figure}[h]
\centering
\begin{tikzpicture}[scale=1.1]
\begin{scope}[decoration={markings,mark=at position 0.5 with {\arrow{>}}}] 
  \tikzstyle{every node}=[font=\small]
  \draw[->](0,0)--(4.5,0) node[anchor=west]{$k$};
  \draw(4,0.05)--(4,-0.05) node[anchor=north]{$K$};
  \draw(0,0)--(0,2.5);
  \draw[dashed](0,1)--(1,1.2)--(2,1.1)--(3,0.6)--(4,0.5);
  \filldraw(0,1) circle (0.05);
  \filldraw(1,1.2) circle (0.05);
  \filldraw(2,1.1) circle (0.05) node[anchor=north east]{$x_k\super{i}$};
  \filldraw(3,0.6) circle (0.05);
  \filldraw(4,0.5) circle (0.05);

  \draw[dashed](0,2)--(1,1.6)--(2,1.4)--(3,1.3)--(4,1.5);
  \filldraw(0,2) circle (0.05);
  \filldraw(1,1.6) circle (0.05);
  \filldraw(2,1.4) circle (0.05) node[anchor=south west]{$x_k\super{j}$};
  \filldraw(3,1.3) circle (0.05);
  \filldraw(4,1.5) circle (0.05);

  \draw(0,1)--(1,1.6)--(2,1.1);
  \draw(2,1.1)[postaction={decorate}] --(3,1.3);
  \draw(3,1.3)--(4,0.5);
\end{scope}
\end{tikzpicture}
\caption{Hopping back and forth (solid line) between two given trajectories (dashed lines).}
\label{fig:L2 distance}
\end{figure}

\label{rem:other semidistances}
\end{remark}

\begin{remark} It should be noted that the semidistances $\nu_K^\meet,\nu_K^\cross$ scale quadratically in space; this becomes even more apparent in the example considered in Section~\ref{ssec:simpexam2}. In the case of the $L^2$-distance~\eqref{eq:ldp L2}, the cost becomes a genuine distance after taking the square root. However, if we take the square roots of $\nu_K^\meet$ and $\nu_K^\cross$, the triangle inequality still fails. We therefore stick to the quadratic scaling as this has the most direct interpretation as large-deviation costs.
\end{remark}

\begin{remark}[Eulerian transport vs Lagrangian mixing] \label{rem:Euler vs Lagrange}
When speaking of transport in this paper, we mean ``\emph{transport} (of probability) \emph{from a trajectory to another}'', to express how the dynamics is mixing up regions these two trajectories come in contact with. This can be seen as a Lagrangian perspective. We express with large-deviation rates the \emph{unlikeliness} of transitions between trajectories, and these are then computed as shortest paths, cf.\ Section~\ref{subsec:discr end ldp}.
Deceivingly similar mathematical constructions show up in~\cite{SGetal15PRE}, where the authors consider ``highly probable paths'' of non-homogeneous Markov chains, which also leads to a time-dependent shortest path problem. Note, however, that this is orthogonal to our concept, as this is quantifying \emph{likeliness}.
A further important distinction is, that their Markov chain is constructed in an Eulerian manner (opposed to our Lagrangian setting), meaning that it describes transport between fixed regions of state space; serving as a discretization of the flow field~\cite{FPET07,FrSaMo10}.
\end{remark}

\subsection{Discretization of the continuous semidistances}
\label{ssec:ldp then discretize}

We now show that the one-way discrete space-time cost $\nu_K$ can also be obtained by discretizing the continuous space-time cost~$\mu_T$. This means that discretization and derivation of the large-deviation principle are interchangeable operations (if done the right way). We will not be precise about the discretization error; of course one needs to assume that the number of floaters is sufficiently large.

We first divide the time interval into subintervals $\lbrack 0,T)=\bigcup_{k=0}^{K-1} \lbrack k\tau , (k+\alpha)\tau)\cup \lbrack(k+\alpha)\tau,(k+1)\tau)$. 
Recall that $\phi_{t_0,t}$ is the flow associated to $v(t,\cdot)$, that is, for any~$t_0,x$,
\begin{equation*}
  \partial_t \phi_{t_0,t}\lbrack x\rbrack = v\big(t,\phi_{t_0,t}\lbrack x\rbrack\big).
\end{equation*}
Note in what follows that~$x_{(\cdot)}$ is some path, not necessarily a trajectory of the flow. In each interval $\lbrack k\tau , (k+\alpha)\tau)$ we approximate by finite differences:
\begin{align*}
  \dot x_t \approx \frac{x_{(k+\alpha)\tau}-x_{k\tau}}{\alpha\tau}
    &&\text{and}&&
  v(t,x_t) \approx \frac{x_{(k+\alpha)\tau} - \phi_{(k+\alpha)\tau,k\tau}\lbrack x_{(k+\alpha)\tau}\rbrack}{\alpha\tau}.
\end{align*}
In each interval $\lbrack (k+\alpha)\tau,(k+1)\tau)$ we approximate:
\begin{align*}
  \dot x_t \approx \frac{x_{(k+1)\tau}-x_{(k+\alpha)\tau}}{(1-\alpha)\tau}
    &&\text{and}&&
  v(t,x_t) \approx \frac{\phi_{(k+\alpha)\tau,(k+1)\tau}\lbrack x_{(k+\alpha)\tau}\rbrack - x_{(k+\alpha)\tau}}{(1-\alpha)\tau}.
\end{align*}
Because of the assumption that the flow is one-to-one, we can always write $x_{(k+\alpha)\tau}=\phi_{k\tau,(k+\alpha)\tau}\lbrack \hat x_k \rbrack$ for some~$\hat x_k$. We thus obtain:
\begin{align*}
  \mfrac12\int_0^T\!\big\lvert \dot x_t - v(t,x_t)\big\rvert^2\,dt \approx \sum_{k=0}^{K-1} \frac{\lvert x_{k\tau} - \hat x_k\rvert^2}{2\alpha\tau} + \frac{\lvert x_{(k+1)\tau} - \phi_{k\tau,(k+1)\tau}\lbrack\hat x_k\rbrack \rvert^2}{2(1-\alpha)\tau}.
\end{align*}
Since the number of floaters $\{x_k\super{i}\}_{k=0,\hdots,K; i=1,\hdots,I}$ is large, we can find an~$x\super{i}_k$ close to~$\hat x_k$, giving
\begin{align*}
  \mu_T(x_0\super{i}\to x_K\super{j}) &\approx \inf\Big\{ \mfrac12\int_0^T\!\lvert \dot x_t + v(t,x_t)\rvert^2\,dt : 
    x_0 = x_0\super{i},
    x_T = x_K\super{j},
    x_{k\tau}\in\A_k,\\
      &\hspace{5.6cm} x_k\super{\ell}:=\phi_{(k+\alpha)\tau,k\tau} \lbrack x_{(k+\alpha)\tau)}\rbrack \in\A_k \Big\} \\
    &\approx \min_{i_{(\cdot)}:i_0=i,i_K=j} \sum_{k=0}^{K-1} \min_{\ell=1,\dots,I} \frac{\lvert x_k\super{i_k} - x_k\super{\ell}\rvert^2}{2\alpha\tau} + \frac{\lvert x_{k+1}\super{i_{k+1}} - x_{k+1}\super{\ell}\rvert^2}{2(1-\alpha)\tau} \\
    &= \nu_K(i\smallto j).
\end{align*}

This shows that we can either derive the large-deviation rate function in continuous space and discretize this to finite trajectories (as done here), or we can restrict the continuous dynamics to finite trajectory data and derive a large-deviation rate function for that (as done above); we obtain consistent results whichever route we take.

\subsection{The simple example revisited}
\label{ssec:simpexam2}

Let us now demonstrate how the results of this section apply to the example of Section~\ref{ssec:simpexam1}.

\paragraph{Discrete time and continuous space.}

Let us first suppose we are given infinitely many ``trajectories'' of the system, one starting at each point~$x\in[0,L]$, and they are sampled at discrete time points~$k\tau$, $k=0,1,\ldots,K$, with $\tau = \frac{T}{K}$.  From Section~\ref{subsec:discr traj ldp} with~$\alpha=1/2$ we obtain, by writing~$\Delta x = \frac{L}{K}$, that
\[
\nu_K(0\smallto L) = \mfrac12 \sum_{k=1}^K \frac{(\Delta x)^2}{\tau} = \mfrac12 K\cdot\frac{(L/K)^2}{T/K} = \frac{L^2}{2T}\,,
\]
where we used that the optimal discrete path in~\eqref{eq:discr rate} is the one making jumps of equal lengths~$\Delta x$. Note that the rate function is identical to that in the fully continuous case. Analogously,~$\nu_K(0\smallto L/2) = \nu_K(L/2 \smallto L) = \frac{L^2}{8T}$, and generally, if~$|x-y|=\delta$, then~$\nu_K(x\smallto y)=\frac{\delta^2}{2T}$. The derived semidistances scale similarly. Note that the semidistances converge to zero as~$T\to\infty$.

\paragraph{Discrete time and space.}

If we are given a finite number~$I$ of equispaced trajectories of this system sampled at the same times as in the previous paragraph, the virtual random walker cannot make arbitrarily small jumps as in the continuous state case, thus
\[
  \nu_K(0\smallto L) = \mfrac12\sum_{k=1}^K \frac{(\Delta x_k)^2}{\tau} \approx \frac{L^2}{2T}\,,\qquad\text{if }K\leq I\,,
\]
since we can take~$\Delta x_k \approx L/K$ with error~$\mathcal{O}(I^{-1})$ as~$I$ grows.
However, if~$K>I$, the smallest jumps are~$\Delta x_k = \frac{L}{I}$, thus
\[
  \nu_K(0\smallto L) = \mfrac12I\cdot\frac{(L/I)^2}{\tau} = \frac{L^2\,K}{2I\,T}\,.
\]
Thus, if the observation time of trajectories grows and they are still observed at the same rate~(i.e., $\tau$ stays constant), the semidistances saturate at~$\frac{L^2}{I \tau}$ and do not converge to zero as in the continuous time case.  Moreover, to reach~$y=L/2$ from~$x=0$, we still cannot make smaller jumps than~$\Delta x = \frac{L}{I}$, but now we only require only~$I/2$ of them, such that we obtain~$\nu_K(0\smallto L/2) = \nu_K(L/2\smallto L) = \frac{I}{2}\cdot \frac{(L/I)^2}{\tau} = \frac{L^2}{4I\tau}$ (for even~$I$, and vanishing error for odd~$I$ as~$I$ grows).

The main lesson is, that while in the continuous space case halving the Euclidean distance makes the semidistance scale by~$\frac{1}{4}$, if the spatial resolution of trajectories is coarse, the discrete semidistance scales only by~$\frac{1}{2}$. In general, if~$|x-y|=\delta$, then on a coarse resolution grid it takes about~$\frac{\delta}{\Delta x}$ jumps to travel between these two points, and we obtain~$\nu_K(x\smallto y) \approx \frac{\delta}{\Delta x}\cdot \frac{(\Delta x)^2}{\tau} = \delta\cdot\frac{\Delta x}{\tau}$. Note that~$\Delta x$ and~$\tau$ are constant quantities, and thus the one-way discrete cost scales \emph{linearly} in the Euclidean distance between the two points, as opposed to \emph{quadratic} scaling in the continuous space case.

\section{Coherence analysis with semidistances}
\label{sec:method}

Let us assume that we are given a set of discrete time and space trajectories $\{x_k\super{j}\}_{j=1,\hdots,I,k=0,\hdots,K}$, and a (semi)distance~$d$. We now describe how such semidistances can be used to distinguish and analyze coherent sets from the finite data. We shall work with an unspecified semidistance~$d$, but of course the semidistances that we have in mind are~$\nu_K^\cross$ and~$\nu_K^\meet$ that we derived in the previous section. Other---not large-deviation based---distance measures could be used just as well, as we discuss below. Nevertheless, the semidistances should not be completely arbitrary; we assume that they share the behavior of~$\nu_K^\meet$ and~$\nu_K^\cross$ that we discuss in sections~\ref{subsec:two mixing subdomains} and~\ref{subsec:double gyre}.

%In this section we mostly work with the meeting distance.

To illustrate the ideas we first analyze the behavior of two one-dimensional prototypical examples. These examples show the difference between two types of regions: ``mixing'' and ``static'' (also known as ``regular''). In these one-dimensional and simple examples, one can easily determine the regions and whether they are mixing or static from the semidistances. One example has two invariant sets under the dynamics, which is (measure-theoretically and topologically) mixing on both of them. The other has two static regions, where the mutual physical distance of trajectories does not change under the dynamics, and these regions are separated by a third, mixing region.

After this we proceed with a more involved model: a two-dimensional periodically forced double gyre flow, where the boundaries of the separate regions are no longer as clear-cut as in the one-dimensional example. Nevertheless we will show that one can identify the separate regions via the tools that we present next.

To this end we introduce the notion of cornerstones, representing possible coherent sets or mixing regions, then discuss how to find them and when to stop searching for them. Finally, to obtain coherent sets, we assign the trajectories to cornerstones. The notion of fuzzy affiliations will be used to express the uncertainty whether a trajectory close to the boundary of a set belongs to it or not. Note that in the case of finite data such an uncertainty is always present.

\subsection{Two illustrative model cases}
\label{subsec:two mixing subdomains}

\paragraph{Two invariant, mixing subdomains.}

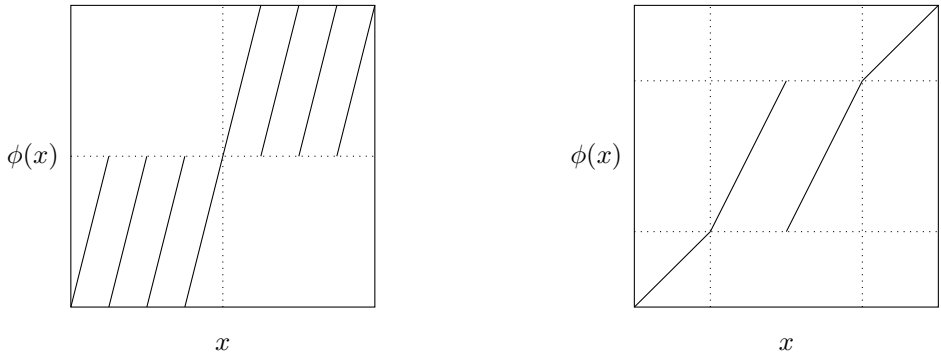
\begin{figure}[h!]
\centering
\begin{tikzpicture}[scale=1]
  \tikzstyle{every node}=[font=\small]
  \draw(0,0) rectangle (4,4);
%  \draw[dotted] (0.5,0)--(0.5,4);
%  \draw[dotted] (1,0)--(1,4);
%  \draw[dotted] (1.5,0)--(1.5,4);
  \draw[dotted] (2,0)--(2,4);
%  \draw[dotted] (2.5,0)--(2.5,4);
%  \draw[dotted] (3,0)--(3,4);
%  \draw[dotted] (3.5,0)--(3.5,4);
  \draw[dotted] (0,2)--(4,2);
  \draw(0,0)--(0.5,2);
  \draw(0.5,0)--(1,2);
  \draw(1,0)--(1.5,2);
  \draw(1.5,0)--(2.5,4);
  \draw(2.5,2)--(3,4);
  \draw(3,2)--(3.5,4);
  \draw(3.5,2)--(4,4);
  \node at (2,-0.5) {$x$};
  \node at (-0.5,2) {$\phi(x)$};
\end{tikzpicture}
\qquad\qquad\qquad
\begin{tikzpicture}[scale=1]
  \tikzstyle{every node}=[font=\small]
  \draw(0,0) rectangle (4,4);
  \draw[dotted] (0,1)--(4,1);
  \draw[dotted] (0,3)--(4,3);
  \draw[dotted] (1,0)--(1,4);
  \draw[dotted] (3,0)--(3,4);
%  \draw[dotted] (2.5,0)--(2.5,4);
%  \draw[dotted] (3,0)--(3,4);
%  \draw[dotted] (3.5,0)--(3.5,4);
%  \draw[dotted] (0,2)--(4,2);
  \draw(0,0)--(1,1);
  \draw(1,1)--(2,3);
  \draw(2,1)--(3,3);
  \draw(3,3)--(4,4);
%  \draw(2.5,2)--(3,4);
%  \draw(3,2)--(3.5,4);
%  \draw(3.5,2)--(4,4);
  \node at (2,-0.5) {$x$};
  \node at (-0.5,2) {$\phi(x)$};
\end{tikzpicture}

\caption{Left: The time-discrete flow map with two invariant mixing subdomains. Right: The time-discrete flow map with two static regions and an invariant mixing subdomains between them.}
\label{fig:two mixing subdomains}
\end{figure}

As mentioned in Section~\ref{subsec:discr semidistances}, we may also apply the techniques developed in this paper to a discrete-time dynamical system. To gain some intuition for the behavior of the semidistances at mixing regions, we consider the discrete-time system on the unit interval~$X = [0,1]$ and one-step flow map, see Figure~\ref{fig:two mixing subdomains} (left),
\[
{
\renewcommand{\arraystretch}{1.25}
\phi(x) = \left\{ \begin{array}{ll}
4x \mod \frac12, & x< \tfrac12 \\
\left(4(x-\tfrac12) \mod \tfrac12\right) + \tfrac12, & x \ge \tfrac12\,.
\end{array}
\right.
}
\]

The sets~$X_1 = [0,\tfrac12]$, $X_2 = (\tfrac12,1]$ are invariant, i.e.,~$\phi^{-1}(X_1) = X_1$ and~$\phi^{-1}(X_2) = X_2$, and~$\phi$ is simply the circle-quadrupling map on each of these sets, i.e., it is mixing on the single components. Consequently\footnote{If a system~$(X,\psi,\mu)$ is mixing, then~$(X\times X,\psi\times\psi,\mu\times\mu)$ is ergodic~\cite[Theorem~1.24]{Wal00}. Thus, for~$\mu\times\mu$-almost every pair~$(x,y)$, the trajectory~$(\psi\times\psi)^t(x,y)$ will enter every set~$A$ of non-zero measure for some~$t\ge 0$. This shows~\eqref{eq:mixmeet} by taking~$\psi = \phi\vert_{(0,1/2)}$ or~$\psi = \phi\vert_{(1/2,1)}$, and $A = \{(x,y)\,\vert\, |x-y|<\ep\}$ for any fixed~$\ep>0$. \label{foot:mix}},
\begin{equation}
\liminf_{t\in\N,\, t\to\infty} \left|\phi^t(x) - \phi^t(y)\right| = 0
\label{eq:mixmeet}
\end{equation}
for (Lebesgue-)almost every pair~$x,y\in X_i$, $i=1,2$.

\begin{figure}[h]
\centering
\includegraphics[width = 0.8\textwidth]{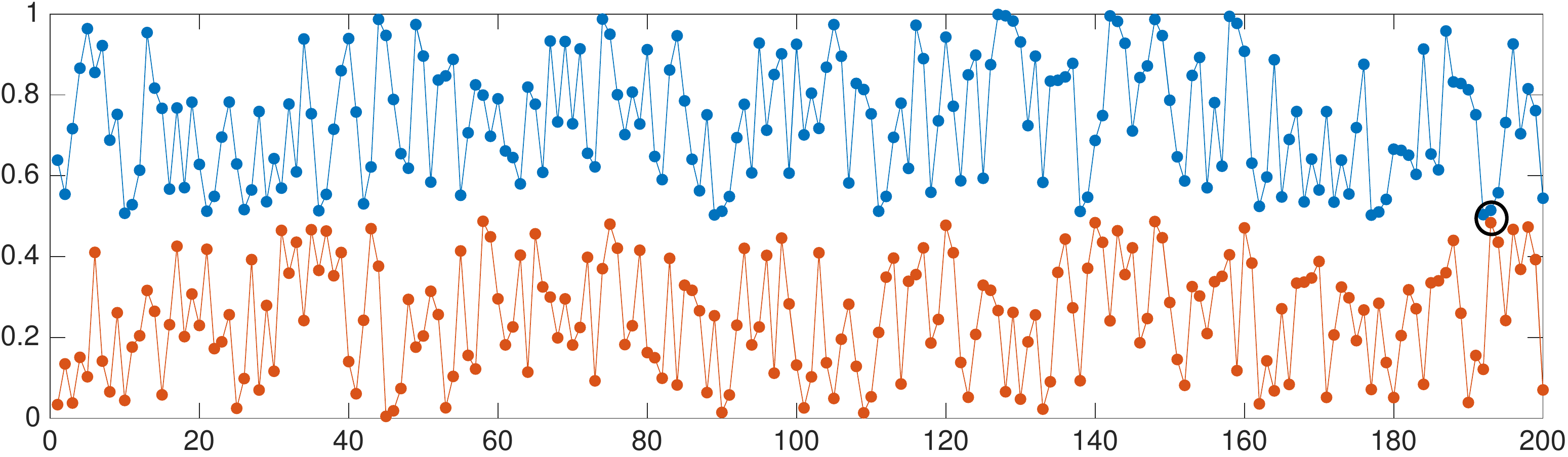}

\caption{Two trajectories of the map~$\phi$ of length 200 steps, starting in $X_1$ and $X_2$, respectively. Theory shows that they come arbitrary close, eventually. Here they get the closest at time step 193, shown by a circle.}
\label{fig:TwoMixing_longtraj}
\end{figure}

\begin{figure}[h]
\centering
\includegraphics[width = 0.45\textwidth]{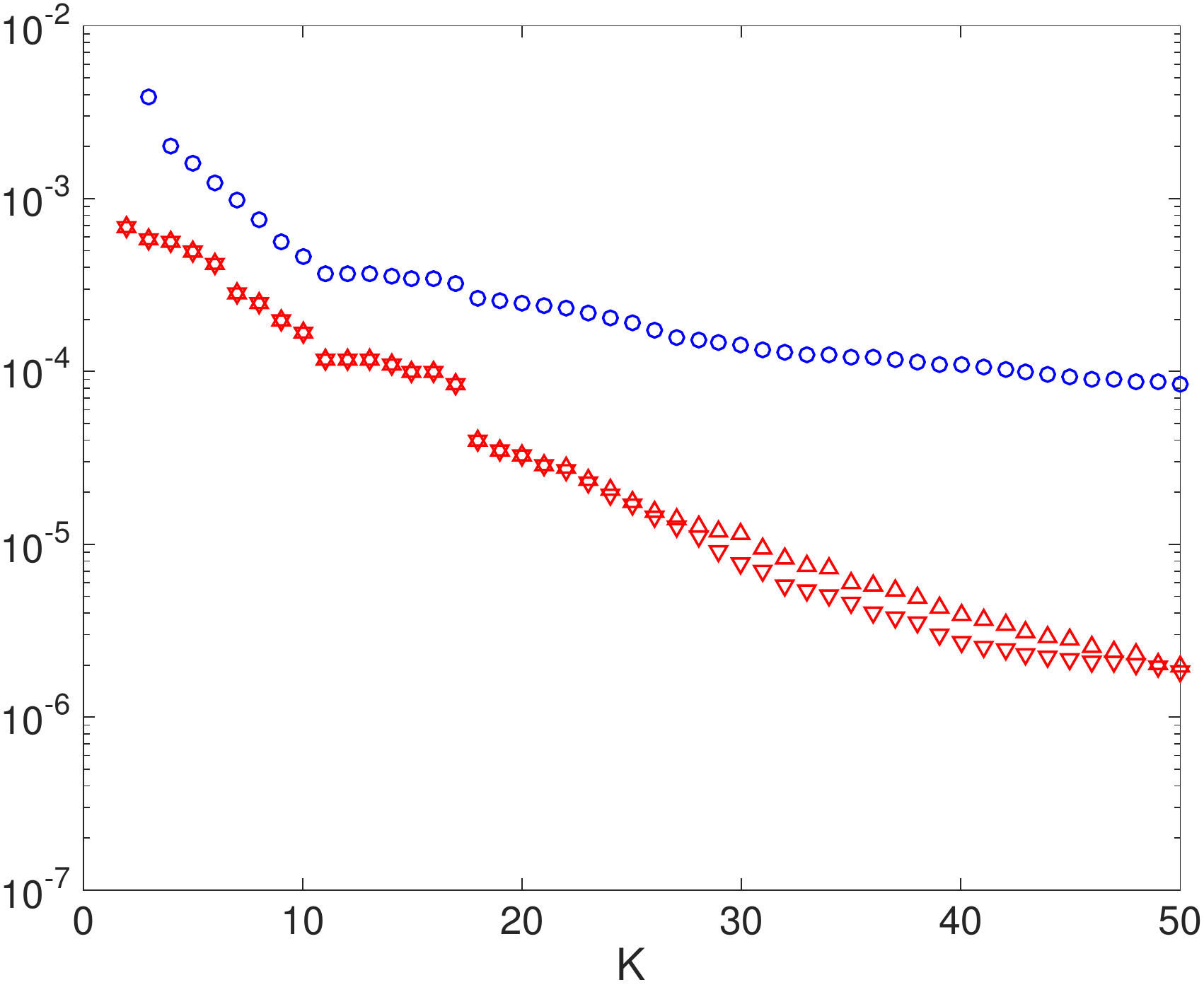}
\hfill
\includegraphics[width = 0.45\textwidth]{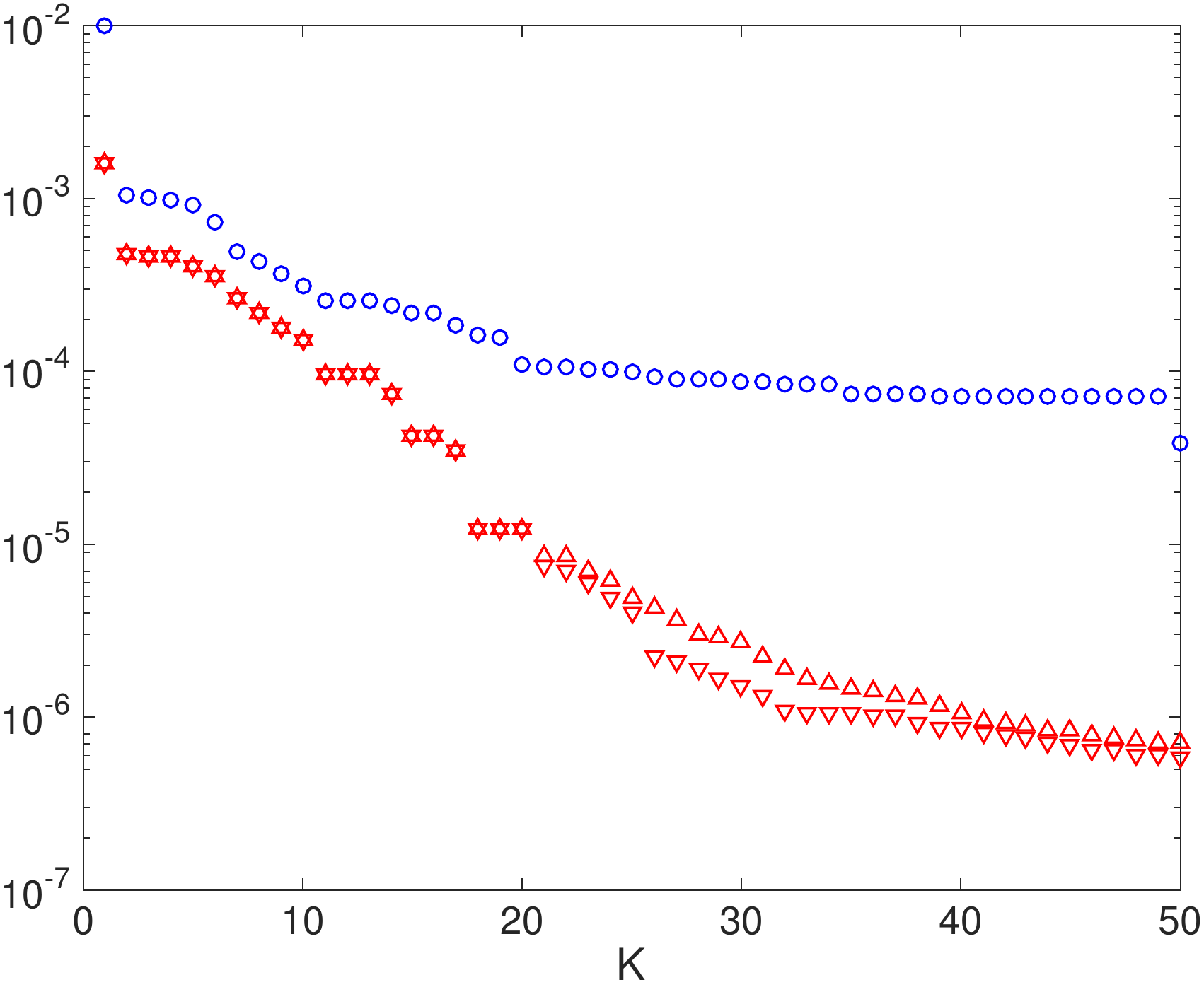}

\caption{Semidistances~$\nu_K^\cross(i\smallto j)$ (left) and~$\nu_K^\meet(i,j)$ (right) for increasing maximal time~$K$, averaged over~$x_i,x_j\in X_1$ (downward-pointing triangles), $x_i,x_j\in X_2$ (upward-pointing triangles), and $x_i,\in X_1,\,x_j\in X_2$ (circles), respectively. Note that the decrease of the distance is much slower for trajectories taken from different invariant sets.}
\label{fig:TwoMixing_avgdists}
\end{figure}

Thus,~$\nu_K(i\smallto j)\to 0$ as~$K\to\infty$, because if~$x_i,x_j$ are both in~$X_1$ or both in~$X_2$, then~\eqref{eq:mixmeet} shows that their trajectories get arbitrarily close eventually. If the trajectories start in different halves of~$[0,1]$, then\footnote{Applying Footnote~\ref{foot:mix} to the case where $\psi=\phi\vert_{(0,1/2)}$, we obtain that $(\phi\vert_{(0,1/2)} \times \phi\vert_{(0,1/2)} )^t(x,y)$ enters $A = \{(x,y)\,\vert\, |x-1/2| + |y| < \ep \}$ eventually. Noting that $\phi\vert_{(1/2,1)} = \phi\vert_{(0,1/2)} (\cdot - \tfrac12) + \tfrac12$, the claim follows. Comparing the different slopes in Figure~\ref{fig:TwoMixing_avgdists}, based on the reasoning in Footnote~\ref{foot:mix} and here we conjecture that the minimal distance of two trajectories decays faster in the case when they both start in the same invariant set, because the set~$\{(x,y)\,\vert\, |x-y|<\ep\}$ is larger in measure than~$\{(x,y)\,\vert\, |x-1/2| + |y| < \ep \}$.}
\begin{equation}
\liminf_{t\in\N,\, t\to\infty} \left| \phi^t(x_i) - \tfrac12 \right| + \left|\phi^t(x_j) - \tfrac12 \right| = 0\,,
\label{eq:mixmeet2}
\end{equation}
thus the jump from one trajectory to another gets arbitrarily cheap. See Figure~\ref{fig:TwoMixing_longtraj}.

The transport semidistances between any two points within the same region are very small, at least if the time window is large enough. This behavior is typical for mixing regions. In fact, since the two mixing regions are only separated by one point, it is relatively cheap to move from one region to the other, and so the semidistances between two points in separate regions converge with increasing time to zero. Nevertheless, the semidistances still detect a difference between the two invariant sets: the semidistance between two trajectories in the same invariant component goes in general quicker to zero than the one between two from different components, as shown in Figure~\ref{fig:TwoMixing_avgdists} for~$I=100$ initially equispaced trajectories. Thus, it is the relative difference between the semidistances that is relevant for the transport-structure of the state space, and not the absolute values.

\paragraph{Two static regions divided by a mixing one.}

To gain some intuition about static regions, let us now consider the discrete-time system on~$X=[0,1]$ given by
\[
{
\renewcommand{\arraystretch}{1.25}
\phi(x) = \left\{ \begin{array}{ll}
x, & x \in [0,\tfrac14) \cup (\tfrac34,1] \\
\left(2(x-\tfrac14) \mod \tfrac12\right) + \tfrac14, & x \in [\tfrac14,\tfrac34]\,,
\end{array}
\right.
}
\]
see Figure~\ref{fig:two mixing subdomains} (right). This map has three invariant sets. The left and right ones are static, such that the mapping restricted to them is the identity, and are meant to model regions of the state space in complicated flows that are ``static'' in the sense that the mutual distance of points is not changed (or just barely) by the dynamics. We will consider these as one kind of prototype for coherent sets. The third region is mixing, and physically separates the other two.

\begin{figure}[h]
\centering
\includegraphics[height = 58mm]{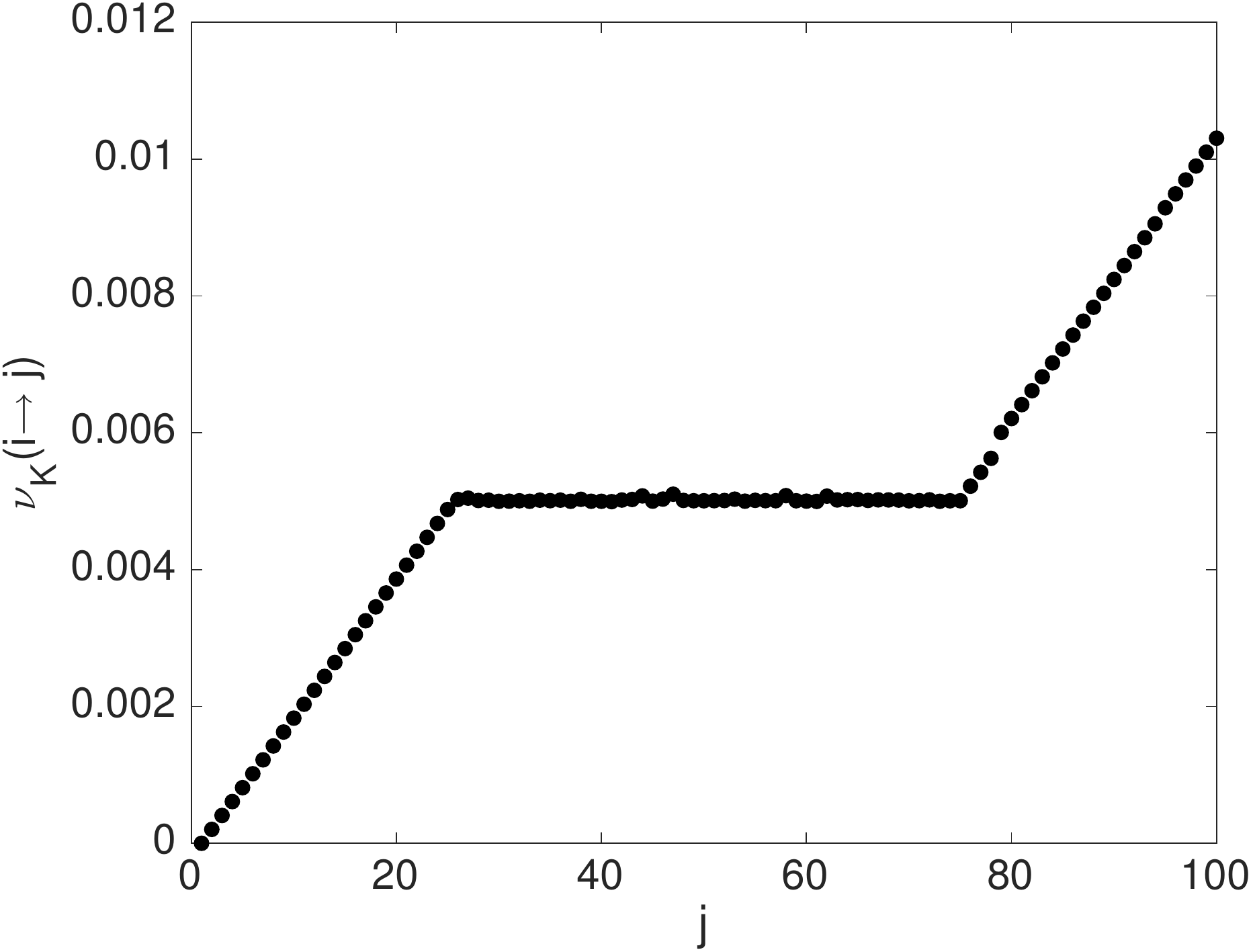}
\hfill
\includegraphics[height = 60mm]{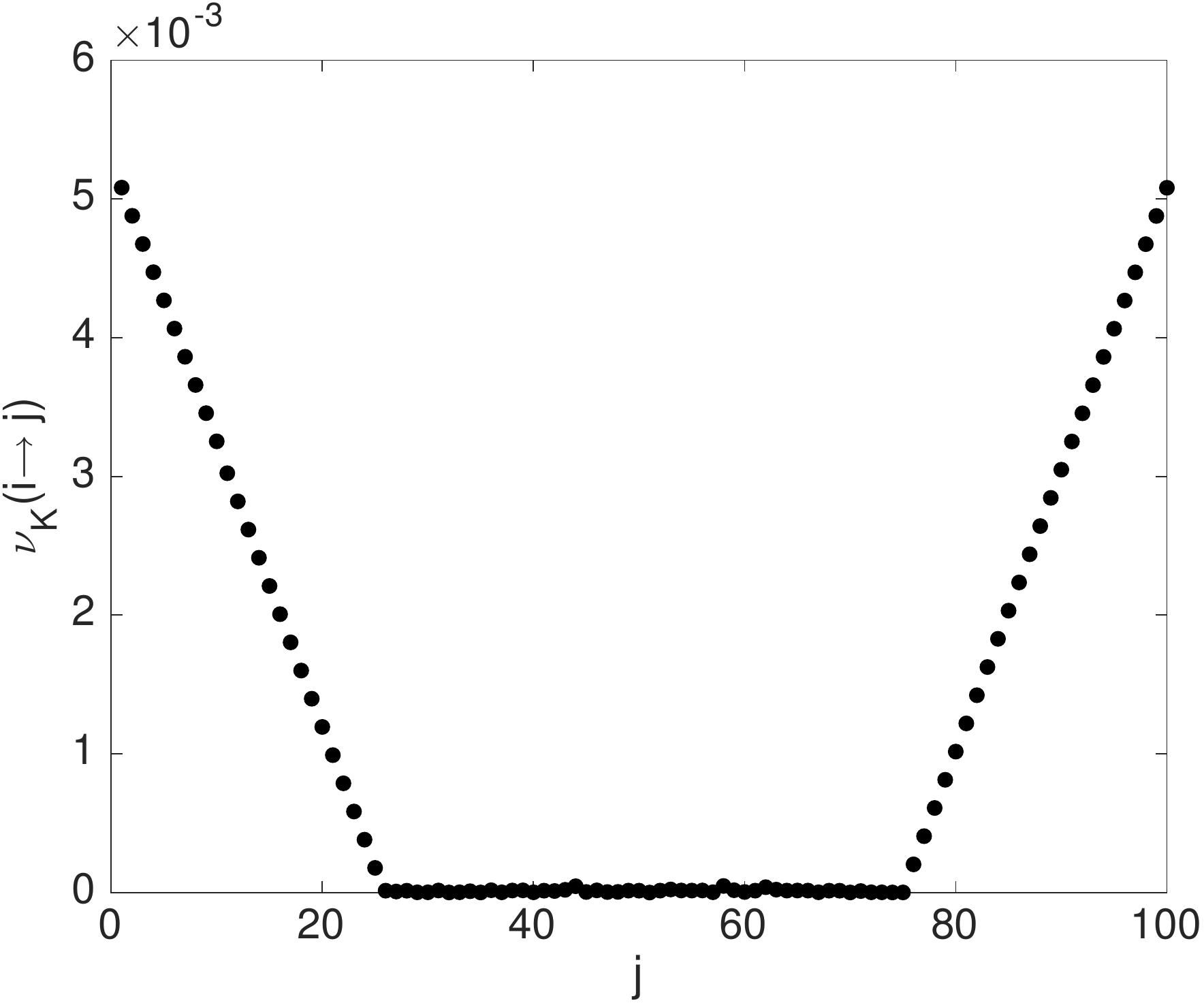}

\caption{One-way cost~$\nu_K(i\smallto j)$ for~$i=1$ (left) and~$i=51$ (right) for the map with two static and one mixing region.}
\label{fig:stay mix stay costs}
\end{figure}

We take~$I=100$ initially equispaced trajectories and compute the one-way costs~$\nu_K(i\smallto \cdot)$ with~$K=50$ for~$i=1$ and~$i=51$, respectively, shown in Figure~\ref{fig:stay mix stay costs}.
From our analysis in Section~\ref{ssec:simpexam2} we would have expected to see quadratic growth of the one-way cost with respect to physical distance in the static regions, but we only observe linear growth. This is due to the finite number of considered trajectories, as also explained in the second paragraph of Section~\ref{ssec:simpexam2}. All points of the mixing region have almost the same cost from any one point in the static regions, and approximately zero cost from one another. To obtain the cost between two points of different static regions, one has to consider the cost to go to the boundary of the static and mixing regions (linear cost in Euclidean distance), travel on a trajectory from there to the boundary of the other static region (at zero cost), and then go from there to the desired point (again, linear cost in Euclidean distance that needs to be covered). Thus, the cost (and semidistance) between these two points is the sum of their one-way cost (and semidistance) to the mixing region, provided the time of consideration is sufficiently large for the mixing to take place.\footnote{Note that for this argument, ergodicity of the dynamics in the ``mixing'' region would be sufficient, since one only needs to travel ``from one static region to the other''. The crucial additional property we get from mixingness is that the mutual semidistances of points in this region go to zero.} Since our fictive random walker uses trajectories of the mixing region to travel from one static region to the other, we will also call it \emph{transition region} henceforth.

To conclude, from Figure~\ref{fig:stay mix stay costs} we can easily identify three separate regions, and from the steepness of the slopes (linear/quadratic or flat), we can determine wether a region is static or mixing. As the next example shows, this distinction is usually not as clear as in these constructed examples, but the main ideas will be based on this observation.

\subsection{The periodically forced double gyre}
\label{subsec:double gyre}

\begin{figure}[h!]
\centering
\begin{tikzpicture}[scale=1.2]
\tikzstyle{every node}=[font=\small]
  \draw(0,0) rectangle (4,2);
%  \draw(0.7,1) ellipse (0.4 and 0.8);
  \draw[->](0.4+0.7,1-0.14) arc (350:0:0.4 and 0.8);
 %  \draw(2.7,1) ellipse (1 and 0.5);
  \draw[->](1.715,0.915) arc (-170:180:1 and 0.5);
  \draw(5,0) rectangle (9,2);
%  \draw(9-2.7,1) ellipse (1 and 0.5);
  \draw[->](7.29,0.91) arc (350:0:1 and 0.5);
%  \draw(9-0.7,1) ellipse (0.4 and 0.8);
  \draw[->](7.91,0.86) arc (-170:180:0.4 and 0.8);
  \end{tikzpicture}
\caption{Sketch of the velocity field of the periodically forced double gyre flow at two different times. The horizontal axis is~$x_1$, the vertical is~$x_2$.}
\label{fig:per double gyre}
\end{figure}
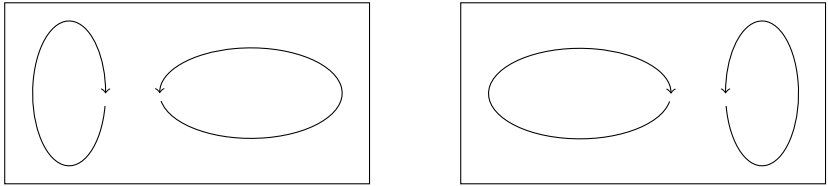

Let us now consider the non-autonomous system $\dot x_t = v(t,x_t)$ on $X = [0,2]\times [0,1]$ with~\cite{FrPa14}
\begin{equation}
\begin{aligned}
 v(t,x):=
  \begin{bmatrix}
    -\pi A \sin\left(\pi f(t,x_1)\right) \cos(\pi x_2) \\
    \pi A \cos\left(\pi f(t,x_1)\right) \sin(\pi x_2)\frac{df}{dz}(t,x_1)
  \end{bmatrix},
\end{aligned}
\label{eq:DoubleGyre}
\end{equation}
where $f(t,z)=\beta \sin(\omega t)z^2+(1-2\beta\sin(\omega t))z$. We fix the parameter values $A=0.25$, $\beta=0.25$ and~$\omega=2\pi$, hence the vector field has time period~1. The system preserves the Lebesgue measure on~$X$. Equation~\eqref{eq:DoubleGyre} describes two counter-rotating gyres next to each other (the left one rotates clockwise), with the vertical boundary between the gyres oscillating periodically, see Figure~\ref{fig:per double gyre}.

We choose a uniform $50\times 25$ grid as initial conditions for the floaters at time~$t=0$; i.e.,~$I=1250$. We sample the trajectories of these floaters at times~$t_k = k\tau$, $k=0,1,\ldots,K$, where~$K=100$ and~$\tau = 0.2$. That means, the length of trajectories in consideration is $20$ times the period of the forcing.

Employing our large-deviation based distance computations on this data set using Algorithm~\ref{algo1} and $\alpha=1/2$, we get the one-way costs~$\nu_K(i\smallto j)$, $i,j=1,\ldots,I$, from which we compute $\nu_K^\cross(i,j)$ and~$\nu_K^\meet(i,j)$.

\begin{figure}[h]
\centering
\includegraphics[width = 0.4\textwidth]{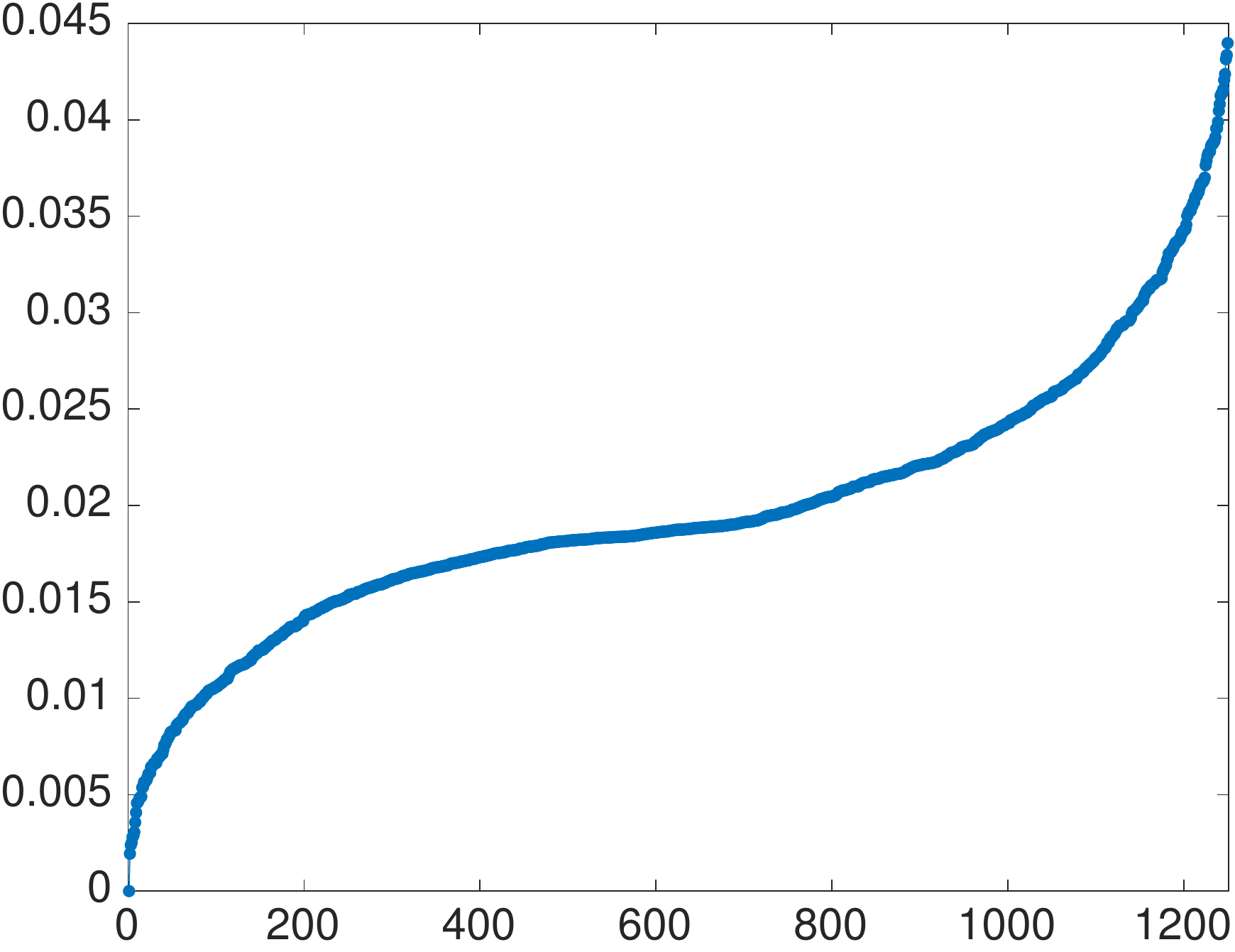}
\hfill
\includegraphics[width = 0.4\textwidth]{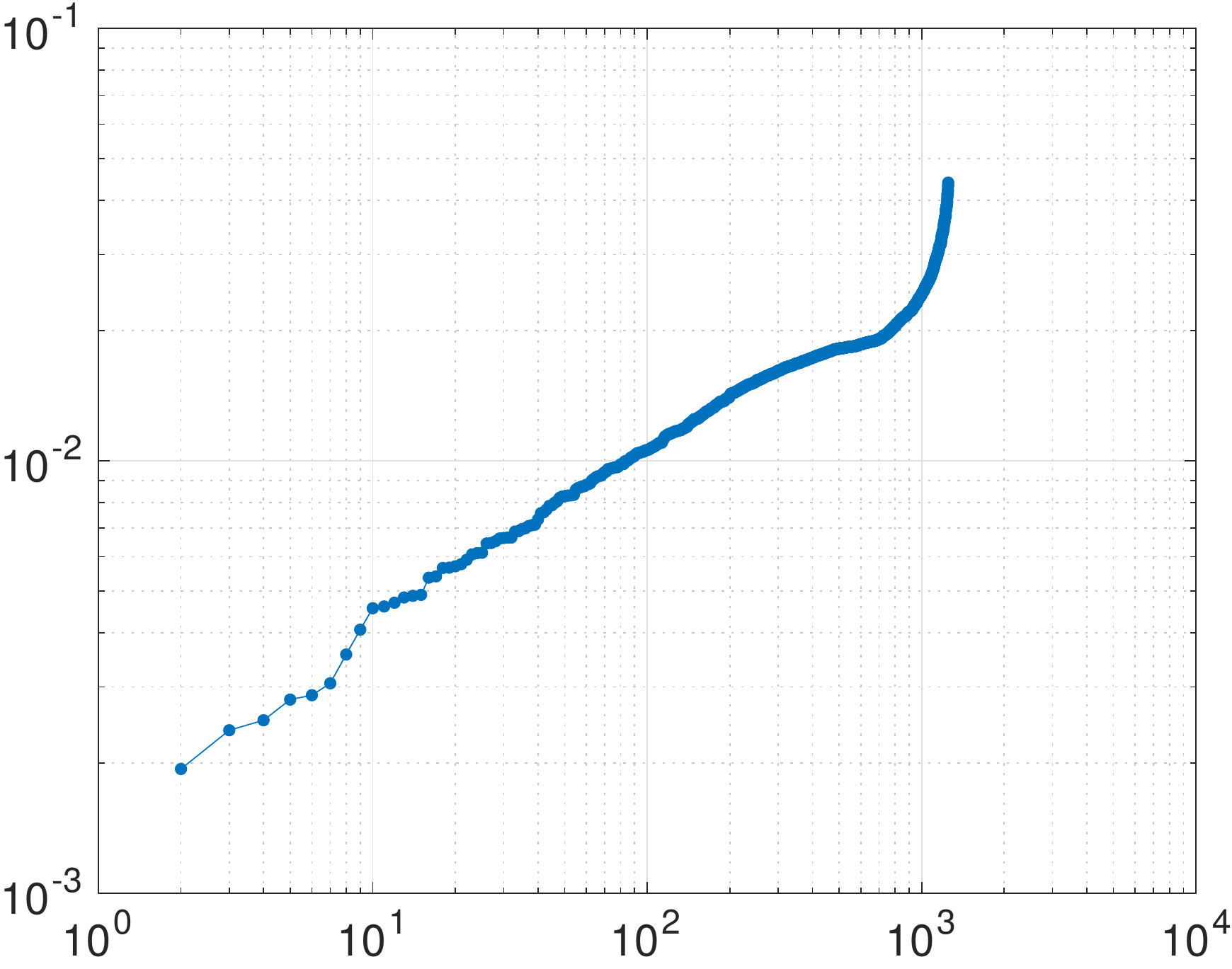}
\caption{Left: $\nu_K^\meet(c_1,\cdot)$ sorted in ascending order. Right: the same as left, on a log-log scale. The horizontal axis shows the rank, the vertical shows the semidistance value.}
\label{fig:DG_meetdist_loglog}
\end{figure}

As a first simple analysis, we can order the points by their semidistances to the center of one gyre, see Figure~\ref{fig:DG_meetdist_loglog}. Here and in the following, the rates and semidistances will be always given in units~$1/\tau$. On a log-log scale, the slope~$1/2$ (square-root-type behavior) indicates that most trajectories in the gyre are approximately concentric circular regions around the center.\footnote{As on a regular grid there are~$\mathcal{O}(\delta^2)$ points not further than Euclidean distance~$\delta$ from a reference point, the $r$-th closest point to the reference point has distance~$\mathcal{O}(r^{1/2})$.} Since the semidistances grow linearly in the Euclidean distance inside the gyre, we see that we are in the low-resolution regime discussed in Section~\ref{ssec:simpexam2}. Analogously to Figure~\ref{fig:stay mix stay costs}, we can again (vaguely) distinguish three regions: a steep (square-root-type) region, a flat region, and another steep (flipped square-root-type) region. As before, the flat region is typically strongly mixing, and the steep regions are static. We shall make this distinction more precise in the next sections.

\begin{remark}[Three-dimensional flows] \label{rem:3d case}
Clearly, the scaling behavior shown in Figures~\ref{fig:stay mix stay costs} and~\ref{fig:DG_meetdist_loglog} are dimension-dependent. For a three-dimensional static region one would see a slope~$1/3$ on a log-log scale. The question remains, how does a typical coherent set behave there; can it be modeled by a static region? If an incompressible flow rotates uniformly in a plane, it necessarily has a constant shifting motion in the, perpendicular axial direction, leading to cylindrical vortices. If the cylindrical rings of a vortex rotate at different angular frequencies, the flow speed along axial directions is nonuniform, and mixing-type behavior occurs in the vortex~\cite{HalEtAl07}. We leave the analysis of such systems to future work, and proceed with analyzing different aspects of prototypical two-dimensional flows here.
\end{remark}

\subsection{Cornerstones}

To start the analysis of the state space under a semidistance $d$, we randomly choose a trajectory, represented by a label $c_0\in\{1,\hdots,I\}$, and compute the trajectory furthest from it, i.e, we set
\begin{equation*}
  c_1 = \argmax_{i=1,\ldots,I} d(i,c_0)\,.
\end{equation*}
To find a set of points that ``spans'' the state space, we identify successively further trajectories that are far away from \emph{all} the other already identified ``cornestones''~$\{c_{q}\}_{q=1,\hdots,Q}$, as in~\cite{RuSaSch17}:
\begin{equation}
  c_{Q+1} = \argmax_{i=1,\ldots,I}\min_{q=1,\ldots,Q} d(i,c_q)\,.
\label{eq:coresearch}
\end{equation}
Observe that in this optimization problem we ignore the first, randomly chosen trajectory~$c_0$; hence the set of cornerstones~$\{c_q\}_{q=1,\hdots,Q}$ will be less dependent on this randomness. Moreover, even if the first trajectory~$c_0$ would represent a coherent set, the algorithm will eventually provide a new cornerstone in that set, which lies closer to the semidistance center of that set.

\begin{figure}[h]
\centering
\includegraphics[width = 0.49\textwidth]{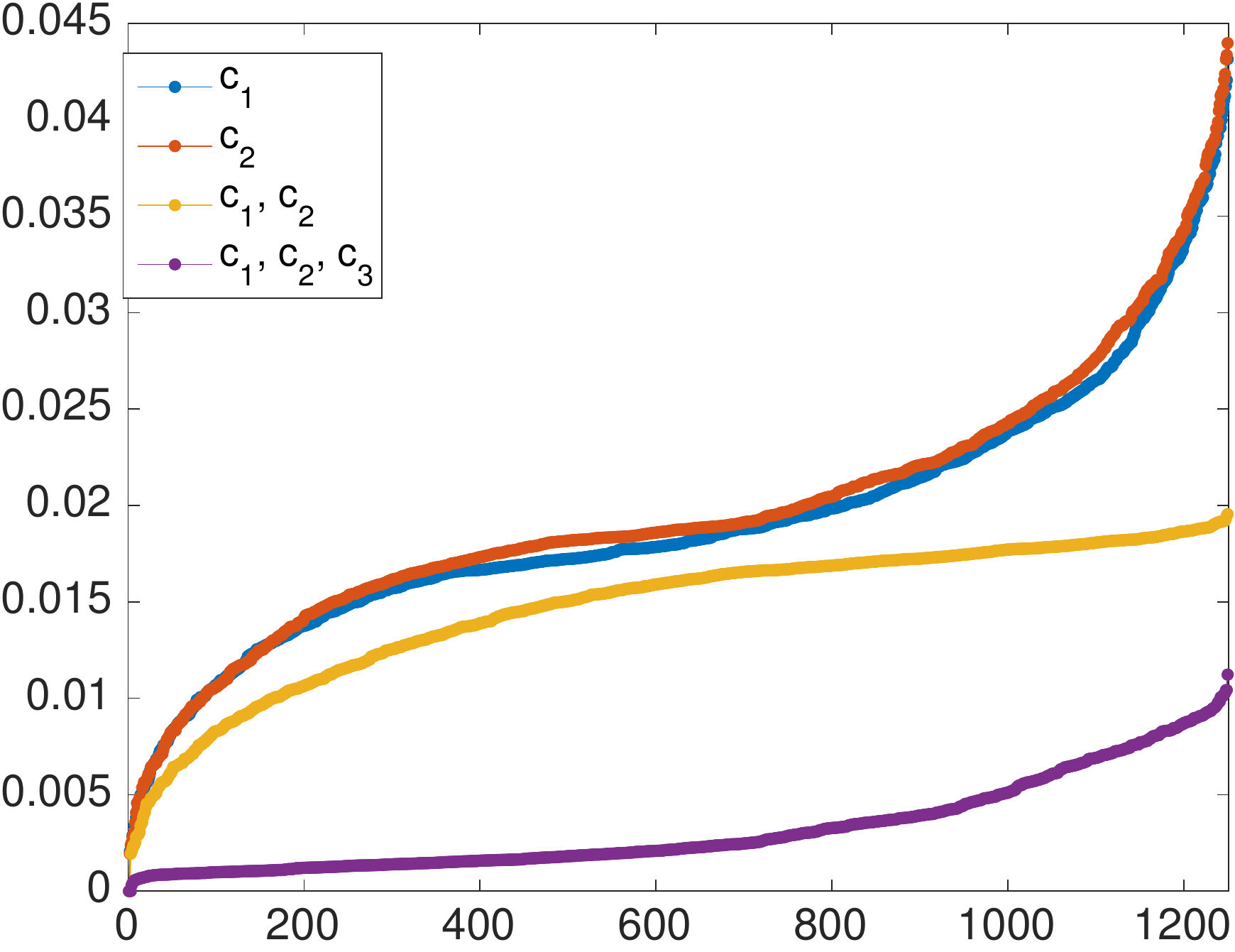}
\caption{The objective functions of the maximization problem in~\eqref{eq:coresearch}, sorted in ascending order (yellow and purple). Blue and red: $\nu_K^\meet(c_i,\cdot)$, $i=1,2$.}
\label{fig:DG_meetdist_fromcores}
\end{figure}

For the double gyre and the meeting distance, we identified three cornerstones. The objective function of the maximization problem~\eqref{eq:coresearch} is plotted in Figure~\ref{fig:DG_meetdist_fromcores}; this yields a similar but more detailed picture as Figure~\ref{fig:DG_meetdist_loglog}. Note that the chaotic, well mixed transition region appears as flat region in these distance graphs, and the gyres appear as steep regions towards the maxima of the respective graphs. That the chaotic region is well mixed, and has no stratification (invariant rings as the gyres), can be seen from its flat behavior towards its maximum. The forth cornerstone is part of a gyre, it starts to stratify it. Nevertheless, its distance to the other corners is much smaller.

\begin{figure}[h!]

\begin{minipage}[b]{0.5\textwidth}
\includegraphics[height = 35mm]{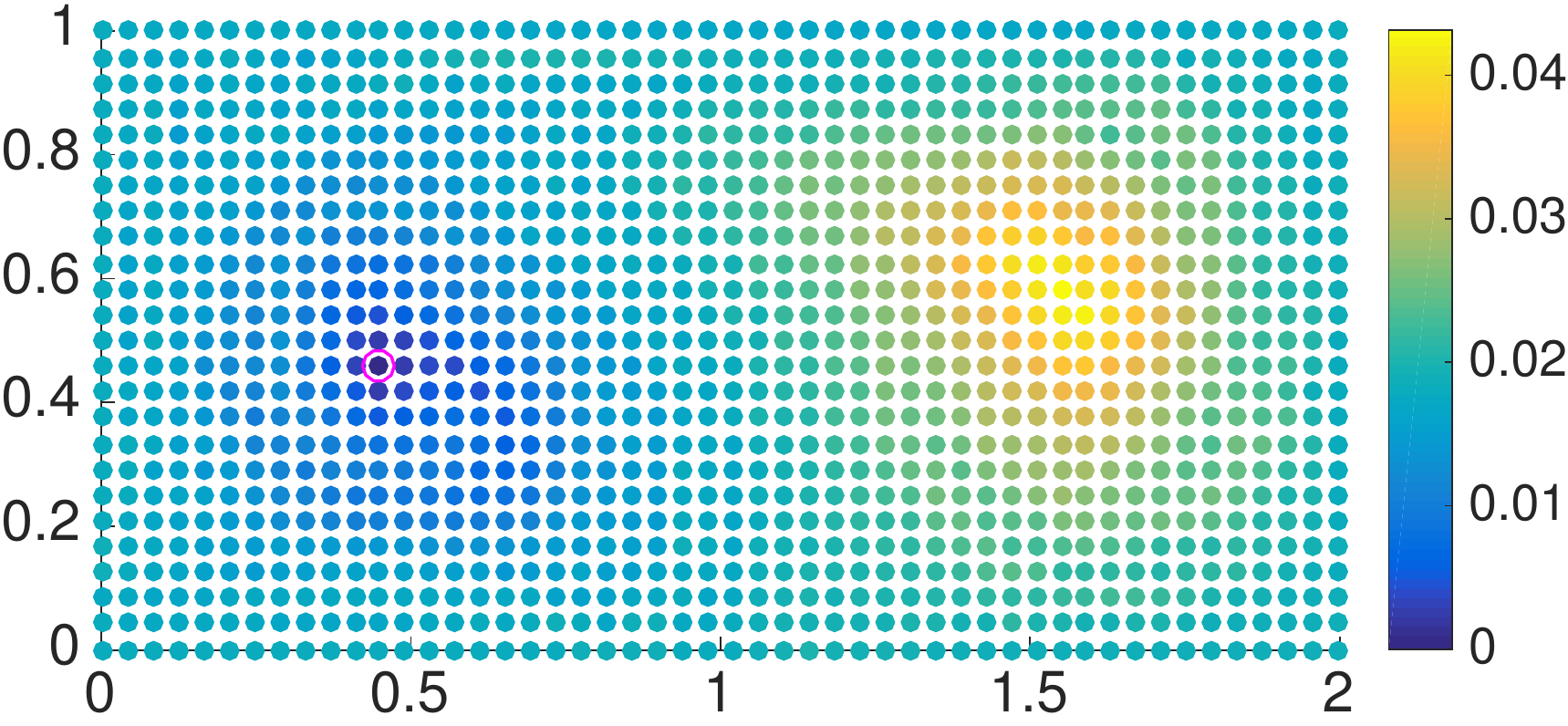}
\end{minipage}
\
\includegraphics[height = 35mm]{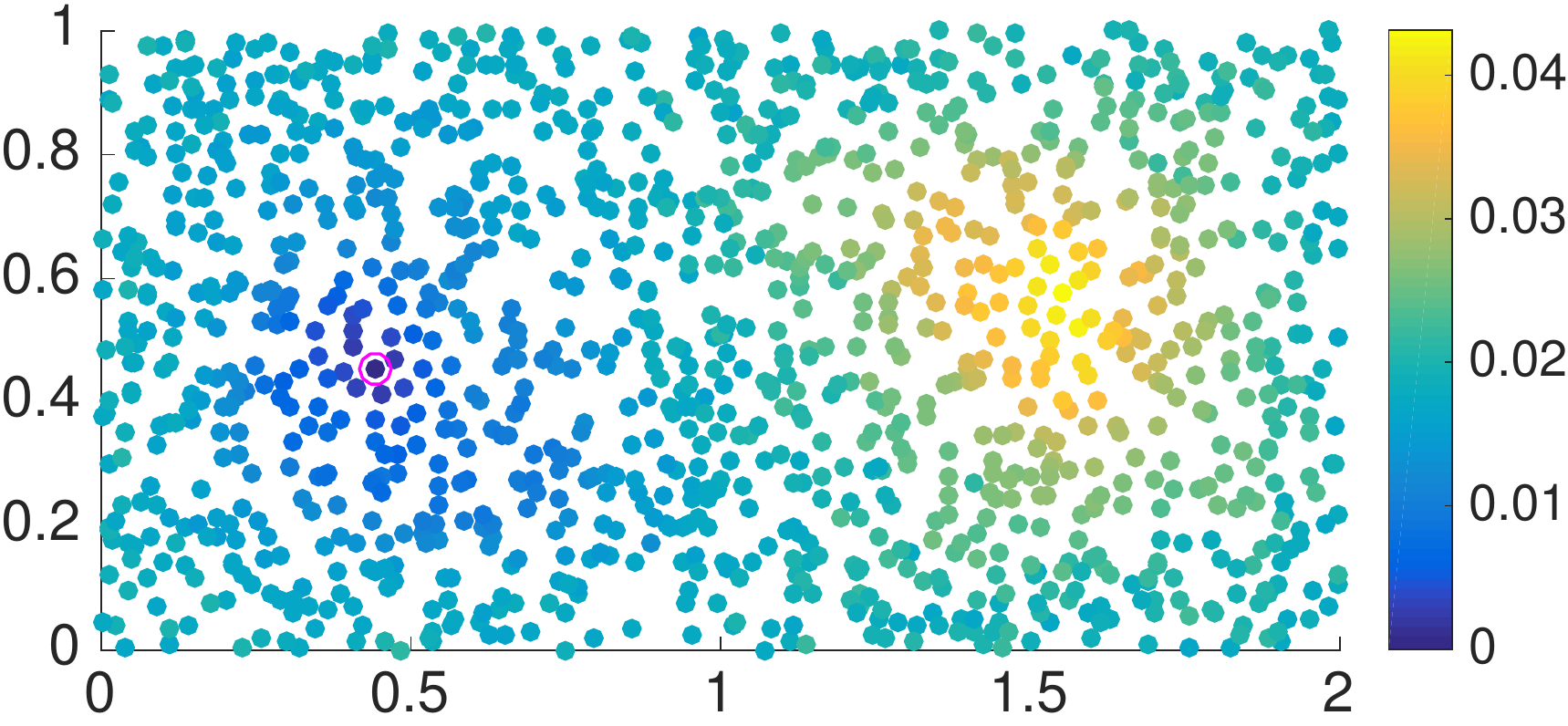}

\begin{minipage}[b]{0.5\textwidth}
\includegraphics[height = 35mm]{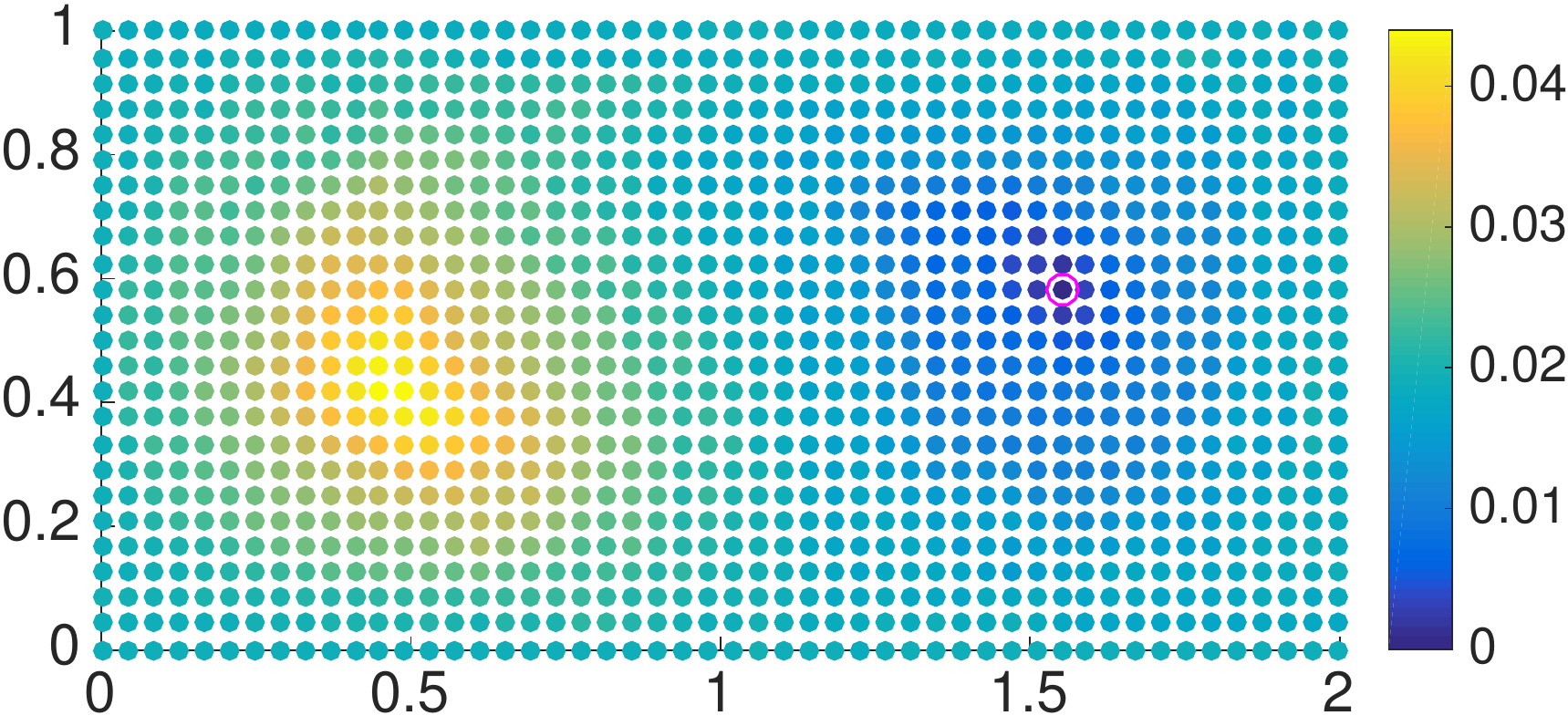}
\end{minipage}
\
\includegraphics[height = 35mm]{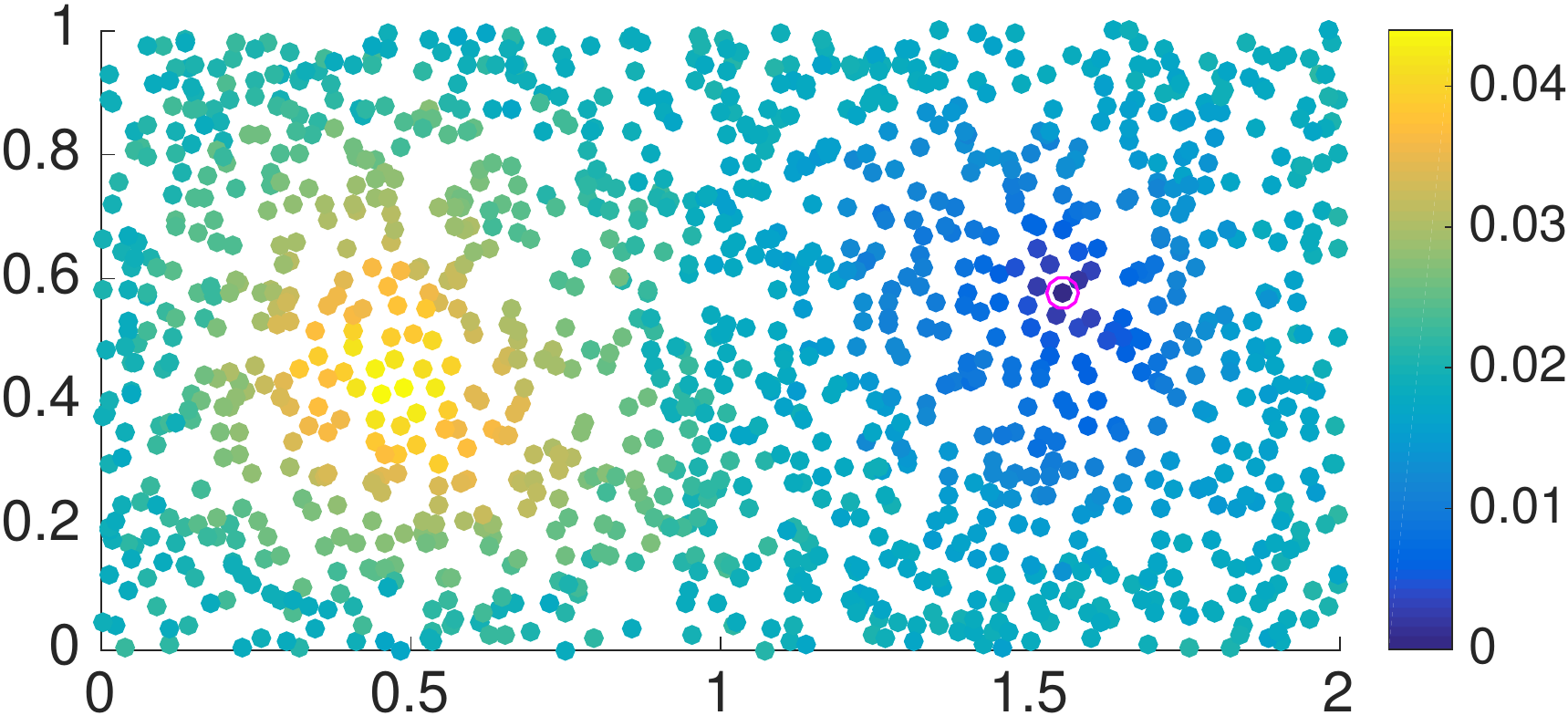}

\begin{minipage}[b]{0.5\textwidth}
\includegraphics[height = 35mm]{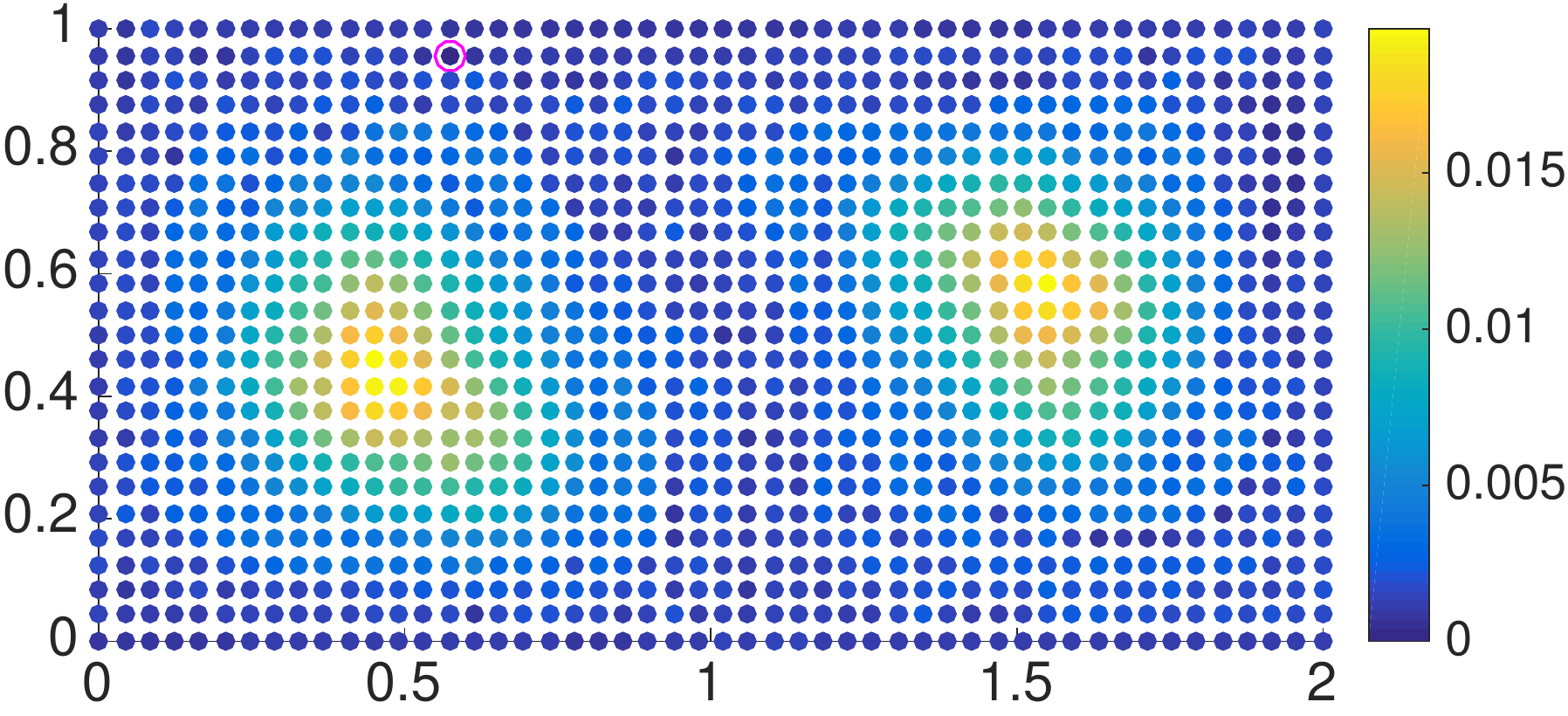}
\end{minipage}
\
\includegraphics[height = 35mm]{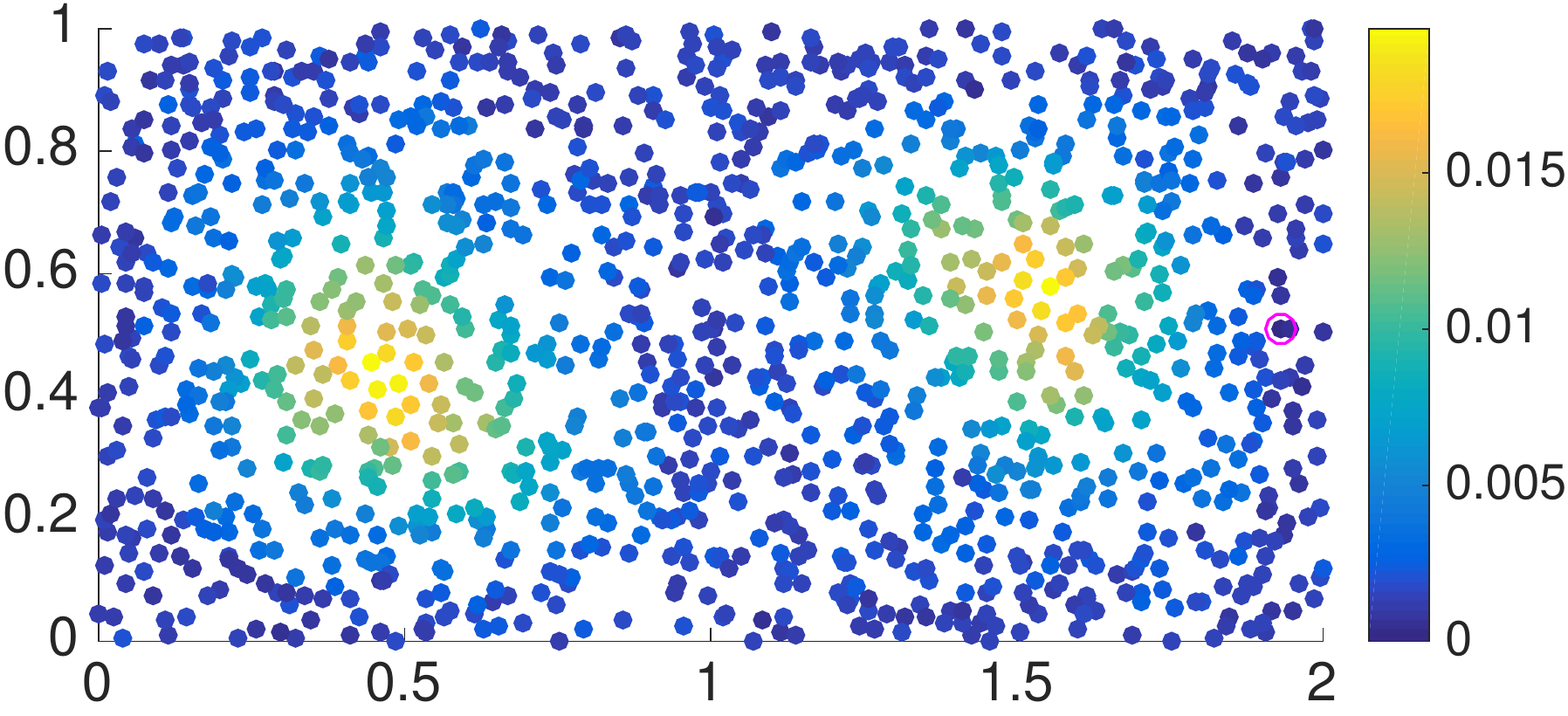}
\caption{Distances $\nu_K^\meet$ from trajectory $c_1$ (top), $c_2$ (middle) and $c_3$ (bottom), marked by the magenta circle, at initial (left) and final times (right). The semidistances are given in units~$1/\tau$. The horizontal axis is~$x_1$, the vertical is~$x_2$.}
\label{fig:DG_meetdist_c1c2c3}
\end{figure}

To get a first glimpse of the separate regions in the state space we have plotted the semidistances from each cornerstone in Figure~\ref{fig:DG_meetdist_c1c2c3}, both at the initial and final times. Note that since we work with trajectory labels rather than physical positions, the semidistances are invariant in time, whereas the physical positions of the floaters change over time. From these figures, one can approximately identify the two static (gyre) regions, being very close and very far from~$c_1$ and $c_2$ respectively, and the chaotic transition region in between, having approximately constant distance from~$c_1$ and~$c_2$, cf.~\cite[Figure~1]{FrPa14}.

\subsection{Number of cornerstones}
\label{ssec:num cornerstones}

How to determine the number of cornerstones that should be used? Is there an optimal number, or is it up to our liking? In the case of the double gyre, as noted above, a fourth cornerstone would be part of one of the gyres, and assigning affiliations would thus split one gyre into two sets. If the gyres would consist of a continuum of periodic orbits, then we could proceed and split them this way into as many rings as we like. The same situation in an idealized framework appears in Section~\ref{subsec:two mixing subdomains} for the static regions: since they are static, arbitrary subsets are perfectly coherent (even invariant in this case).

A good place to stop searching for further cornerstones would be when they would start to subdivide ``maximal coherent'' sets, as the gyres in the double-gyre example, or the static sets in the second, and the invariant sets in the first example of Section~\ref{subsec:two mixing subdomains}. To this end we make an idealized assumption: \emph{Coherent sets appear as for the second example in Section~\ref{subsec:two mixing subdomains}, i.e., multiple static regions divided by one mixing region}. Here, ``static'' is meant in the sense that the mutual distances between points in the set barely change. Such an assumption was also utilized in~\cite{HaEtAl16}.

Note that if there are~$C\ge 2$ coherent sets, the first~$C$ corner stones are going to be in them, one in each. This is due to the fact that to move from the center of one static region to another, the shortest path in~\eqref{eq:discr rate} needs to move out of one set, travel in the transition region to the other set, and move to its center, hence maximizing the minimal distance to all other cornerstones. After finding all static regions, the next cornerstone is to be found in the transition region, if all static regions are approximately of the same size---which we assume here. The crucial observation is, that this~$(C+1)$-st cornerstone is half as far from the other cornerstones\footnote{Here we assume that we are in the coarse spatial resolution case, where the semidistances scale linearly and not quadratically, cf.\ Section~\ref{ssec:simpexam2}. Otherwise, the drop in the distance is more than a factor two (towards factor four).}, as they are from one another. In other words,~$d(c_i,c_{C+1}) + d(c_j,c_{C+1}) \approx d(c_i,c_j)$, $i,j\le C$.

To summarize, our simple check when to stop searching for cornerstones is going to be, when the value of the objective function in~\eqref{eq:coresearch} drops by at least a factor two compared with the previous value. Observe how nicely this works in the periodically forced double gyre case: the rightmost points of the curves in Figure~\ref{fig:DG_meetdist_fromcores} are the optimizers, and the corresponding value of the yellow curve is less than half of the values for the first two cornerstones. This indicates to stop with three cornerstones, as they will represent both the gyres and the transition region.

\subsection{Clustering and fuzzy affiliations}

\begin{figure}[h]
\centering
\includegraphics[height = 35mm]{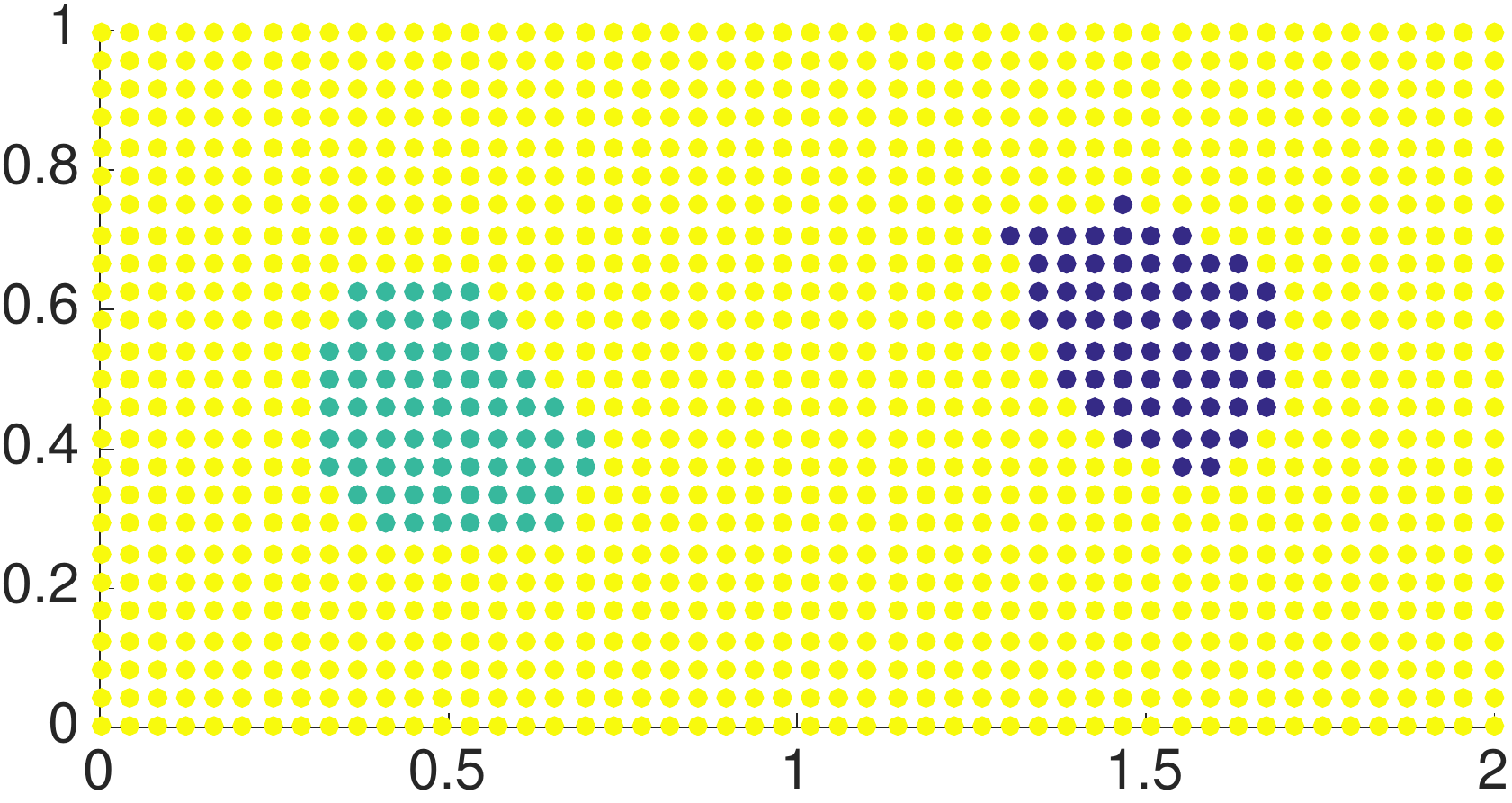}
\qquad
\includegraphics[height = 35mm]{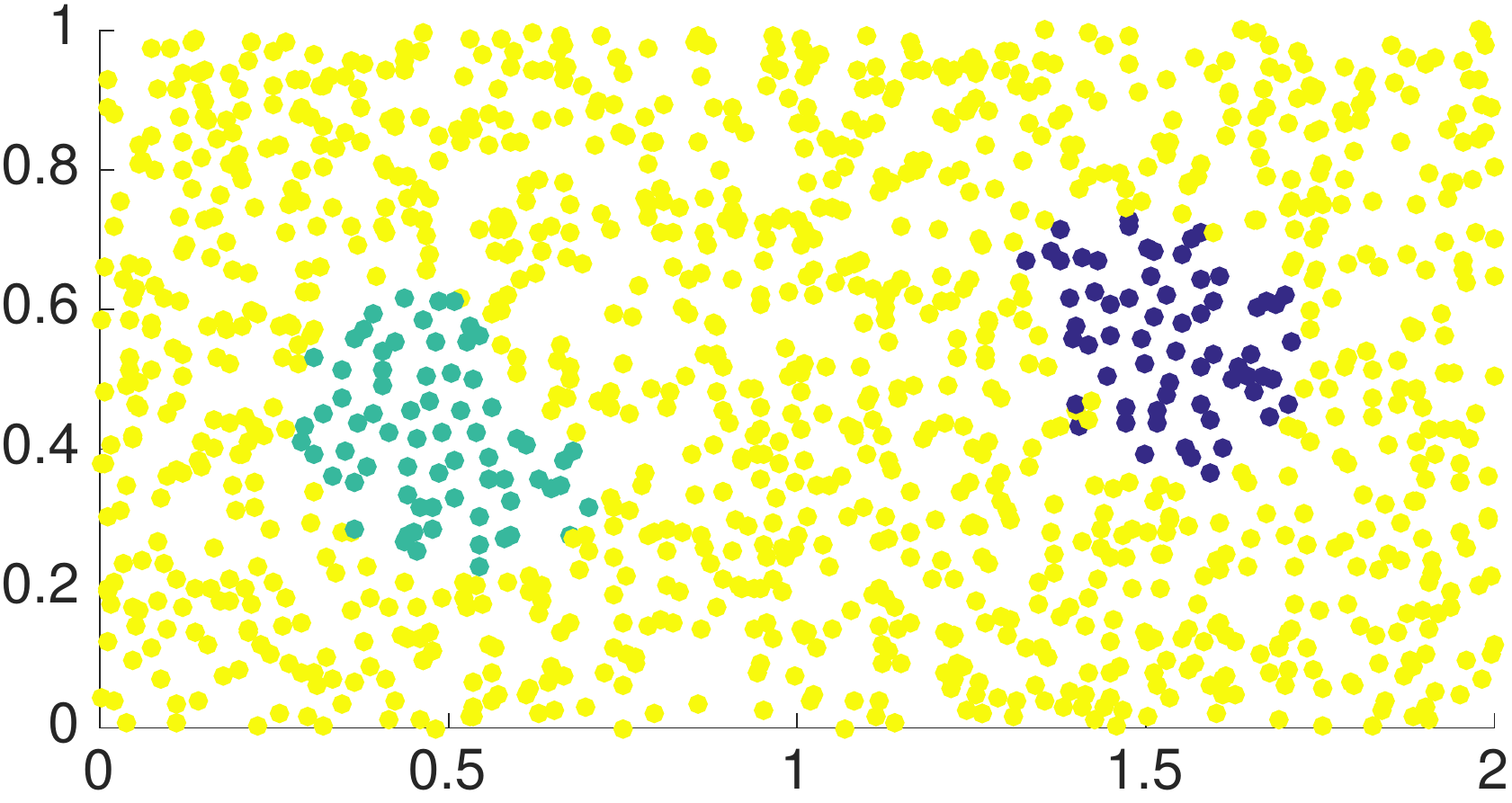}
\caption{The trajectories closest in terms of~$\nu_K^\meet$ to one of the cornerstones than to the others. Left: initial time, right: final time. The horizontal axis is~$x_1$, the vertical is~$x_2$.}
\label{fig:DG_meetdist_nearestcluster_3}
\end{figure}

To get an even more crisp picture of the subdivision of the state space into regions which are far away in terms of the semidistance $d$, we assign to each cornerstone~$c_1,c_2,c_3$ the trajectories that are closer to them than to the two other cornerstones, respectively. For the periodically forced double gyre and the meeting distance this is shown in Figure~\ref{fig:DG_meetdist_nearestcluster_3}.

Comparing this picture with the typical trajectories of the time-1 Poincar\'e map of the system (again, see~\cite[Figure~1]{FrPa14}), it appears that the gyre regions in our figure are smaller. This is due to the nature of the transport distance at hand: the gyres are partly made up of so-called ``regular regions'' of the Poincar\'e map, meaning that typical trajectories move on periodic orbits that are approximately concentric circular lines. Transport between these trajectories is only possible through diffusion, and the price one has to pay for this transport in radial direction is reflected by the rate function (recall, this is what we model by the static regions in Section~\ref{subsec:two mixing subdomains}). The cost to get from the center of the gyre (the cornerstone~$c_1$ or~$c_2$) to a regular trajectory in the same gyre is proportional to the ``radial distance'' between them (compare with the static part of the second example in Section~\ref{subsec:two mixing subdomains}). This behavior is not characteristic for the well-mixed transition region, because there the dynamics (eventually) brings any two trajectories close to each other. The effect is most prominent if the time frame of consideration grows infinitely large, and on our finite time horizon it appears as a flattening of the curve. This brings us back to why the blue and green regions in Figure~\ref{fig:DG_meetdist_nearestcluster_3} are smaller than gyres in the Poincar\'e map. The answer is simply, because the outer periodic orbits are closer to the transition region than to the center of the gyre, hence also closer to the cornerstone~$c_3$ that is in the transition region, because the points in the transition region have very small distance from one another.

%This is underlined in Figure~\ref{fig:DG_meetdist_loglog}. If we were to take a regular grid of points, and plot their Euclidean distances from a fixed point~$x$ ranked in ascending order, the distances would follow an~$r^{1/2}$ asymptotics, where $r$ denotes the rank. This is because there are~$\mathcal{O}(\delta)$ points $\delta$-far from $x$, or, equivalently, there are $\mathcal{O}(\delta^2)$ points closer to $x$ than any point that has distance~$\delta$ to~$x$, and thus~$\delta^2 \sim r$. This shows up as a slope $1/2$ on the double-logarithmic plot Figure~\ref{fig:DG_meetdist_loglog} (right) for the close-by points, hence the distance around the cornerstone $c_1$ grows proportionally to the radius, indicating concentric periodic trajectories. The same holds for~$c_2$ (now shown). This behavior is not characteristic for the well-mixed transition region, because the dynamics (eventually) brings any two trajectories close to each other. This effect is most prominent if the time frame of consideration grows infinitely large, and on our finite time horizon it appears as a flattening of the curve. This brings us back to why the blue and green regions in Figure~\ref{fig:DG_meetdist_nearestcluster_3} are smaller than gyres in the Poincar\'e map. The answer is simply, because the outer periodic orbits are closer to the transition region than to the center of the gyre, hence also closer to the cornerstone $c_3$ that is in the transition region, because the points in the transition region have very small distance from one another.

\begin{figure}[h]
\centering
\includegraphics[width = 0.32\textwidth]{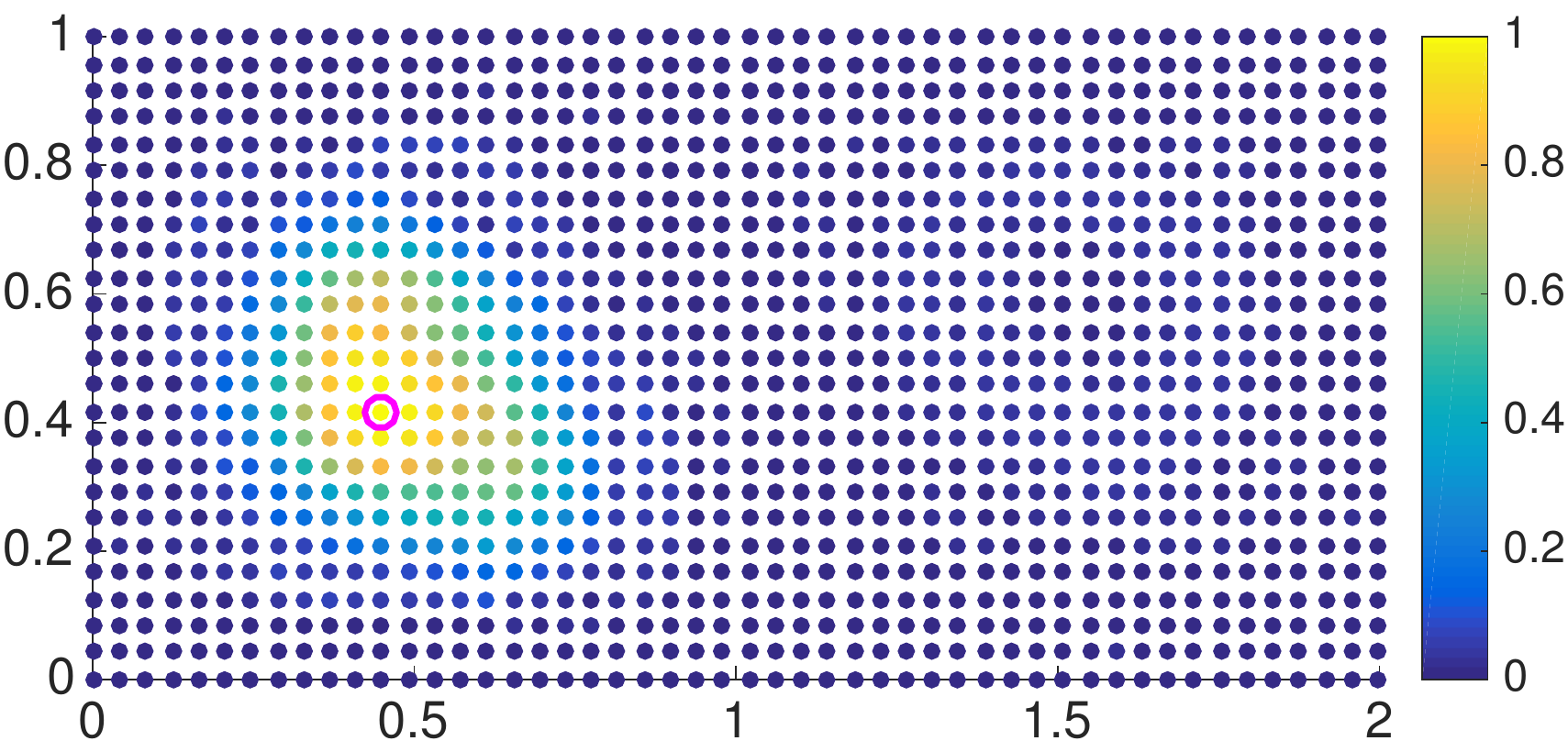}
\hfill
\includegraphics[width = 0.32\textwidth]{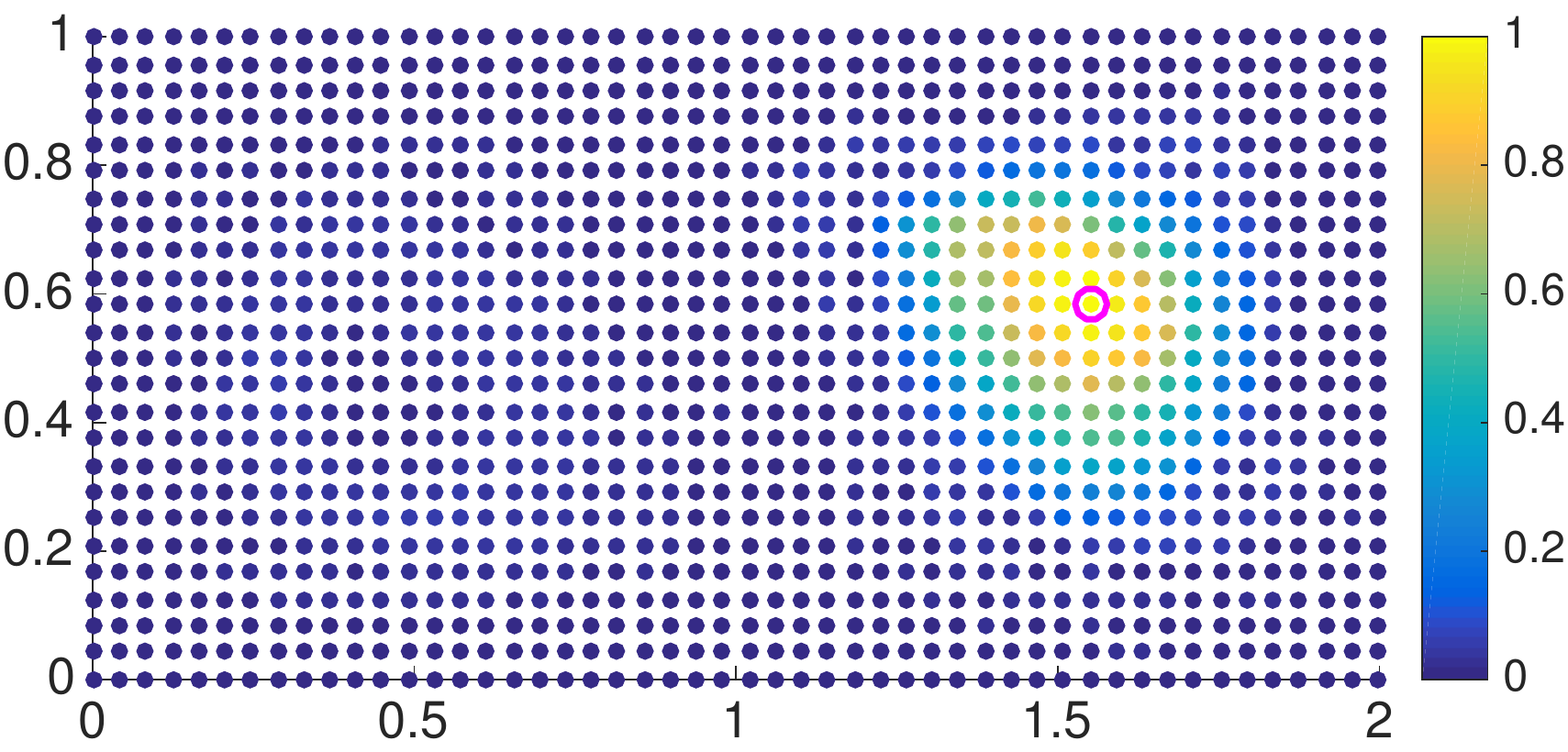}
\hfill
\includegraphics[width = 0.32\textwidth]{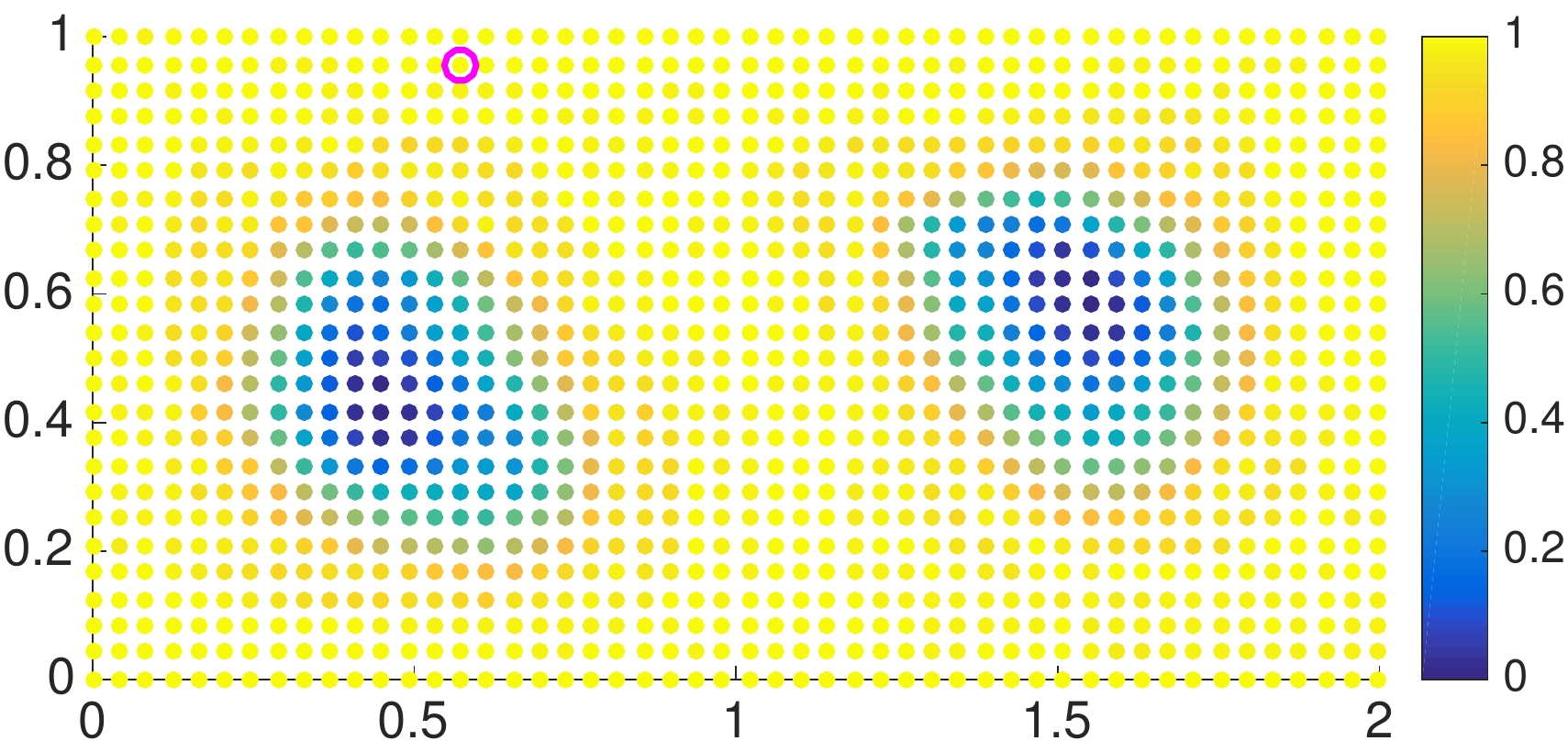}
\hfill
\caption{Fuzzy affiliations~$q_{c_i}(\cdot)$ of the trajectories to the three cornerstones, $c_1,c_2,c_3$ (from left to right) for fuzziness exponent~$m=2$, shown at initial time. The horizontal axis is~$x_1$, the vertical is~$x_2$.}
\label{fig:DG_meetdist_fuzzy}
\end{figure}

Instead of a hard clustering we can assign the trajectories to the cornerstones by \emph{fuzzy affiliations}~$q_{c_i}(\cdot)$, to obtain more refined information on coherence. For instance, let~$m>1$, and minimizing the affiliation-weighted penalty function
\[
\sum_{j=1}^I\sum_{i=1}^{\ell} q_{c_i}(j)^m d(c_i,j)^2
\]
subject to the constraints~$0\le q_{c_i}$ for~$i=1,\ldots,\ell$ and~$\sum_{i=1}^{\ell}q_{c_i}(j)=1$ for every~$j=1,\ldots,I$, yields
\begin{equation}
q_{c_i}(j) = \frac{1}{\sum_{k=1}^{\ell}\left(\frac{d(c_i,j)}{d(c_k,j)}\right)^{\frac{2}{m-1}}}\,.
\label{eq:fuzzyaffiliation}
\end{equation}
This is the affiliation function in the \emph{fuzzy c-means algorithm}~\cite{Bez81}, giving~$q_{c_i}(j)=1$ $\Leftrightarrow$ $d(c_i,j)=0$, i.e., affiliation is maximal if the distance is minimal. Further, the parameter~$m$ controls the fuzziness of the clustering: large~$m$ gives soft clusters, while~$m$ approaching~$1$ gives more and more ``crisp'' clusters as the affiliations converge either to~$0$ or to~$1$~\cite{BezEtAl87}. The resulting affiliations (indicated at initial time) for~$m=2$ are shown in Figure~\ref{fig:DG_meetdist_fuzzy}. For~$m$ close to~$1$ we obtain affiliations very similar to the hard clusters in Figure~\ref{fig:DG_meetdist_nearestcluster_3}.

%\begin{figure}[h!]
%\centering
%\includegraphics[width = 0.49\textwidth]{meetdist_maxdist_reldiff}
%\hfill
%\includegraphics[width = 0.49\textwidth]{meetdist_maxdist_reldiff_tfinal}
%\caption{Relative difference of distances $\nu_K^\meet(c_1,\cdot)$ and~$\nu_K(c_1\smallto\cdot)$ at initial (left) and final times (right). \red{We never refer to this figure! What did you want to say with this?}}
%\label{fig:DG_meetdist_maxdist}
%\end{figure}

\section{Numerical results}
\label{sec:numerics}

In the previous section we already presented numerical results for the double gyre system, where we used the results to motivate and develop the analysis tools. In this section we apply these tools to two other well-analyzed test cases: the perturbed Bickley Jet and the rotating (transitory) double gyre. They are different paradigmatic examples, as the Bickley Jet has a non-vortex coherent set (the jet core), and the transitory double gyre is not a periodically forced system, thus genuinely living on a finite time interval.

Let us also point out that the examples presented here and in the previous section are all one- or two-dimensional. Although we expect the analysis in higher dimensions to be at least qualitatively not very different from the two-dimensional case, dealing with the nevertheless arising subtle differences (see Remark~\ref{rem:3d case}) is beyond the scope of this conceptual work.

It turns out that the choice between the two semidistances~$\nu_K^\cross$ or~$\nu_K^\meet$ has only marginal influence on the results. In this section we shall mostly work with the cross semidistance.

\subsection{The Bickley Jet}
\label{ssec:Bickley}

We consider a perturbed Bickley Jet as described in~\cite{RypEtAl07}. This is an idealized zonal jet approximation in a band around a fixed latitude, assuming incompressibility, on which three traveling Rossby waves are superimposed, see Figure~\ref{fig:Bickley Jet}.
The dynamics is given by~$\dot x_t = v(t,x_t)$ with~$v(t,x)=(-\frac{\partial\Psi}{\partial x_2}, \frac{\partial\Psi}{\partial x_1})$ and stream function
\begin{equation*}
  \Psi(t,x_1,x_2) = -U_0 L \tanh\big(x_2/L\big) + U_0L\,\mathrm{sech}^2\big(x_2/L\big) \sum_{n=1}^3 A_n\cos\left(k_n\left(x_1- c_n t\right)\right).
\end{equation*}
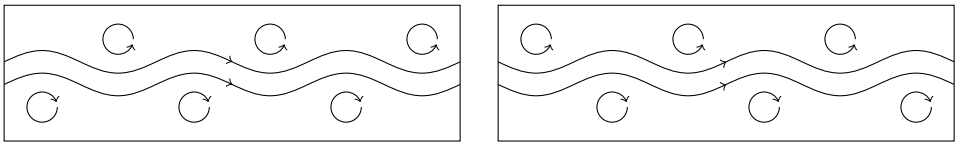
\begin{figure}[h!]
\centering
\begin{tikzpicture}[scale=0.5]
\tikzstyle{every node}=[font=\small]
\begin{scope}[decoration={markings,mark=at position 0.5 with {\arrow{>}}}] 
  \draw(0,0) rectangle (12,3.6);
%  \draw(1,0.9) circle (0.4);
  \draw[->](1.4,0.9) arc(360:20:0.4);
  \draw[->](5.4,0.9) arc(360:20:0.4);
  \draw[->](9.4,0.9) arc(360:20:0.4);
%  \draw(3,2.7) circle (0.4);
  \draw[->](3.4,2.7) arc(0:340:0.4);
  \draw[->](7.4,2.7) arc(0:340:0.4);
  \draw[->](11.4,2.7) arc(0:340:0.4);
  \draw[postaction={decorate}](0,2.1) sin (1,2.4) cos (2,2.1) sin (3,1.8) cos (4,2.1) sin (5,2.4) cos (6,2.1) sin (7,1.8) cos (8,2.1) sin (9,2.4) cos (10,2.1) sin (11,1.8) cos (12,2.1);
  \draw[postaction={decorate}](0,1.5) sin (1,1.8) cos (2,1.5) sin (3,1.2) cos (4,1.5) sin (5,1.8) cos (6,1.5) sin (7,1.2) cos (8,1.5) sin (9,1.8) cos (10,1.5) sin (11,1.2) cos (12,1.5);

  \draw(13,0) rectangle (25,3.6);
%  \draw(1,0.9) circle (0.4);
  \draw[->](16.4,0.9) arc(360:20:0.4);
  \draw[->](20.4,0.9) arc(360:20:0.4);
  \draw[->](24.4,0.9) arc(360:20:0.4);
%  \draw(3,2.7) circle (0.4);
  \draw[->](14.4,2.7) arc(0:340:0.4);
  \draw[->](18.4,2.7) arc(0:340:0.4);
  \draw[->](22.4,2.7) arc(0:340:0.4);
  \draw[postaction={decorate}](13,2.1) sin (14,1.8) cos (15,2.1) sin (16,2.4) cos (17,2.1) sin (18,1.8) cos (19,2.1) sin (20,2.4) cos (21,2.1) sin (22,1.8) cos (23,2.1) sin (24,2.4) cos (25,2.1);
  \draw[postaction={decorate}](13,1.5) sin (14,1.2) cos (15,1.5) sin (16,1.8) cos (17,1.5) sin (18,1.2) cos (19,1.5) sin (20,1.8) cos (21,1.5) sin (22,1.2) cos (23,1.5) sin (24,1.8) cos (25,1.5);
\end{scope}
\end{tikzpicture}
\caption{Sketch of the Bickley Jet flow field at two different times. The flow pattern travels from left to right on the horizontally periodic domain. The horizontal axis is~$x_1$, the vertical is~$x_2$.}
\label{fig:Bickley Jet}
\end{figure}
The constants are chosen as in~\cite[Section~4]{RypEtAl07}. In particular, we set~$k_n = 2n/r_e$ with~$r_e = 6.371$, $U_0 = 5.414$, and~$L = 1.77$. The phase speeds~$c_n$ of the Rossby waves are~$c_1 = 0.1446U_0$, $c_2 = 0.205U_0$, $c_3 = 0.461U_0$, their amplitudes~$A_1 = 0.0075$, $A_2 = 0.15$, and~$A_3 = 0.3$, as in~\cite{HaEtAl16}. The system is considered on a state space~$X = [0,\pi r_e]\times[-3,3]$ which is periodic in the horizontal~$x_1$ coordinate.
 
We choose a uniform $60\times 18$ grid as initial conditions for the floaters at time~$t=0$; i.e.,~$I=1080$. We sample the trajectories of these floaters at times~$t_k = k\tau$, $k=0,1,\ldots,K$, where~$K=80$ and~$\tau = 0.5$. In this time interval, typical trajectories cross the cylindrical state space horizontally~4-5 times, trajectories in the jet core (the wavy structure in Figure~\ref{fig:Bickley Jet}) up to~9 times.

Employing our large-deviation based distance computations on this data set using Algorithm~\ref{algo1} and $\alpha=1/2$, we get the rates~$\nu_K(i\smallto j)$, $i,j=1,\ldots,I$. From these rates we readily obtain the~$\nu_K^\cross(i,j)$ via
\[
\nu_K^\cross(i,j) = \nu_K(i\smallto j) + \nu_K(j\smallto i)\,.
\]

We repeat the cornerstone finding analysis from the previous section. The optimal values of the objective function in the cornerstone finding problem~\eqref{eq:coresearch} are for~8 cornerstones, in order:
\[
2.06,\ 3.21,\ 2.45,\ 2.33,\ 2.30,\ 2.14,\ 1.42,\ 0.70.
\]
Recall, that the first value is with respect to a random cornerstone~$c_0$ that we discard. These numerical values with our previous analysis shed light on the topological structure of the state space with respect transport and mixing. Note, that our assumption from Section~\ref{ssec:num cornerstones}, that all coherent sets are divided by \emph{one} mixing region, is not satisfied: the jet core is a coherent set itself, dividing two mixing regions (below and above it), each containing 3 further coherent sets (the gyres). Thus,~$c_1$ and~$c_2$ have maximal distance~($\nu_K^\cross(c_1,c_2)=3.21$), because the random walker needs to cross the jet core. Every further cornerstone~$c_3,\ldots,c_6$ can be reached from either~$c_1$ or~$c_2$ through one of the mixing regions, and thus have a very similar cost. The deviation of these costs,~$2.14-2.45$, shows that we did not reach the state of full mixing on the chosen time interval.

Now, the seventh cornerstone lies in the jet core, which has to be crossed if traveling between cornerstones that are below and above it, respectively. The corresponding cost~(1.42) is a bit larger than half of the previous cost, because~$c_7$ does not lie on the shortest path between cornerstones below and above the jet core. Intuitively, the ``center line'' of the jet core should be equally far from all cornerstones~$c_1,\ldots,c_6$, if the time interval is large enough such that the regions around the gyres are truly mixing. Since it is not, there are points on the boundary of the jet core which are easier to reach from them, and thus easier to cross there. The cornerstone~$c_7$ represents the position where it is the hardest to cross. The eighth cornerstone has truly half the semidistance to the closest one than~$c_7$, and lies in one of the mixing regions.

\begin{figure}[h]
\centering
\includegraphics[width = 0.48\textwidth]{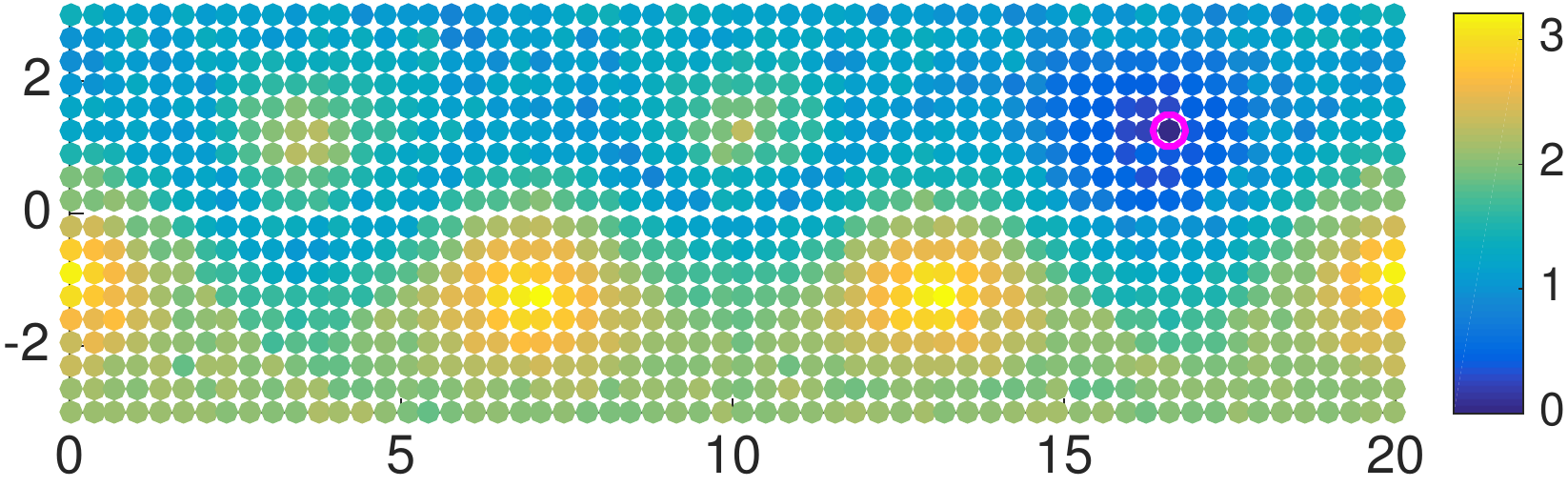}
\hfill
\includegraphics[width = 0.48\textwidth]{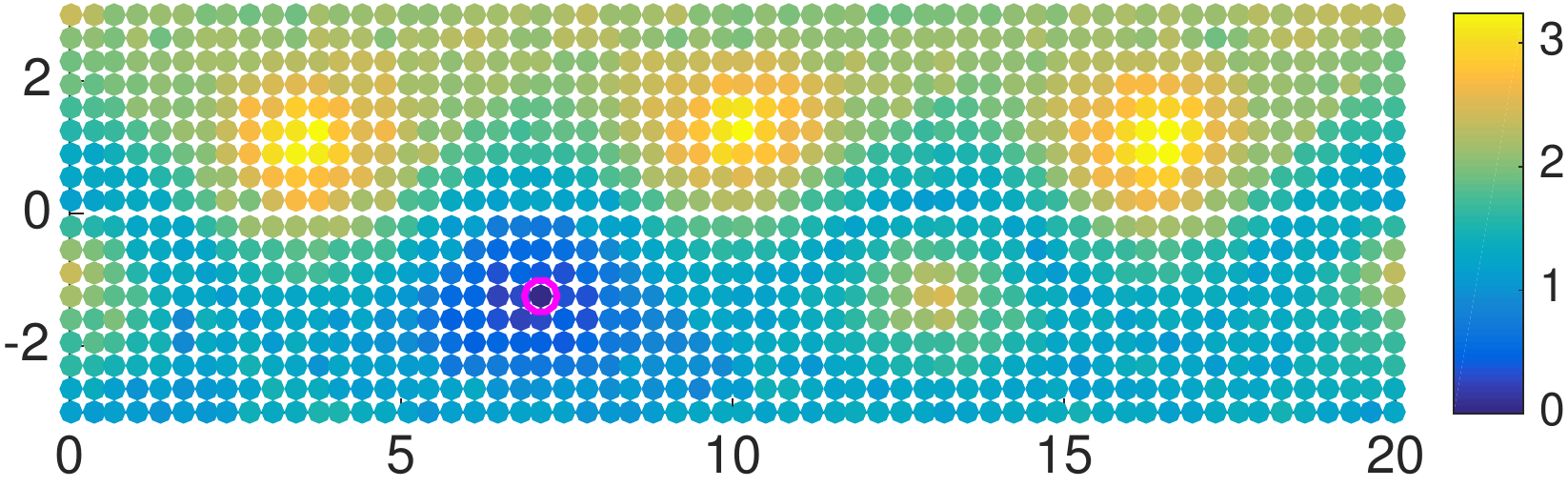}

\includegraphics[width = 0.48\textwidth]{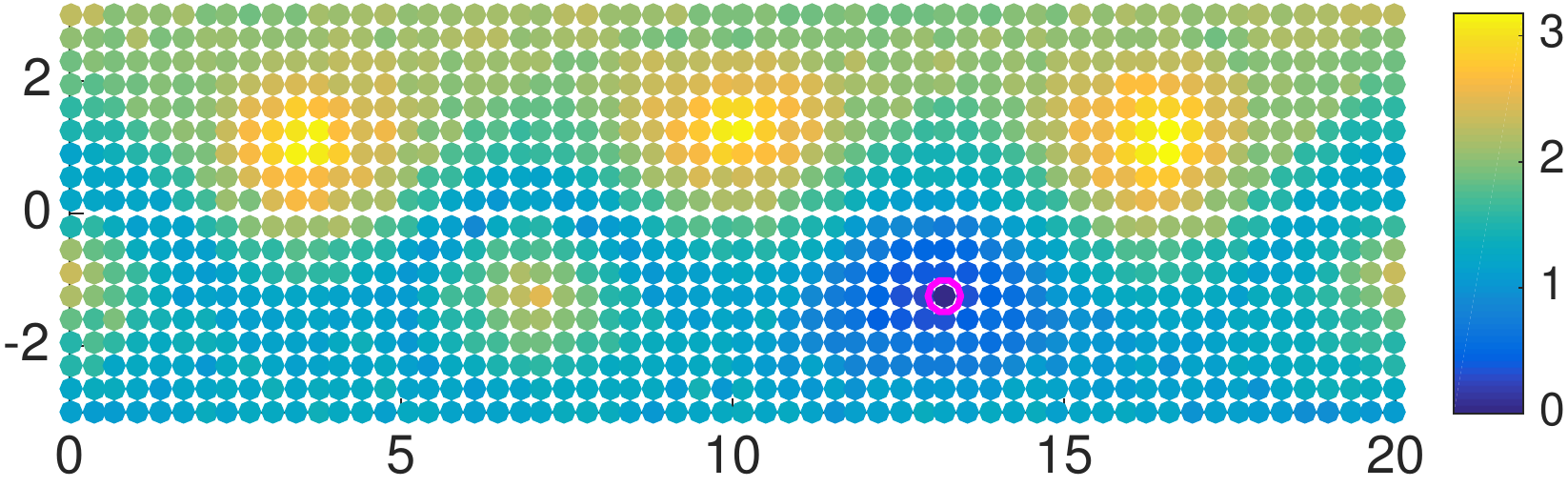}
\hfill
\includegraphics[width = 0.48\textwidth]{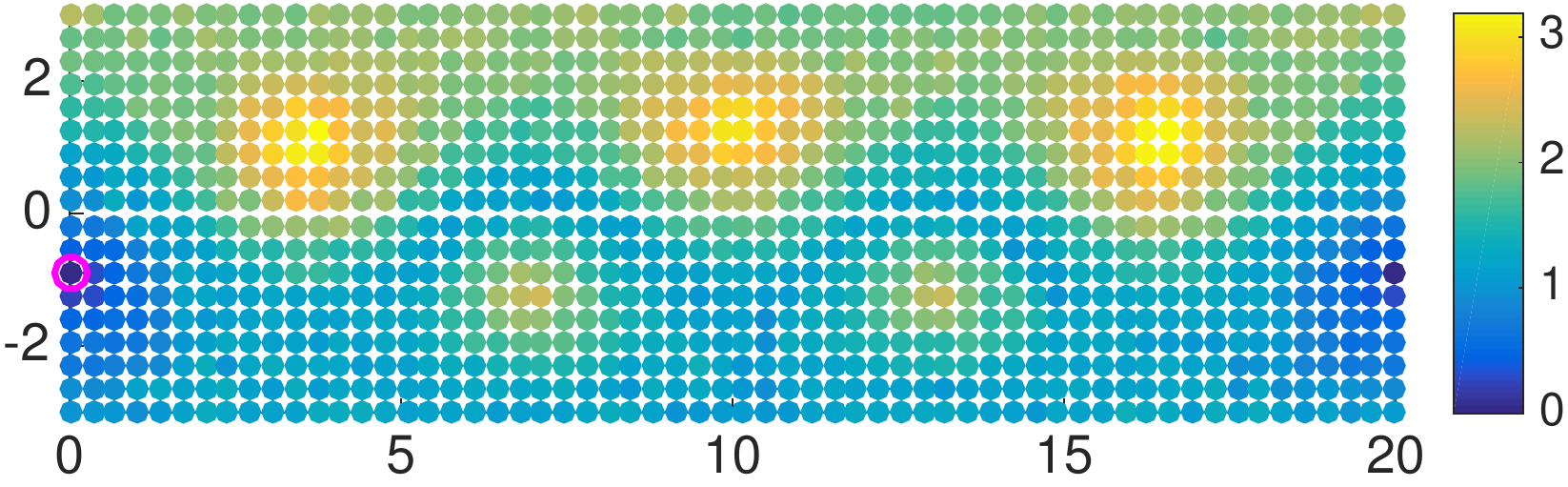}

\includegraphics[width = 0.48\textwidth]{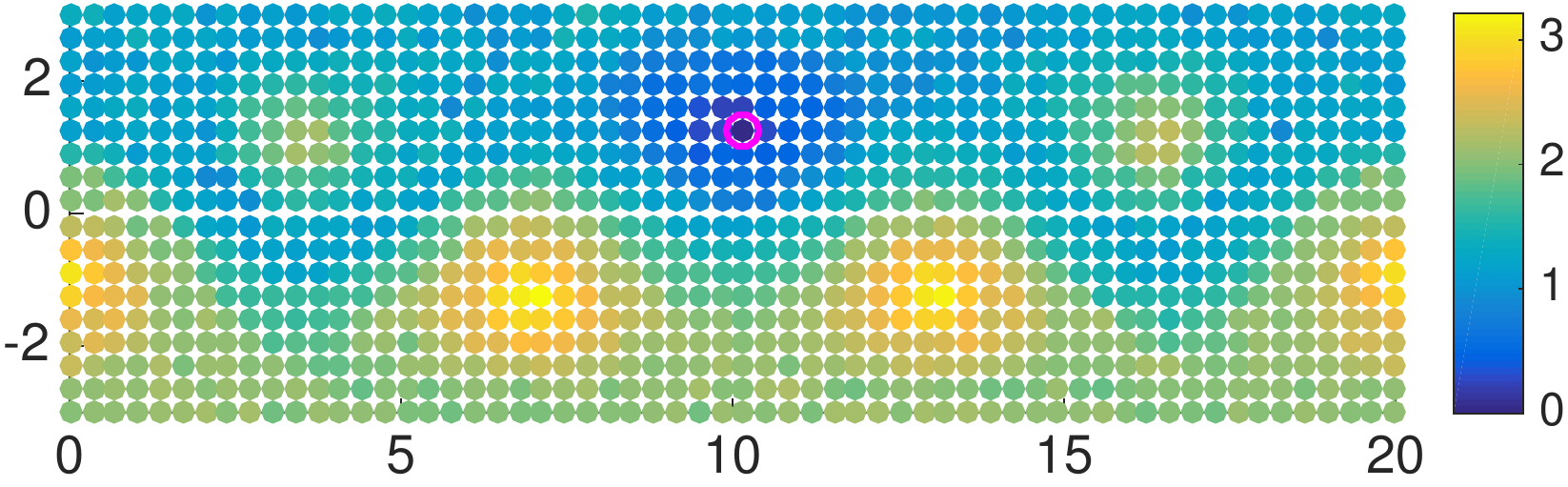}
\hfill
\includegraphics[width = 0.48\textwidth]{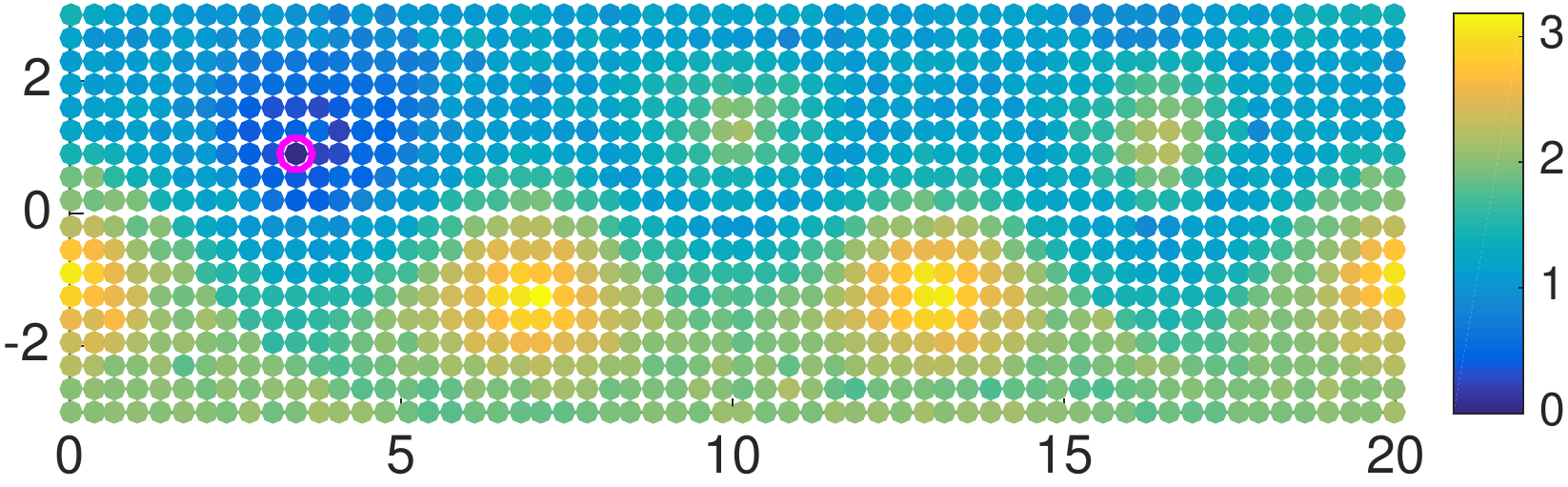}

\includegraphics[width = 0.48\textwidth]{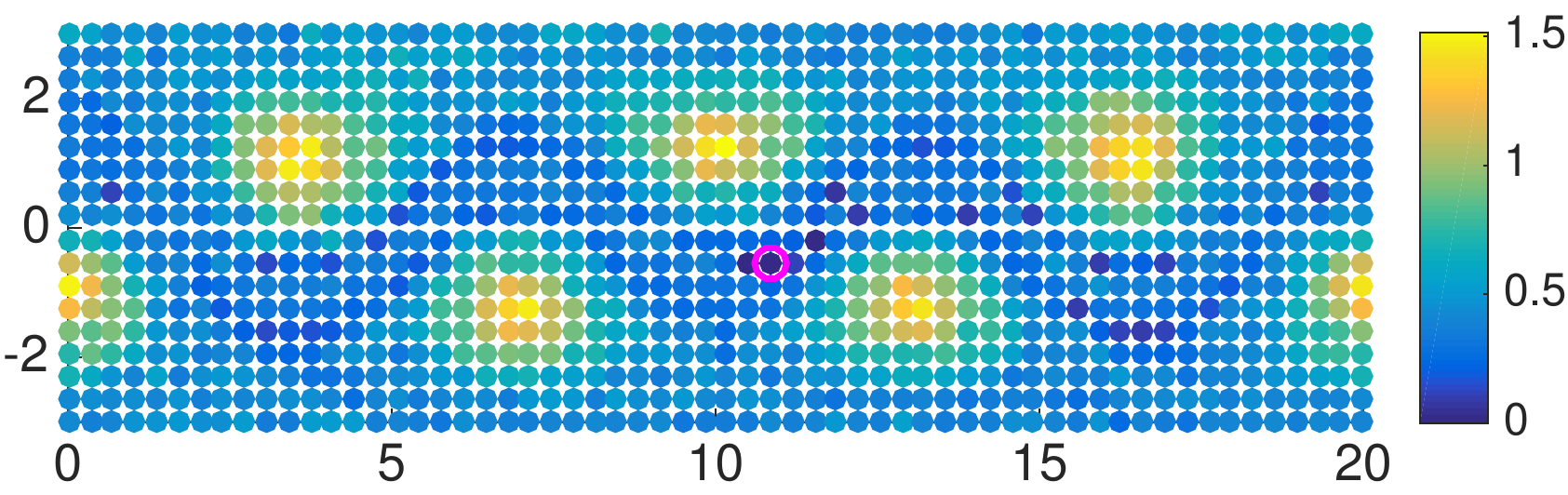}

\caption{Row-wise from top left to bottom: the identified corner stores~$c_i$, $i=1,\ldots,7$, (magenta circles) and their distances~$\nu^{\cross}(c_i,\cdot)$ to the other trajectories, at initial time. The cornerstones are located in the six gyres and the central jet region. The distances are given in units~$1/\tau$. The horizontal axis is~$x_1$, the vertical is~$x_2$.}
\label{fig:Bickley_meetdist_cores}
\end{figure}

\begin{figure}[h]
\centering
\includegraphics[width = 0.48\textwidth]{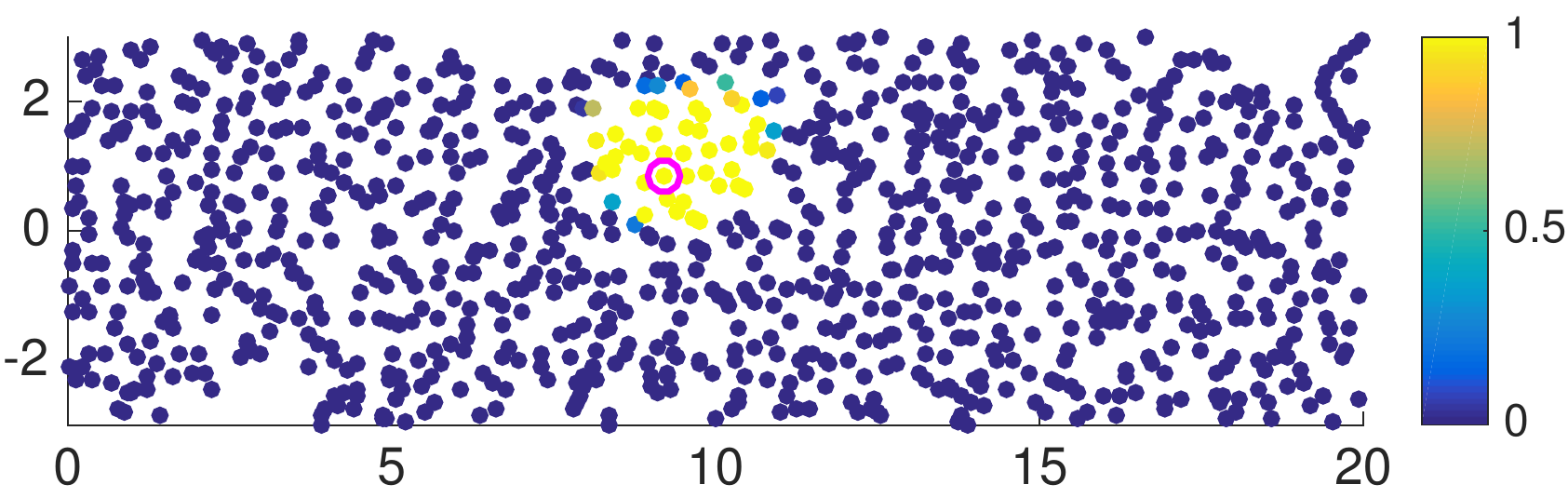}
\hfill
\includegraphics[width = 0.48\textwidth]{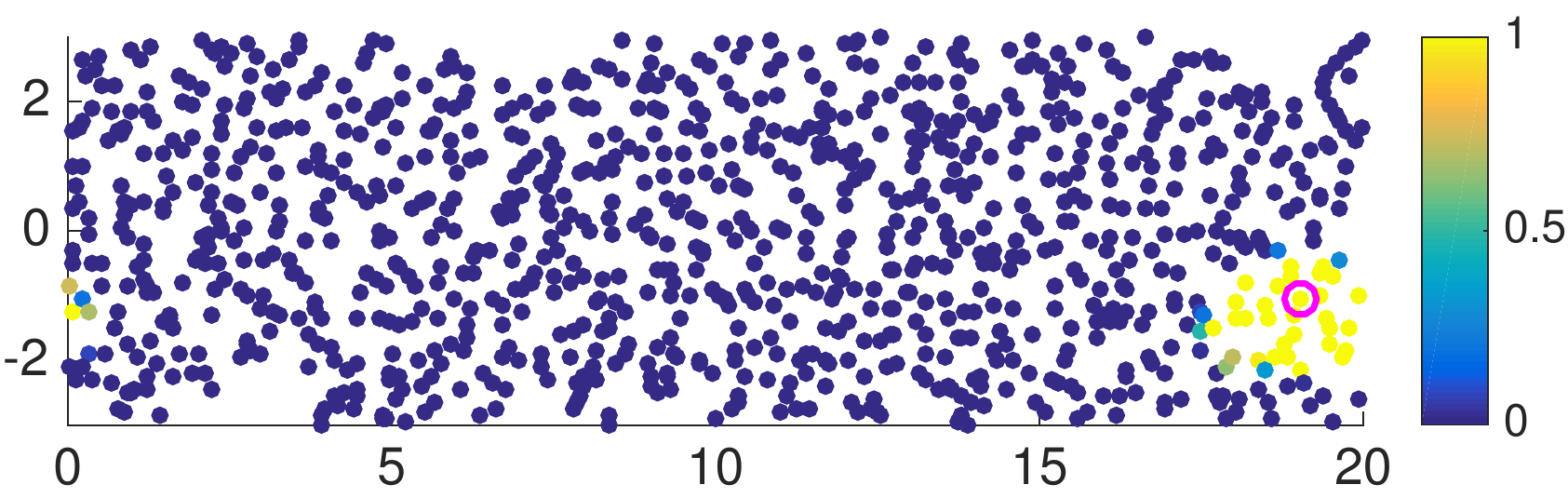}

\includegraphics[width = 0.48\textwidth]{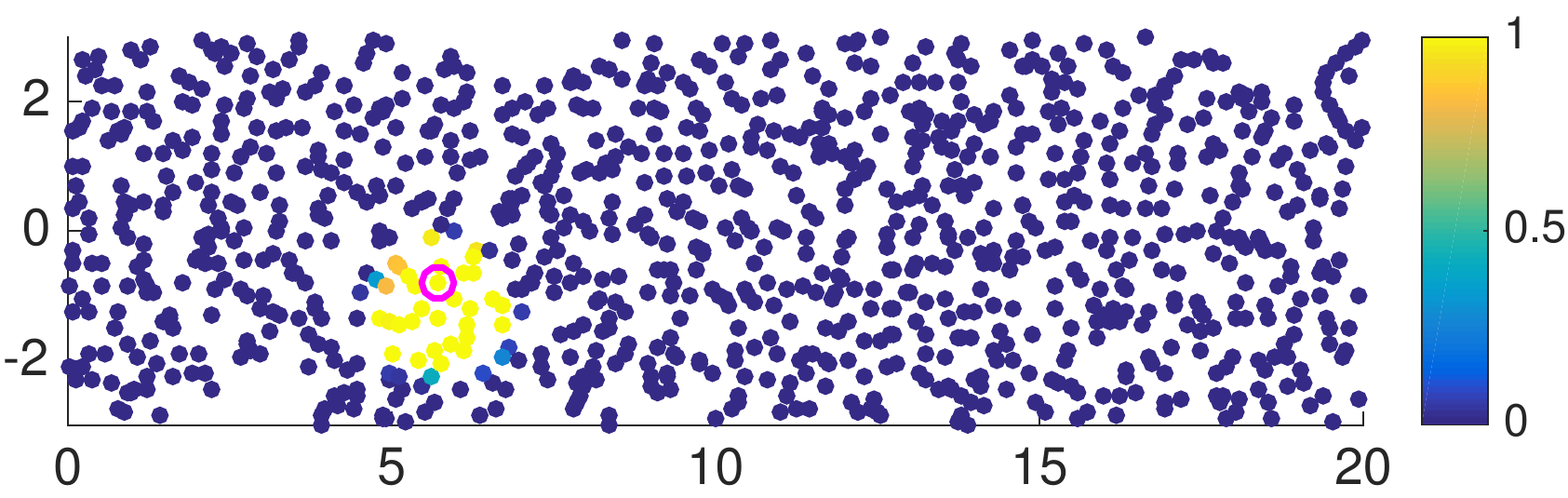}
\hfill
\includegraphics[width = 0.48\textwidth]{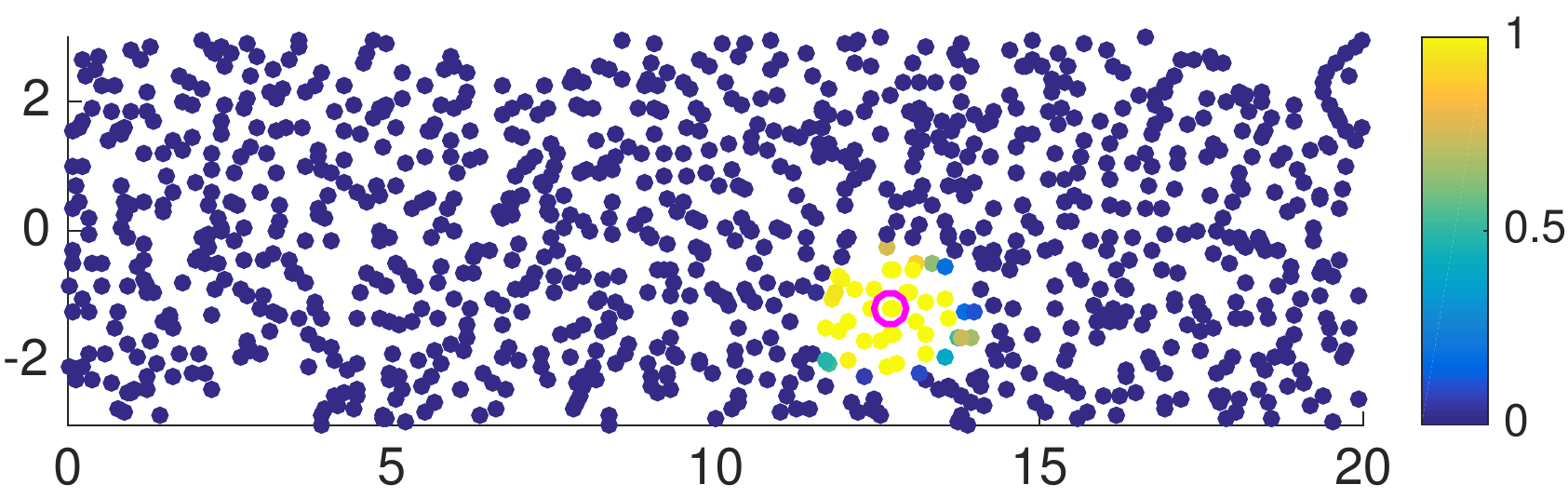}

\includegraphics[width = 0.48\textwidth]{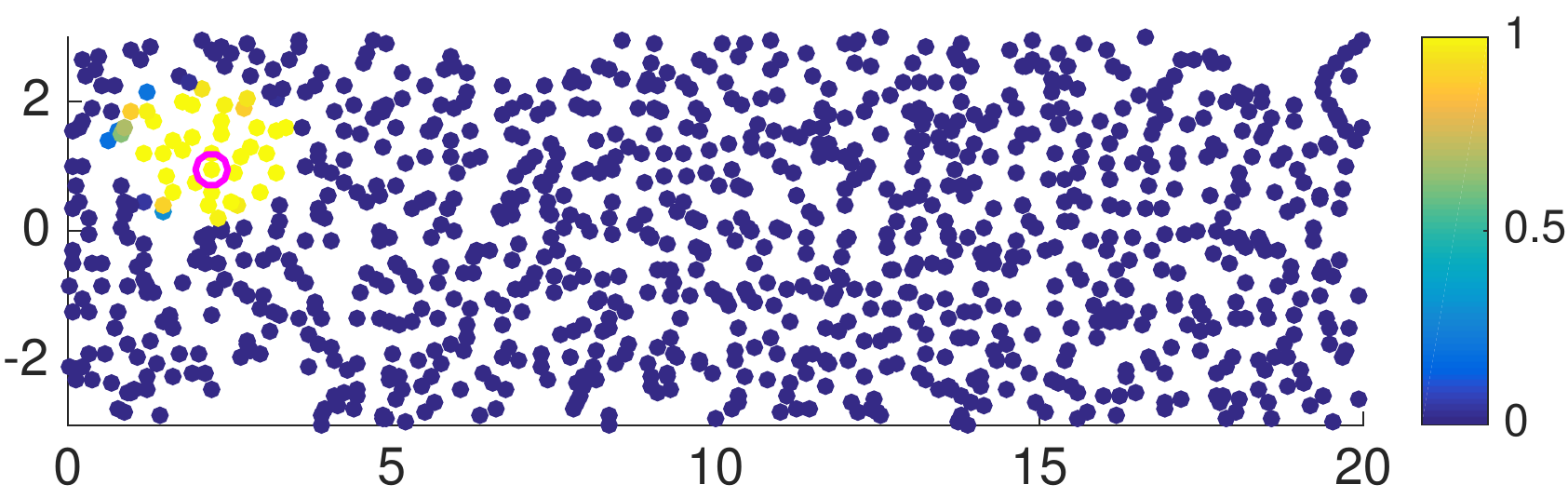}
\hfill
\includegraphics[width = 0.48\textwidth]{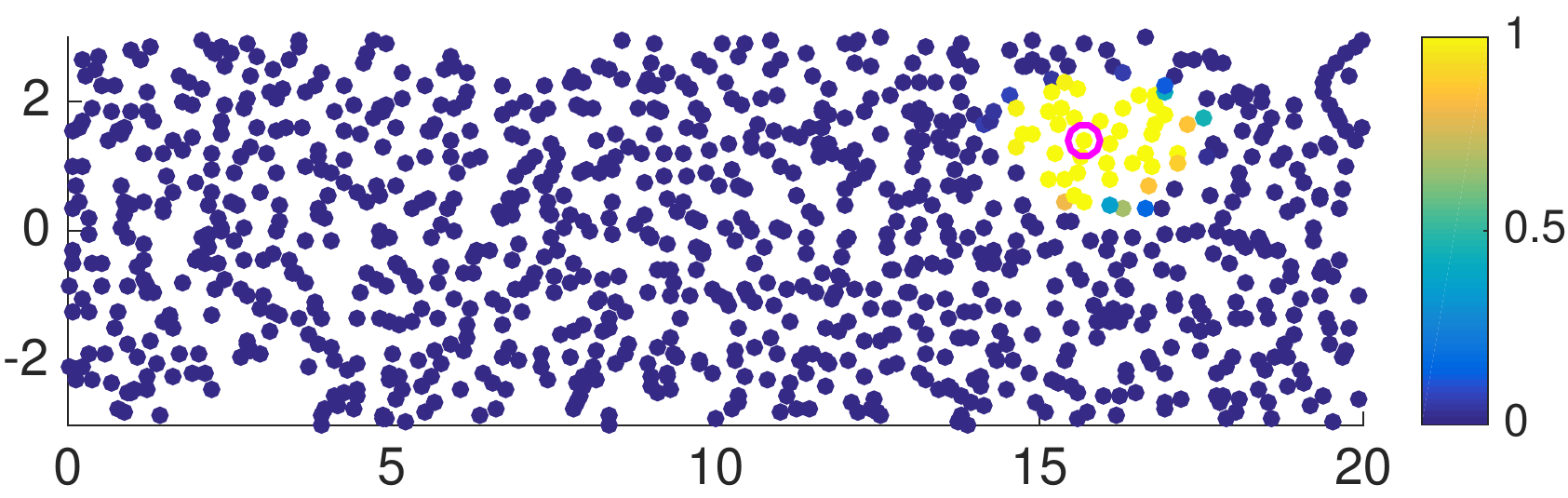}

\includegraphics[width = 0.48\textwidth]{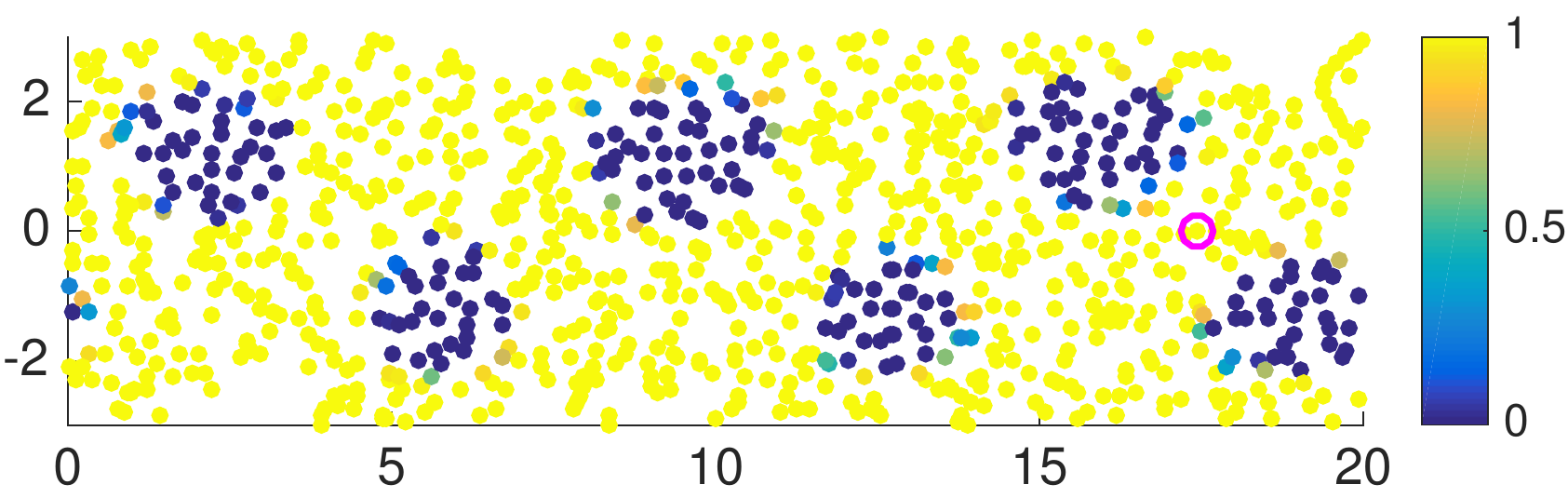}
\hfill
\includegraphics[width = 0.48\textwidth]{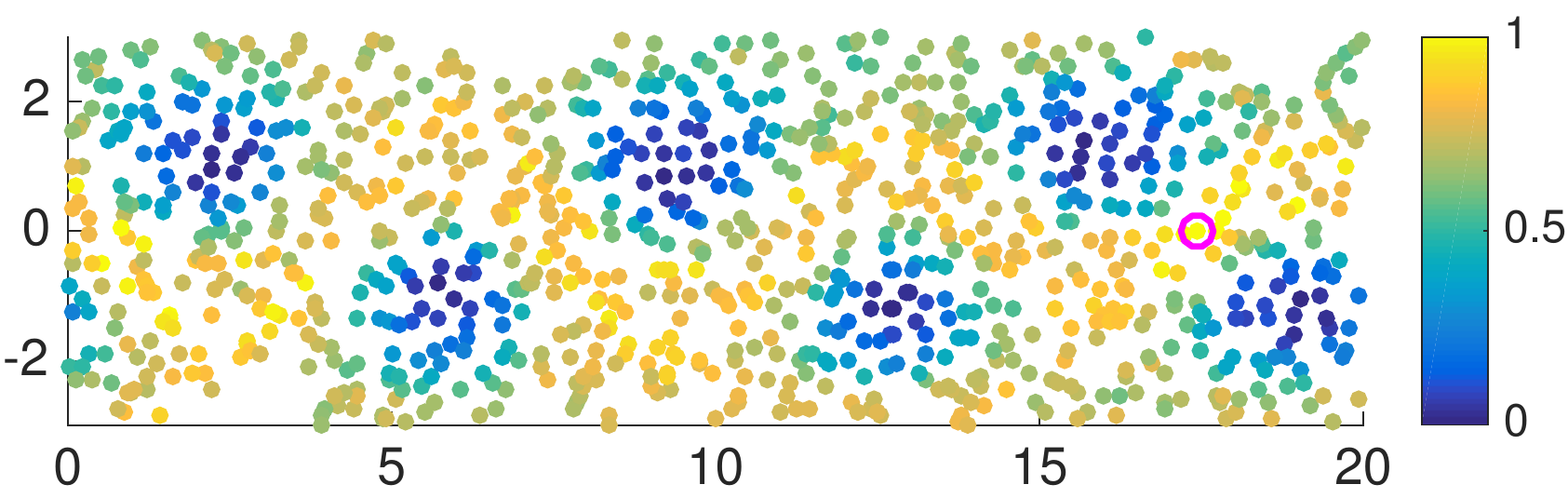}

\caption{Row-wise from top left to bottom: the fuzzy affiliations~\eqref{eq:fuzzyaffiliation} of the trajecories at time~$t=5$ to the cornerstones~$c_1,\ldots,c_7$, respectively (magenta circles). Bottom right: affiliation~$q_{c_7}(\cdot)$ for~$m=1.9$, which suggests that the~$7\text{th}$ coherent region could contain a coherent set itself: the jet core. The horizontal axis is~$x_1$, the vertical is~$x_2$.}
\label{fig:Bickley_meetdist_affiliations}
\end{figure}

We show our results for seven cornerstones\footnote{If we include additional cornerstones, results tend to deteriorate due to the low resolution and because the chosen time interval is not giving full mixing.}. The semidistances are shown in Figure~\ref{fig:Bickley_meetdist_cores}.
The corresponding fuzzy affiliations from~\eqref{eq:fuzzyaffiliation} for~$m=1.1$ are shown in Figure~\ref{fig:Bickley_meetdist_affiliations}. They show a very crisp distinction of the six gyres from the rest of the state space. The bottom right figure shows the affiliation~$q_{c_7}(\cdot)$ for~$m=1.9$, which suggests that the region around the gyres could still be partitioned into coherent sets itself: the jet core appears more strongly affiliated to this cornerstone than the other trajectories. It is not surprising that we could not see this for~$m=1.1$, since the closer~$m$ is to~$1$, the more ``crisp'' the affiliation function is forced to be, and the mixing region is more easily reached from the thin jet core than from the gyres.

\subsection{The rotating double gyre}
\label{ssec:rotating double gyre}

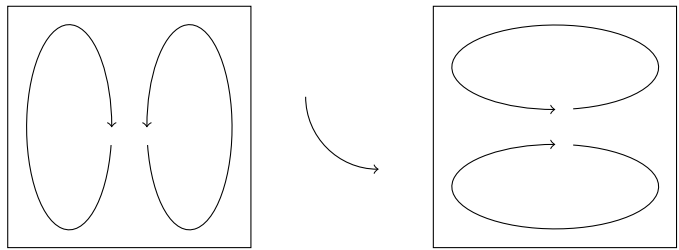
\begin{figure}[h!]
\centering
\begin{tikzpicture}[scale=0.8]
  \draw(0,0) rectangle (4,4);
%  \draw(1,2) ellipse (0.7 and 1.7);
  \draw[->](1.7,1.7) arc (350:0:0.7 and 1.7);
%  \draw(3,2) ellipse (0.7 and 1.7);
  \draw[->](2.3,1.7) arc (-170:180:0.7 and 1.7);
%  \draw[->] (2.4,3.5) sin (2,3.6) cos (1.6,3.5);
%  \draw[->] (1.6,0.5) sin (2,0.4) cos (2.4,0.5);

  \draw(7,0) rectangle (11,4);
%  \draw(7,3) ellipse (1.7 and 0.7);
  \draw[->](9.3,2.3) arc (-80:270:1.7 and 0.7);
%  \draw(7,1) ellipse (1.7 and 0.7);
  \draw[->](9.3,1.7) arc (80:-270:1.7 and 0.7);
%  \draw[->,rotate around={90:(7,2)}] (7.4,3.5) sin (7,3.6) cos (6.6,3.5);
%  \draw[->,rotate around={90:(7,2)}] (6.6,0.5) sin (7,0.4) cos (7.4,0.5);
%  \draw[->](4.5,1.5)--(6.5,1.5);
%  \draw(5.5,2.5) circle (0.8);
  \draw[->](4.9,2.5) arc (180:270:1.2);
\end{tikzpicture}
\caption{Sketch of the flow field of the rotating double gyre at initial (left) and final (right) times. The horizontal axis is~$x_1$, the vertical is~$x_2$.}
\label{fig:rotating double gyre}
\end{figure}

Let us consider a prototype for a system, where transport is considered only on a limited time interval. The rotating double gyre system~\cite{MoMe11} is given by the stream function~$\psi(t,x_1,x_2) = (1-s(t))\psi_P(x_1,x_2) + s(t)\psi_F(x_1,x_2)$, with~$s(t) = t^2(3-2t)$, $\psi_P(x_1,x_2) = \sin(2\pi x_1)\sin(\pi x_2)$, and~$\psi_F(x_1,x_2) = \sin(\pi x_1)\sin(2\pi x_2)$, and is considered on the state space~$X = [0,1]^2$ and time interval~$t\in[0,1]$. The two gyres, which initially occupy the left and right halves of the unit square, turn during this time by~$\pi/2$ to occupy the top and bottom halves at final time, see Figure~\ref{fig:rotating double gyre}.

We choose a uniform~$30\times 30$ grid as initial conditions for the floaters at time~$t=0$; i.e.,~$I=900$. We sample the trajectories of these floaters at times~$t_k = k\tau$, $k=0,1,\ldots,K$, where~$K=100$ and~$\tau = 0.01$. We employ the cross-semidistance, and start our cornerstone search. The first three values of the optimization problem~\eqref{eq:coresearch} are
\[
0.0274,\ 0.0474,\ 0.0262.
\]
We identify the significant drop after two corner stones, hence we expect two coherent sets with one mixing region dividing them. The drop in the distance is by a factor~$0.55$, which is not below one half, the reason for this being again that the time interval of consideration is not sufficient for perfect mixing of the transition region.
\begin{figure}[htb]
\centering
\includegraphics[width = 0.32\textwidth]{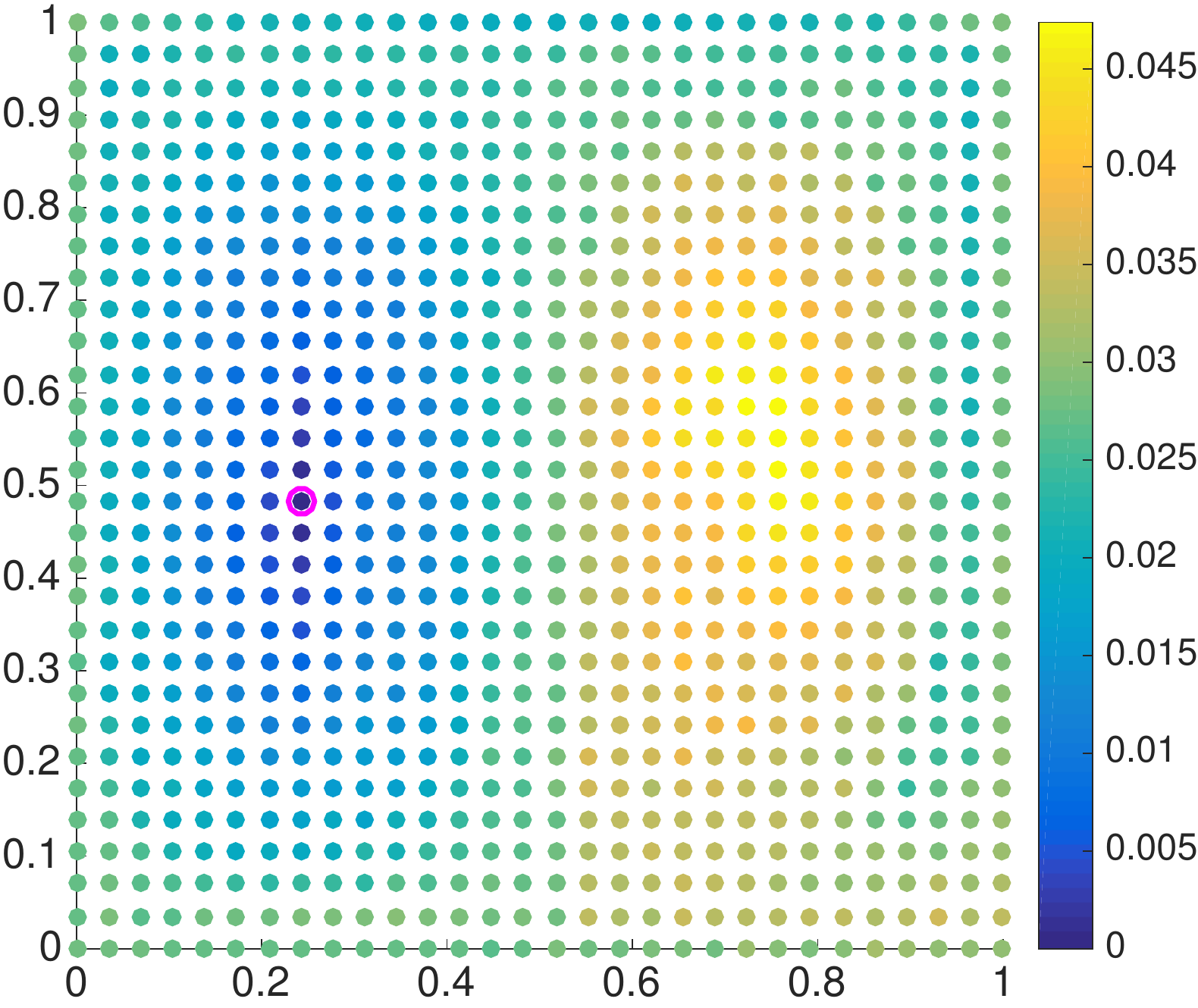}
\hfill
\includegraphics[width = 0.32\textwidth]{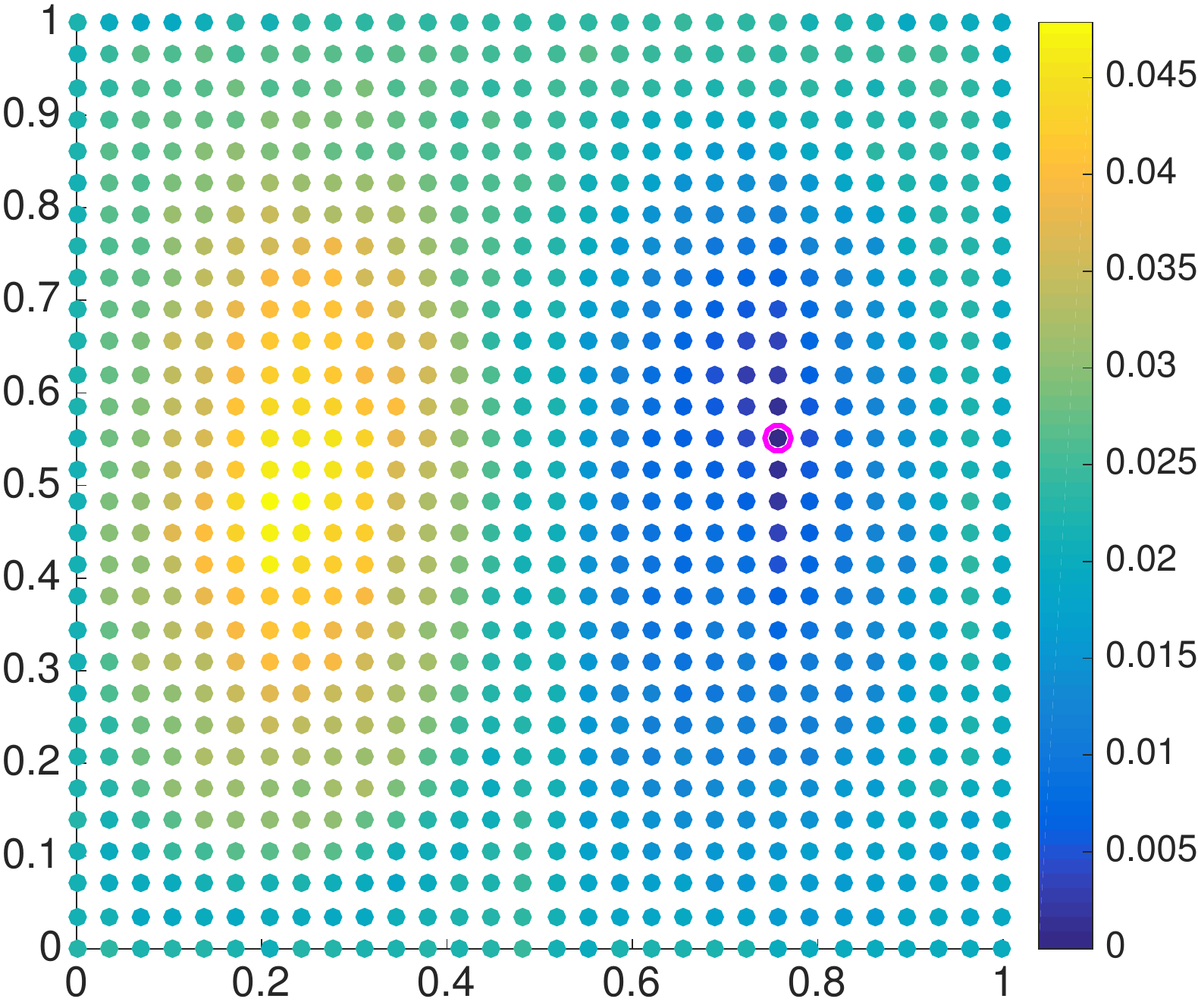}
\hfill
\includegraphics[width = 0.32\textwidth]{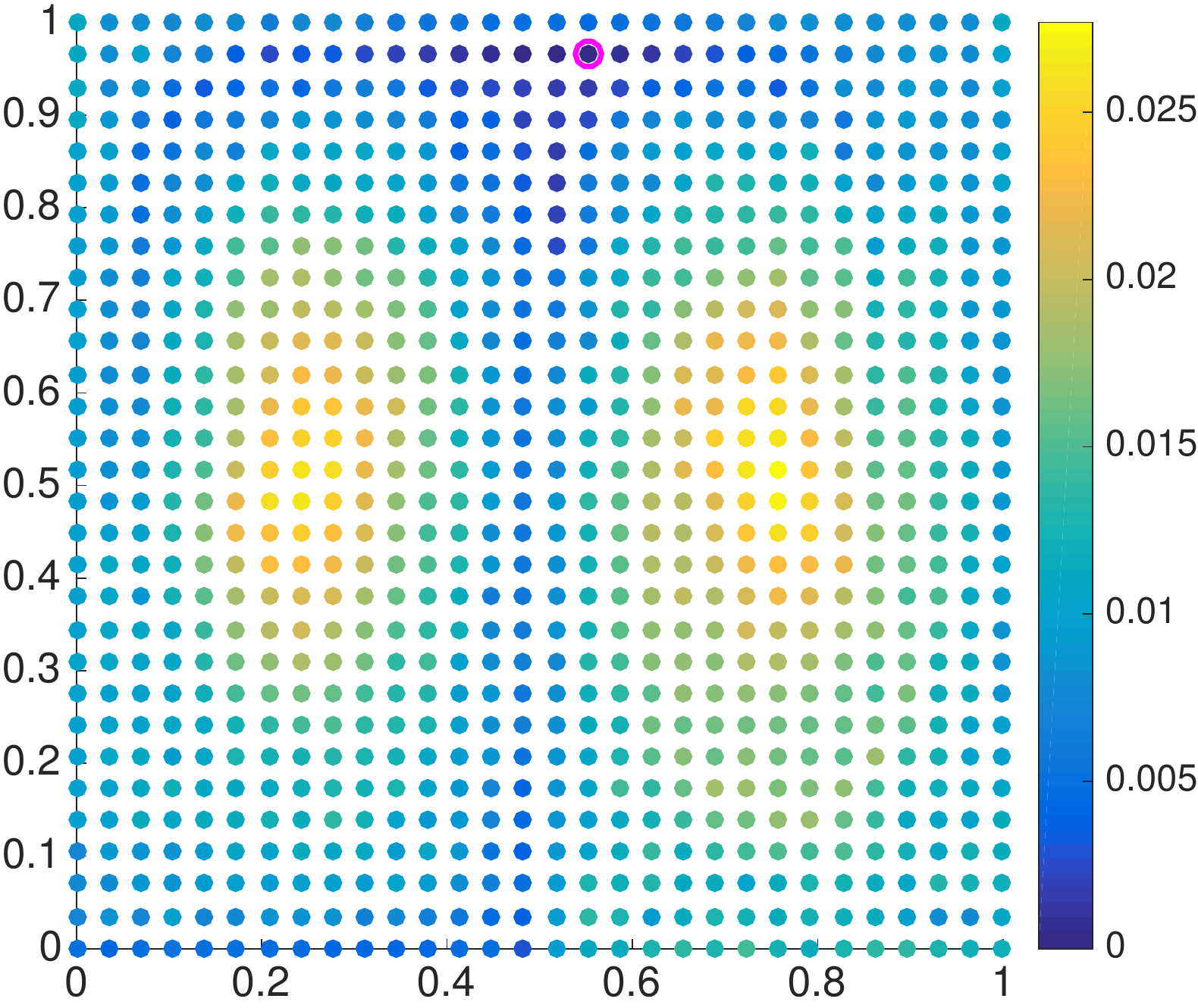}

\includegraphics[width = 0.32\textwidth]{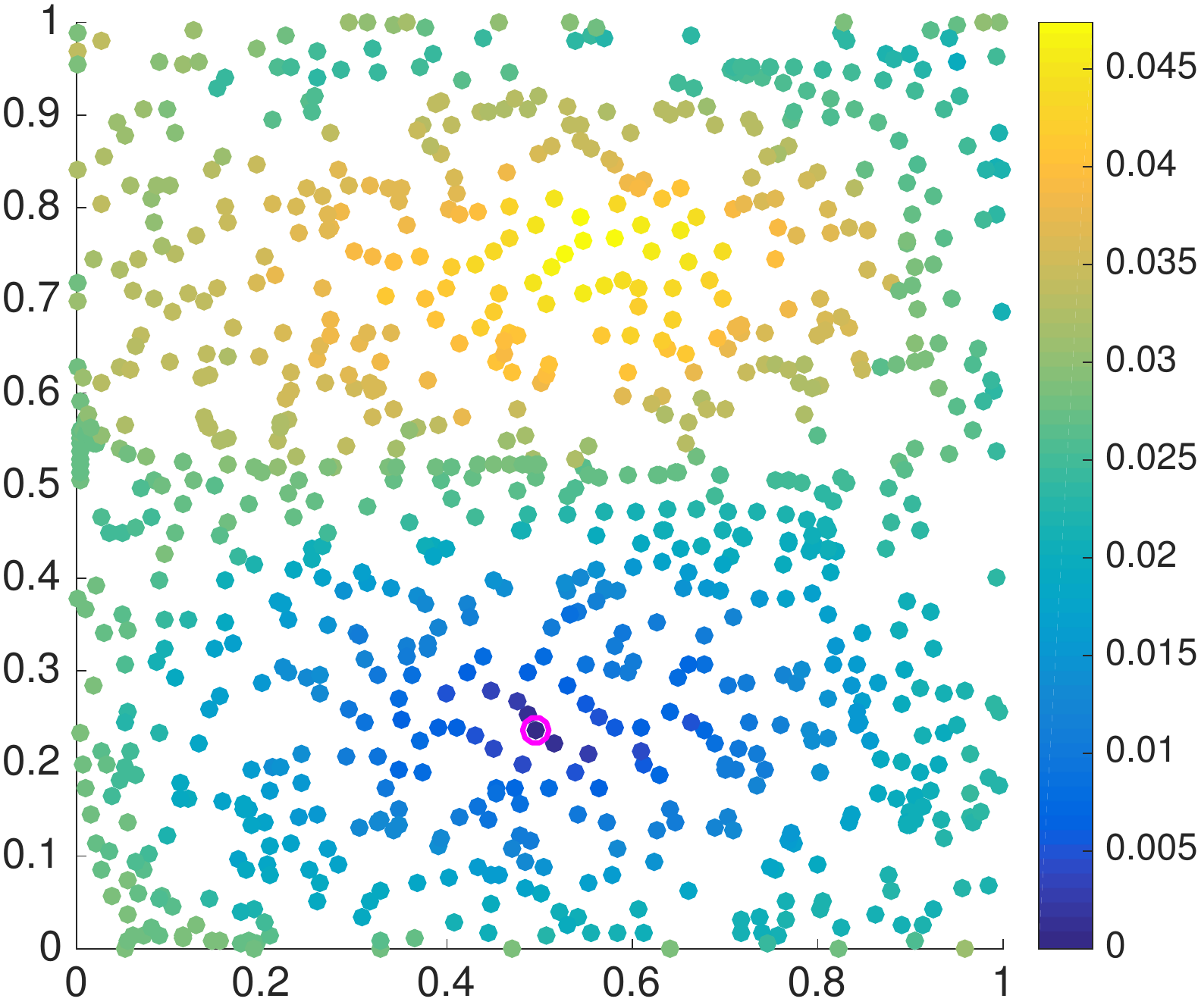}
\hfill
\includegraphics[width = 0.32\textwidth]{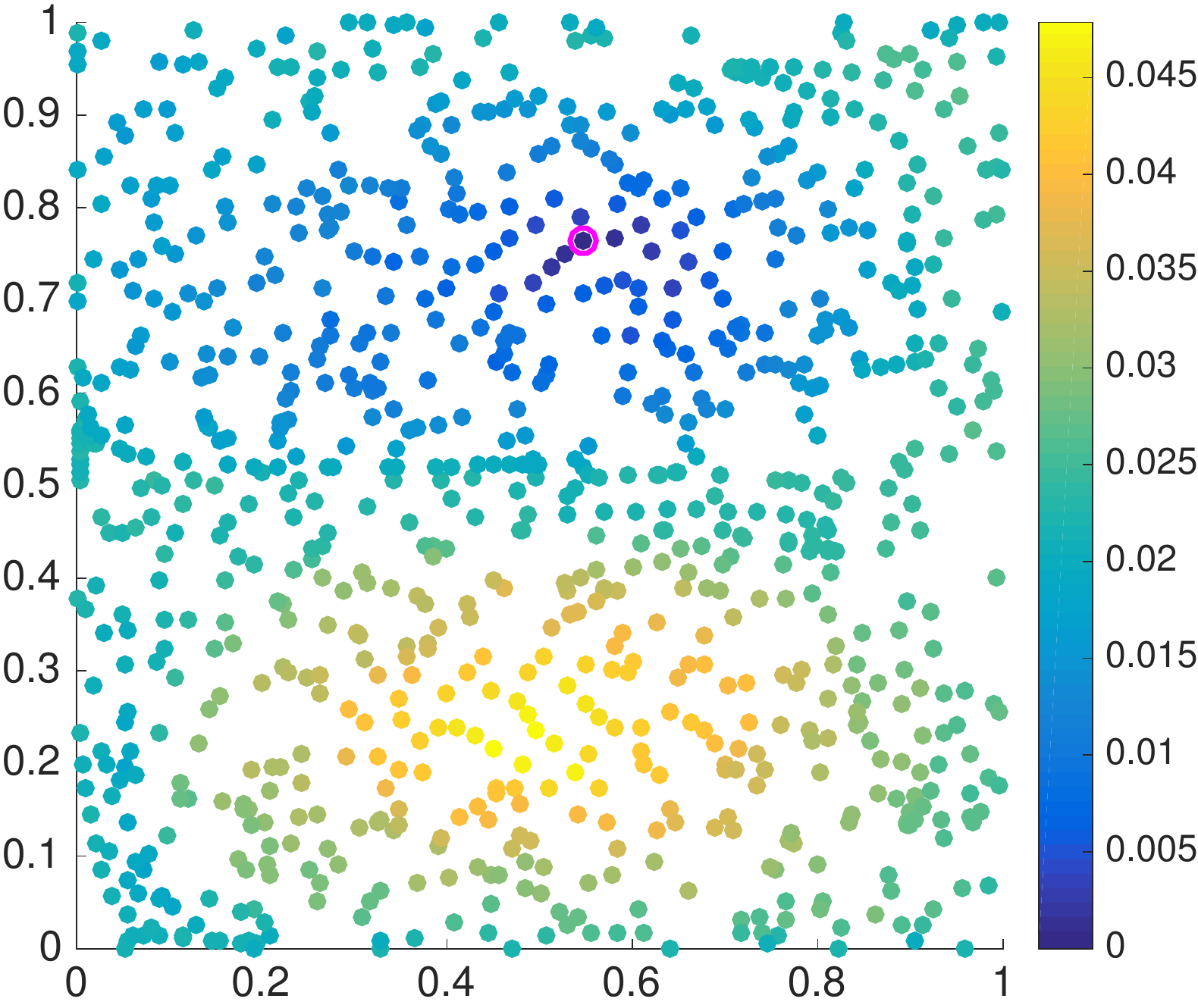}
\hfill
\includegraphics[width = 0.32\textwidth]{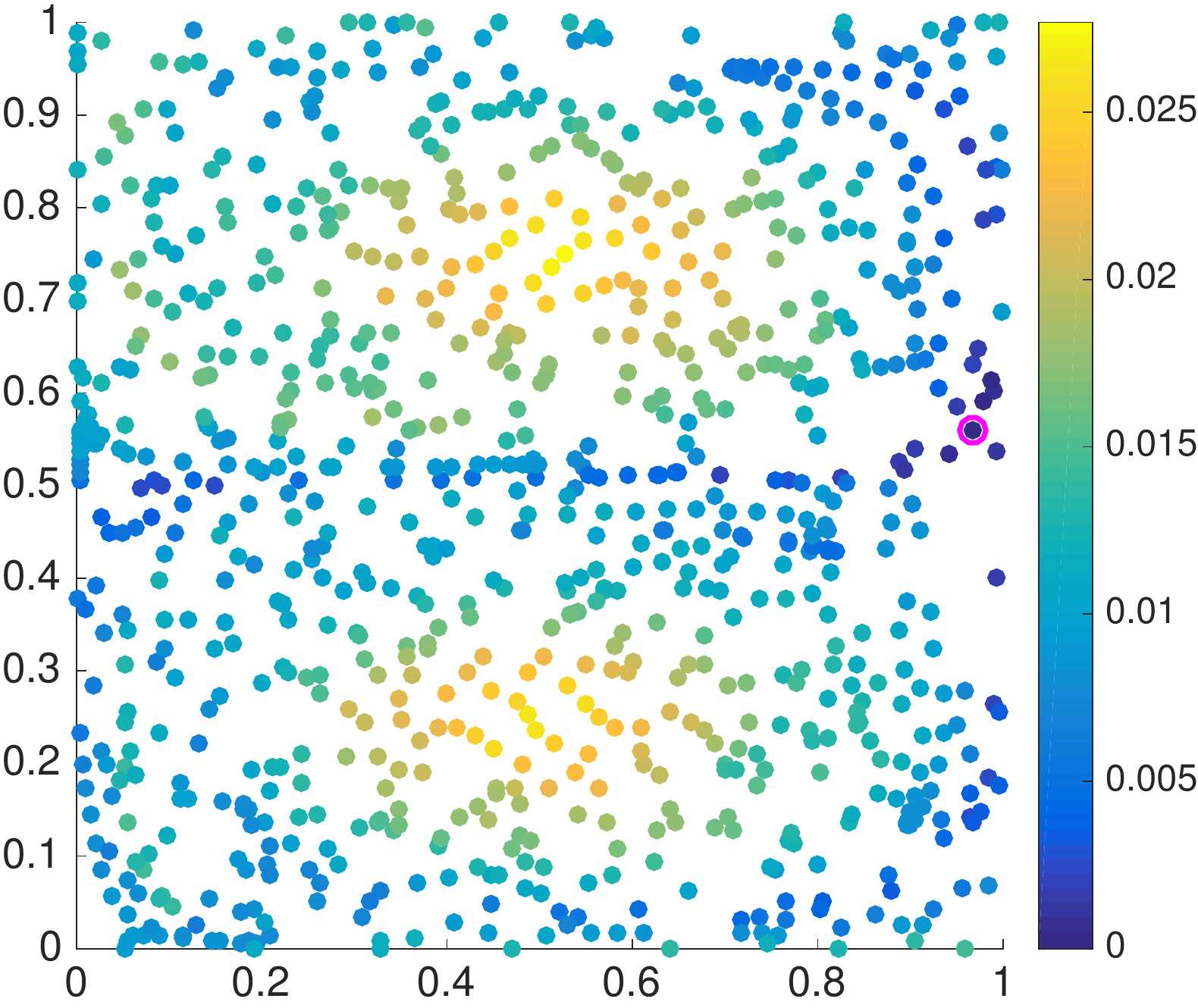}

\caption{From left to right: the identified corner stores~$c_i$, $i=1,\ldots,3$, (magenta circles) and their distances~$\nu^{\cross}(c_i,\cdot)$ to the other trajectories, at initial time (top) and final time (bottom). The distances are given in units~$1/\tau$. The horizontal axis is~$x_1$, the vertical is~$x_2$.}
\label{fig:TransDG_crossdist_cores}
\end{figure}
\begin{figure}[htb]
\centering
\includegraphics[width = 0.32\textwidth]{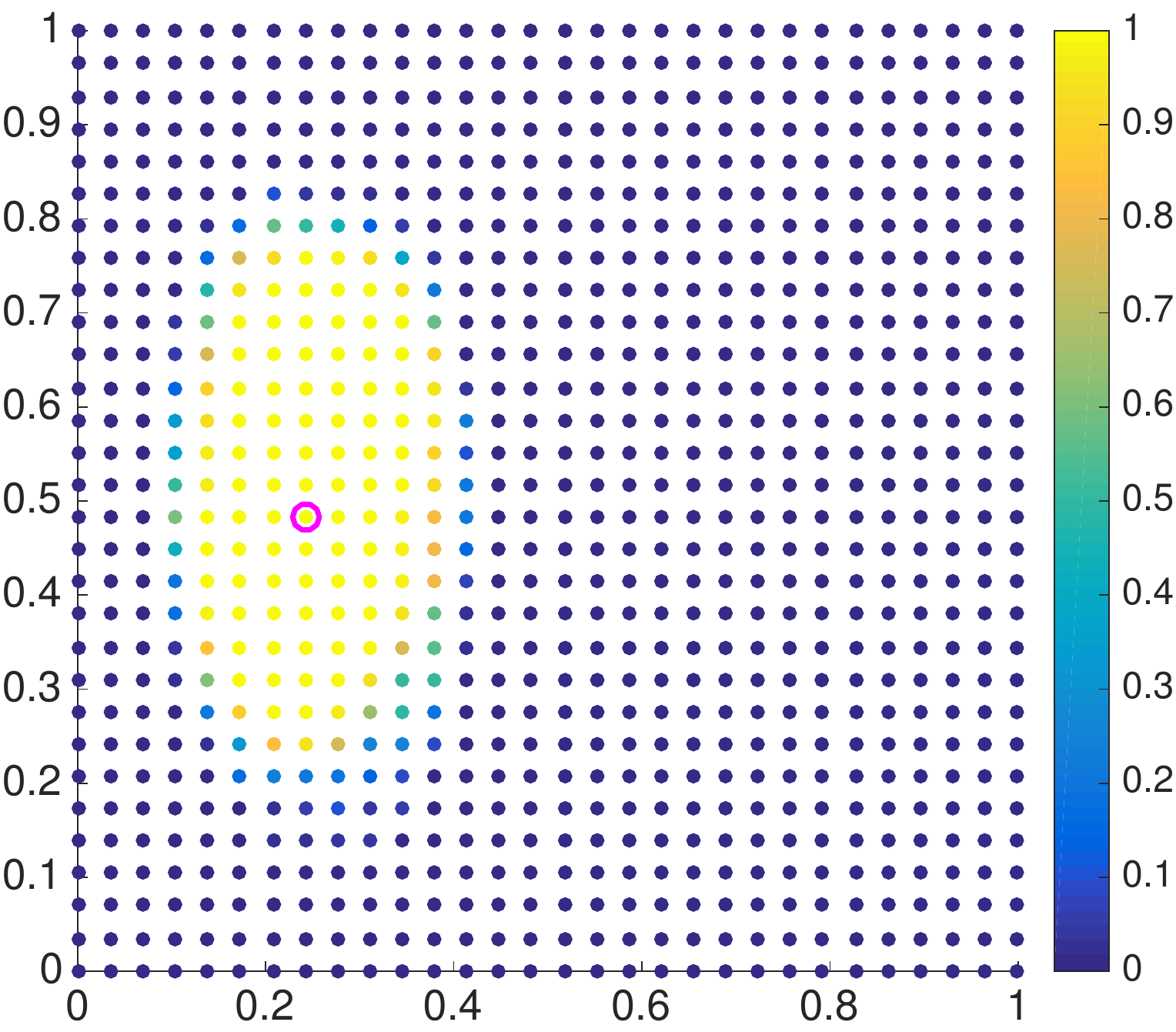}
\hfill
\includegraphics[width = 0.32\textwidth]{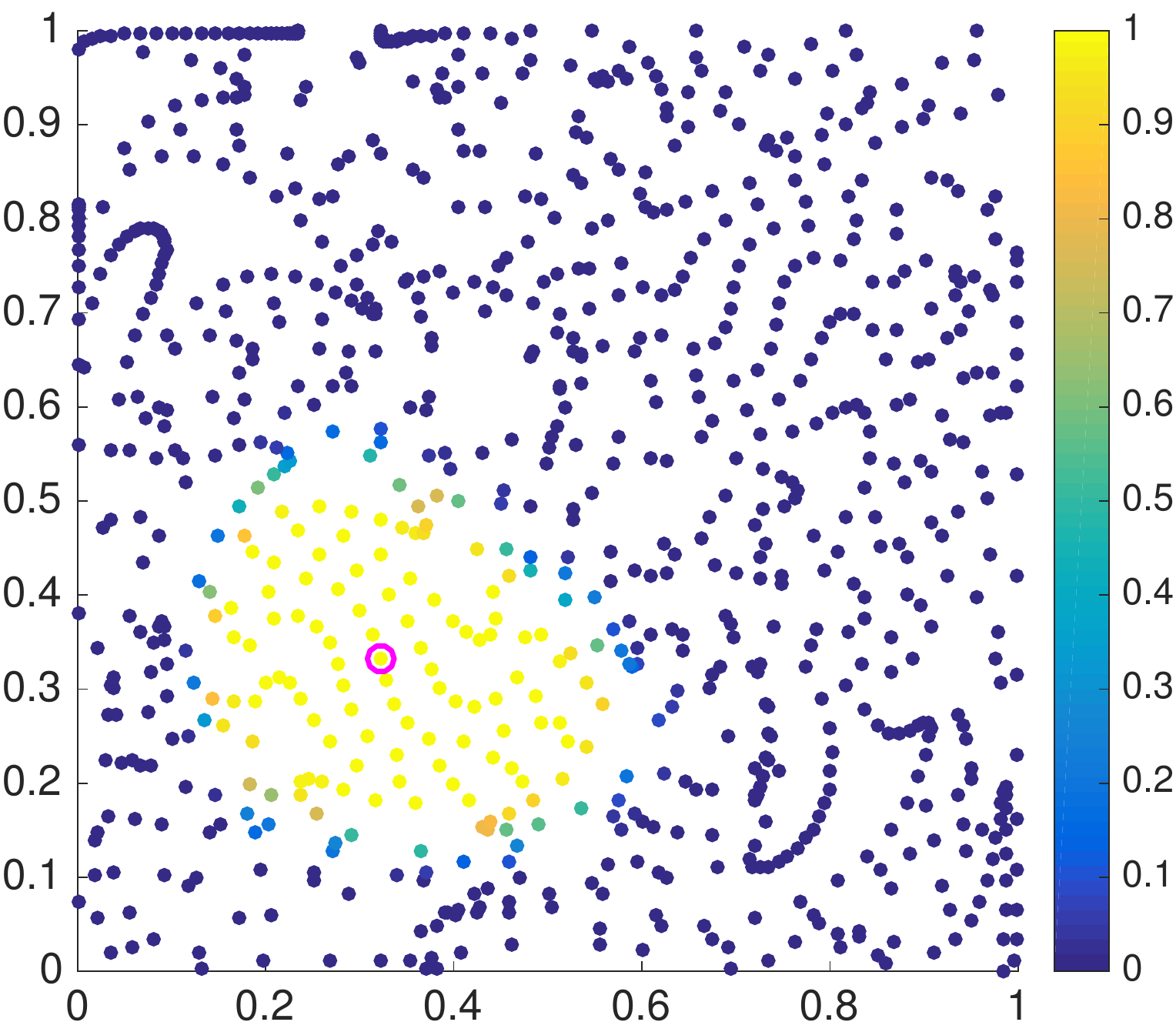}
\hfill
\includegraphics[width = 0.32\textwidth]{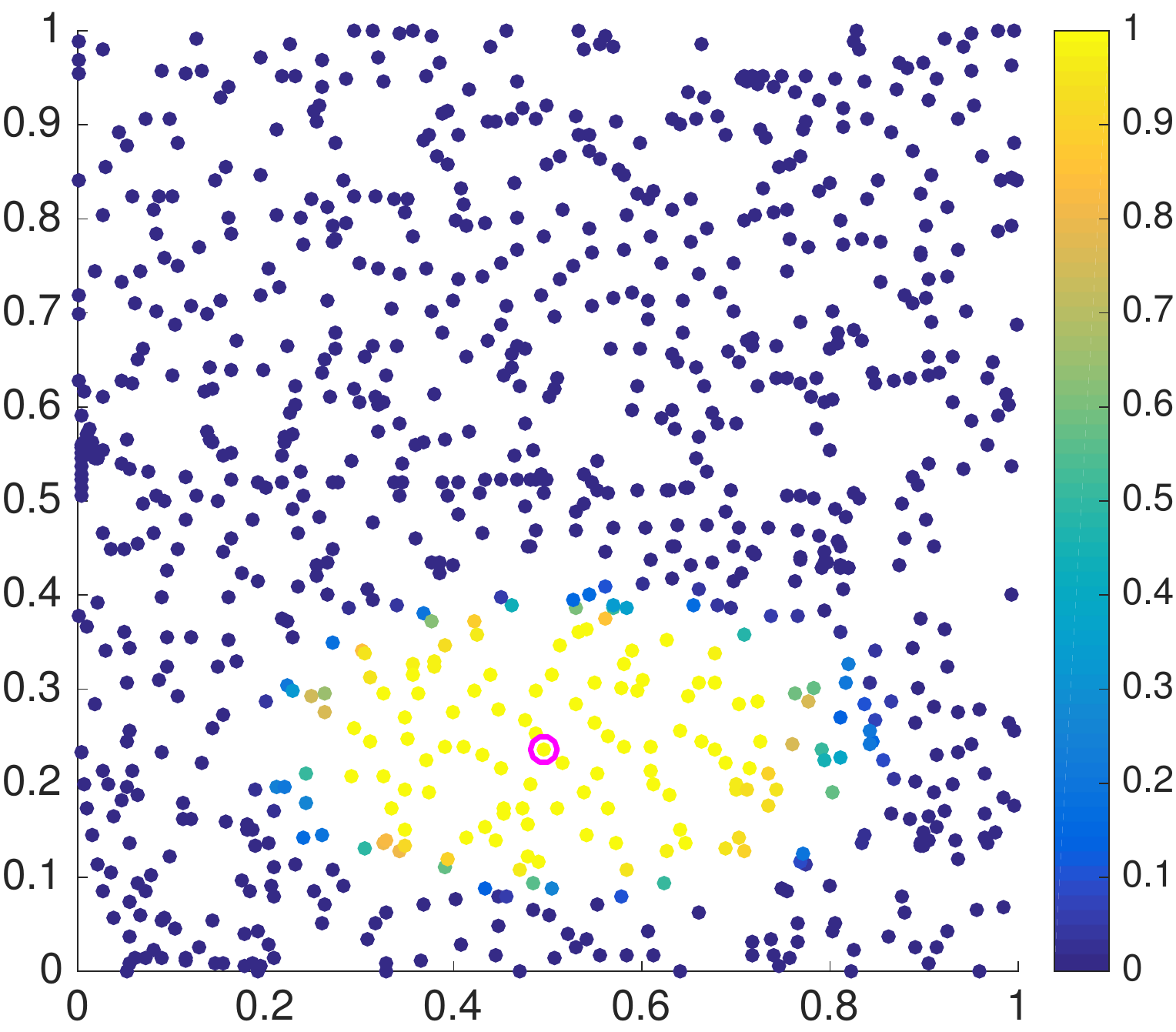}

\includegraphics[width = 0.32\textwidth]{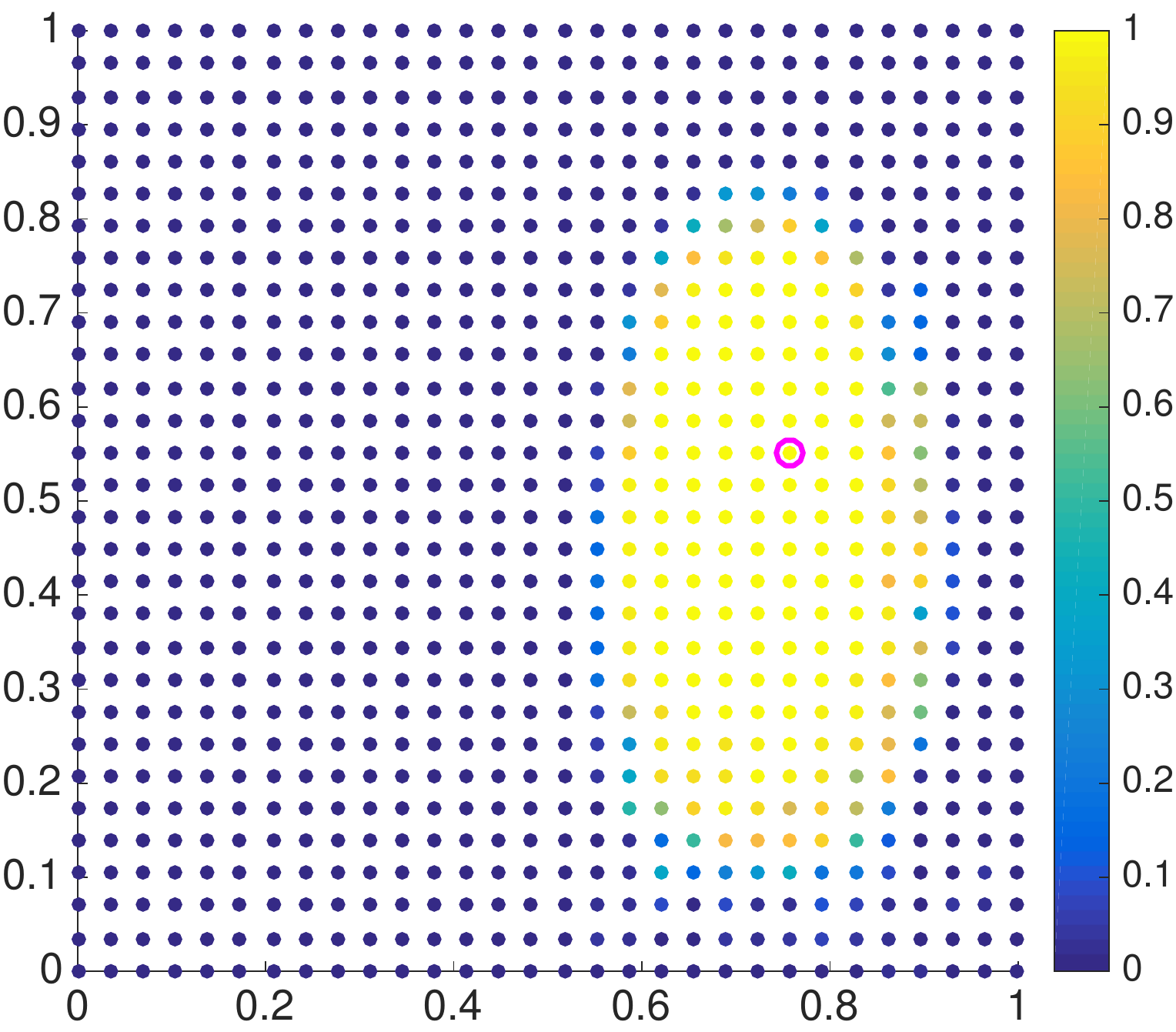}
\hfill
\includegraphics[width = 0.32\textwidth]{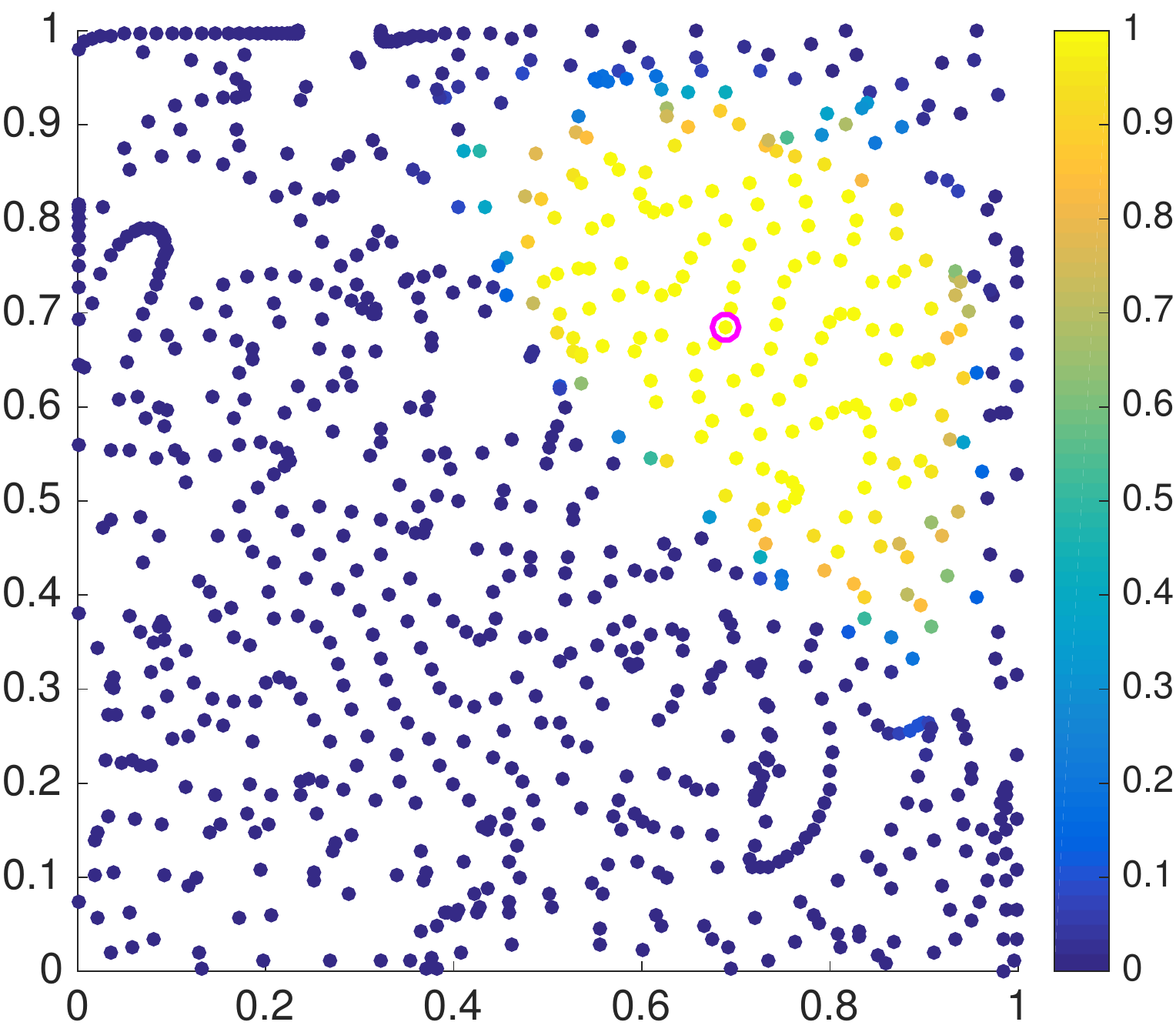}
\hfill
\includegraphics[width = 0.32\textwidth]{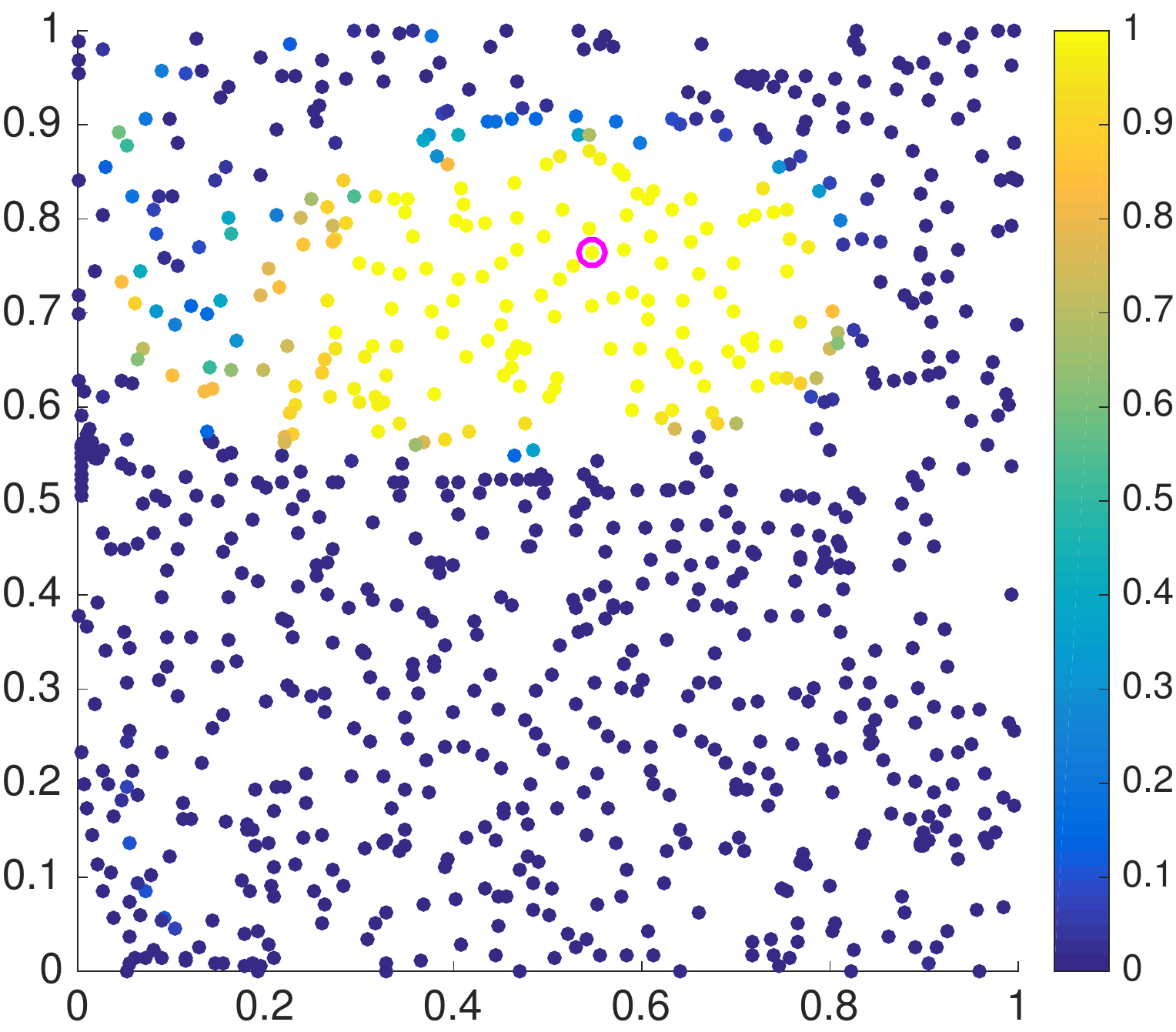}

\caption{The fuzzy affiliations computed with~$m=1.2$ to the cornerstones~$c_1$ (top) and~$c_2$ (bottom), at times $t=0,0.5,1$, from left to right, respectively. The horizontal axis is~$x_1$, the vertical is~$x_2$.}
\label{fig:TransDG_crossdist_affiliations}
\end{figure}
The semidistances from the three identified cores and the affiliations to these cores for exponent~$m=1.2$ are shown in figures~\ref{fig:TransDG_crossdist_cores} and~\ref{fig:TransDG_crossdist_affiliations}, respectively. Although both~$\nu^\cross$ and~$\nu^\meet$ have been shown to be able to detect coherent sets, we demonstrate their different nature by showing the shortest paths in the respective distance in Figure~\ref{fig:TransDG_paths_c1c2}.
\begin{figure}[htb]
\centering
\includegraphics[width = 0.32\textwidth]{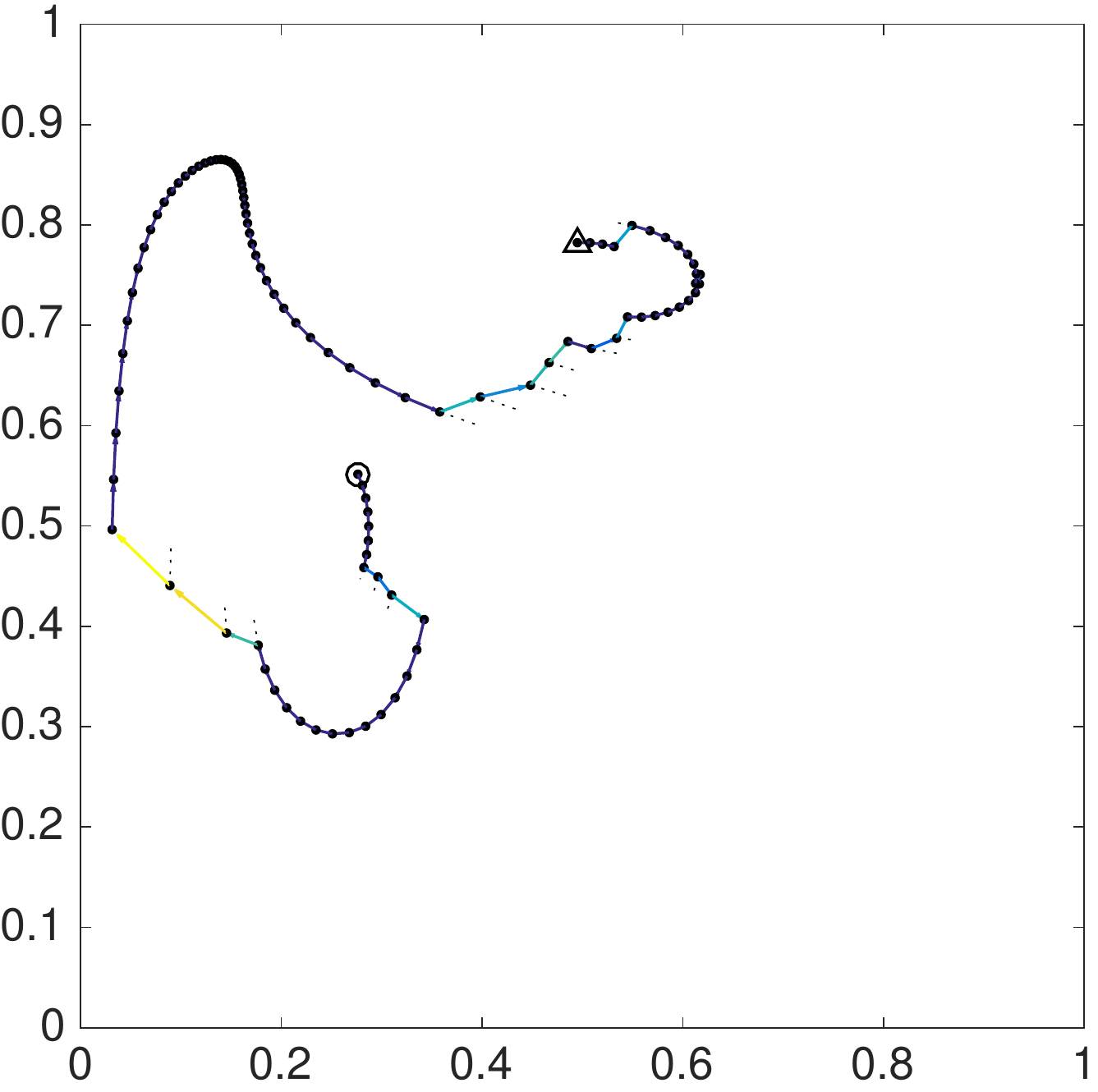}
\hfill
\includegraphics[width = 0.32\textwidth]{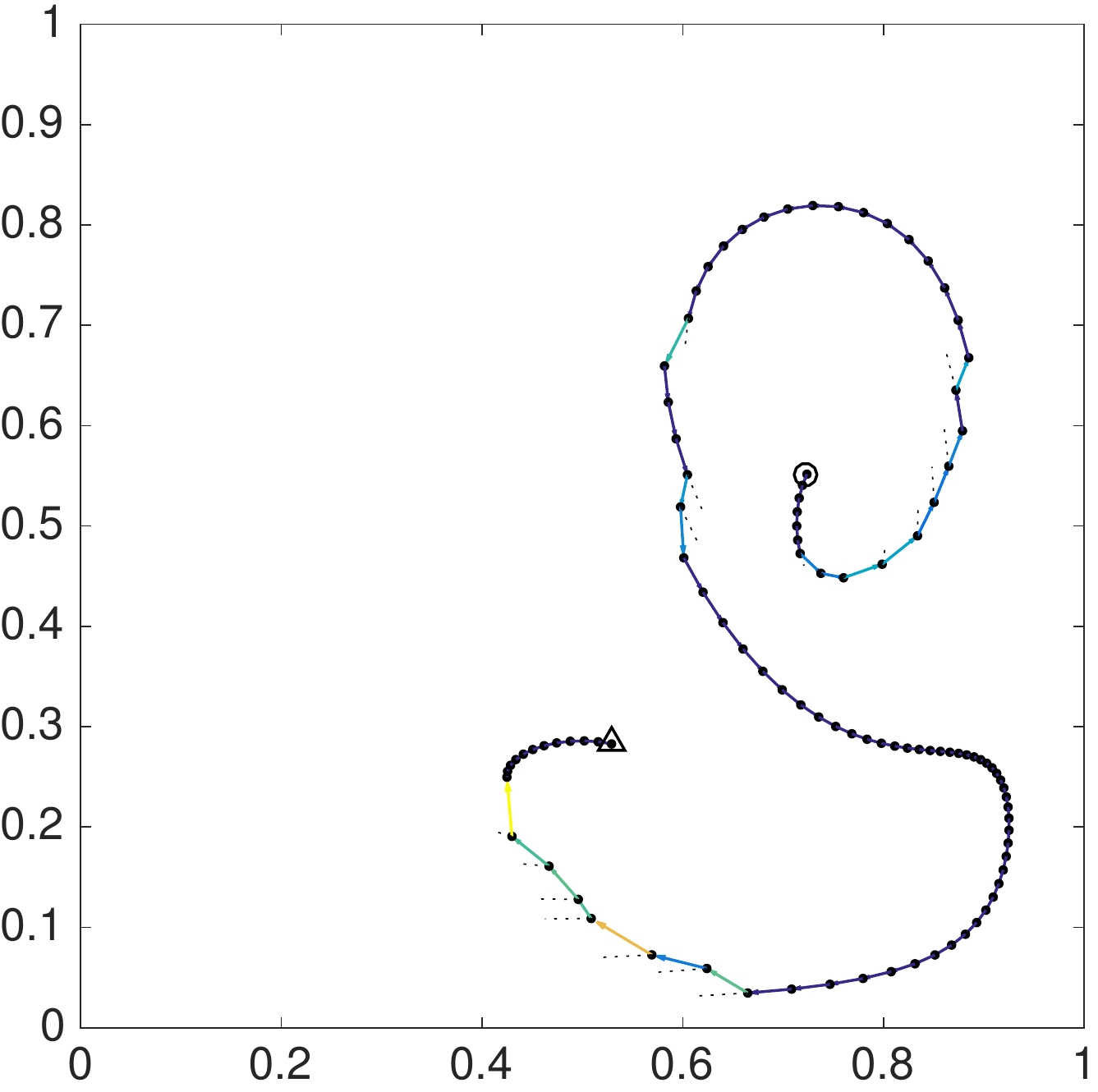}
\hfill
\includegraphics[width = 0.32\textwidth]{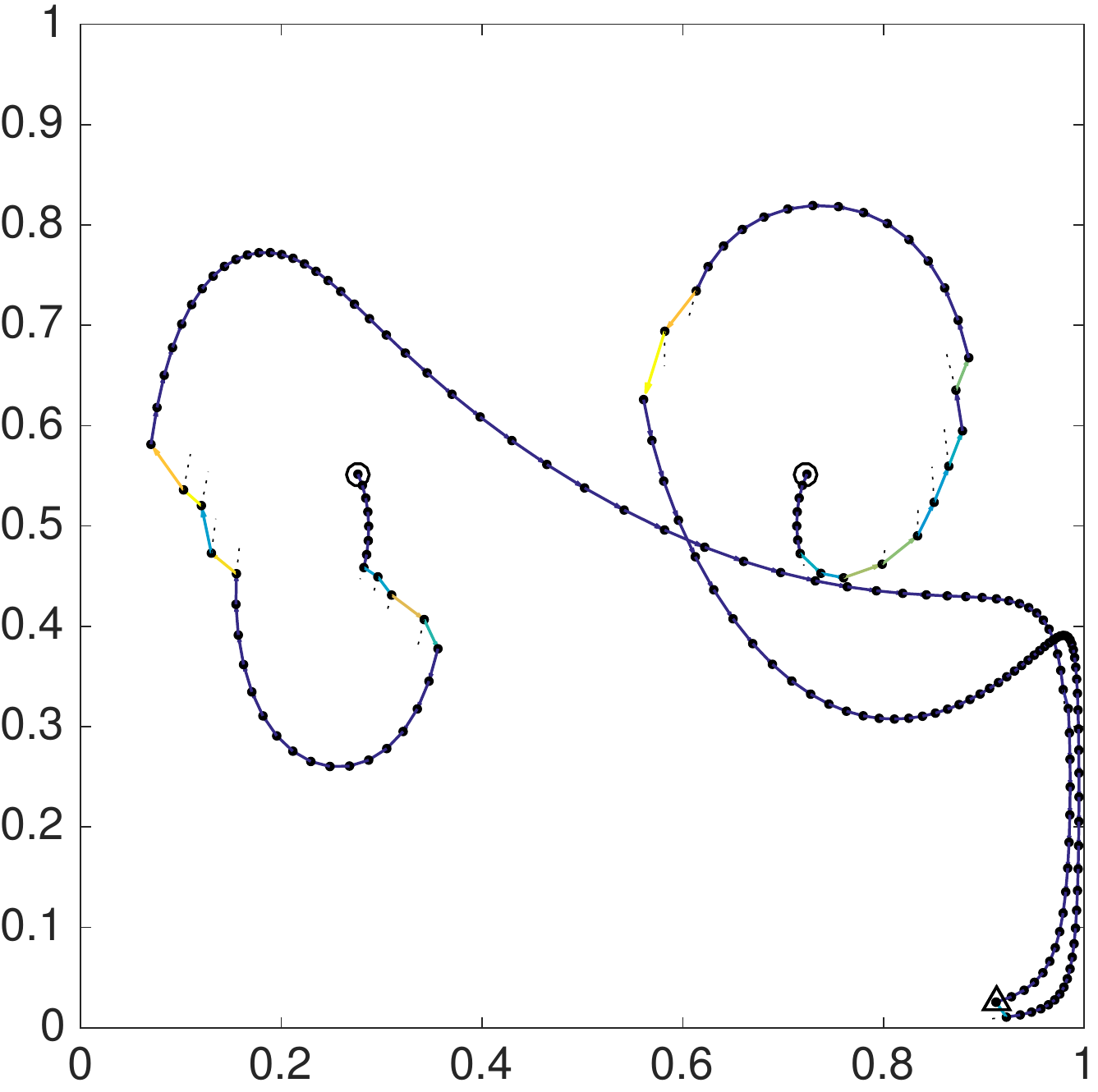}

\caption{Shortest paths from cornerstone~$c_1$ to~$c_2$ (left), from~$c_2$ to~$c_1$ (middle), and the meeting paths with the shortest joint length (right). The brighter the color of a path segment, the larger the cost of that transition. For those segments for which the path crosses from one trajectory to another we show with a dashed segment how the former trajectory would have continued. The starting point of the path is indicated by a circle, and the endpoint by a triangle. The horizontal axis is~$x_1$, the vertical is~$x_2$.}
\label{fig:TransDG_paths_c1c2}
\end{figure}

\begin{figure}[h]
\centering
\includegraphics[width = 0.32\textwidth]{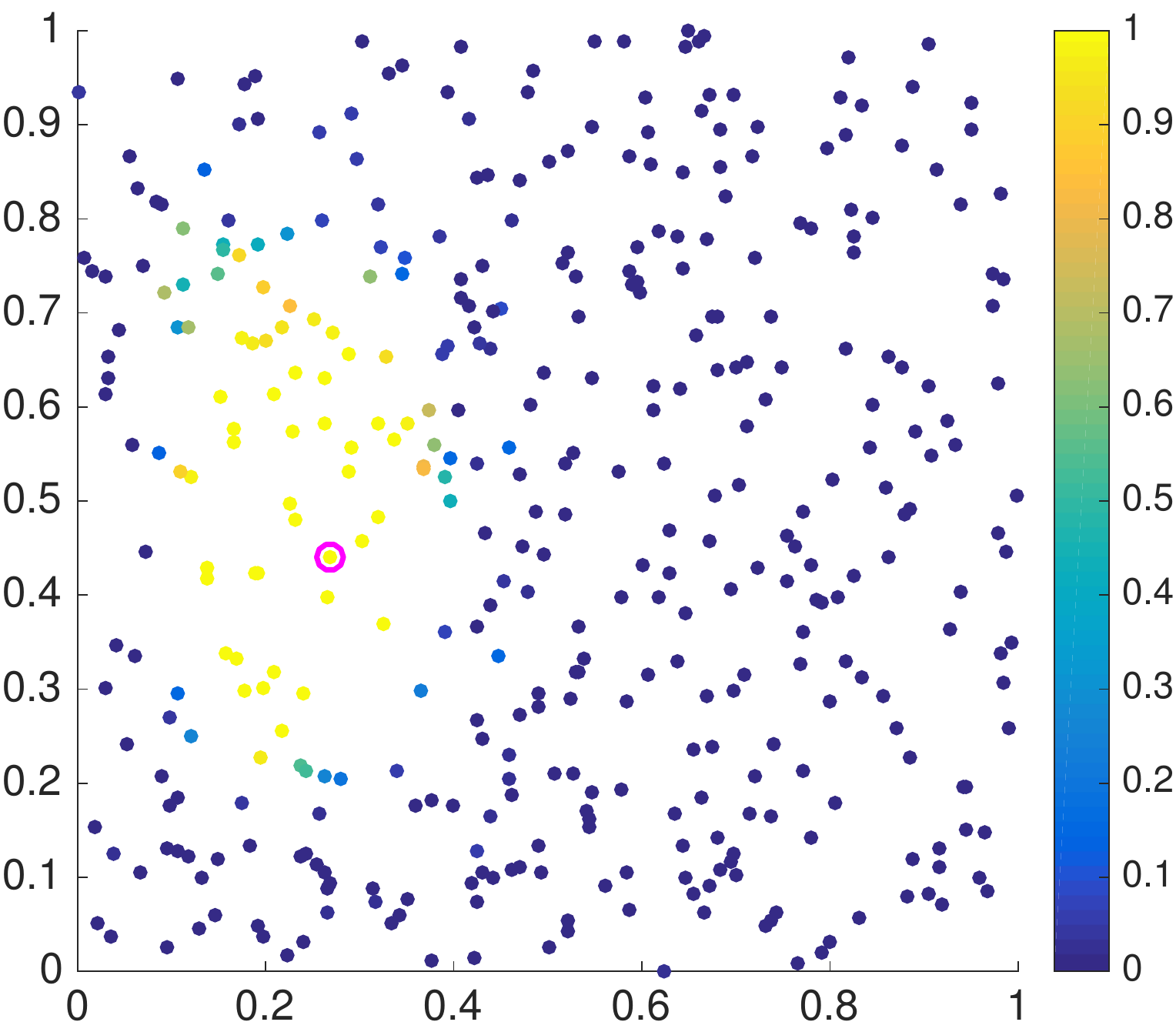}
\hfill
\includegraphics[width = 0.32\textwidth]{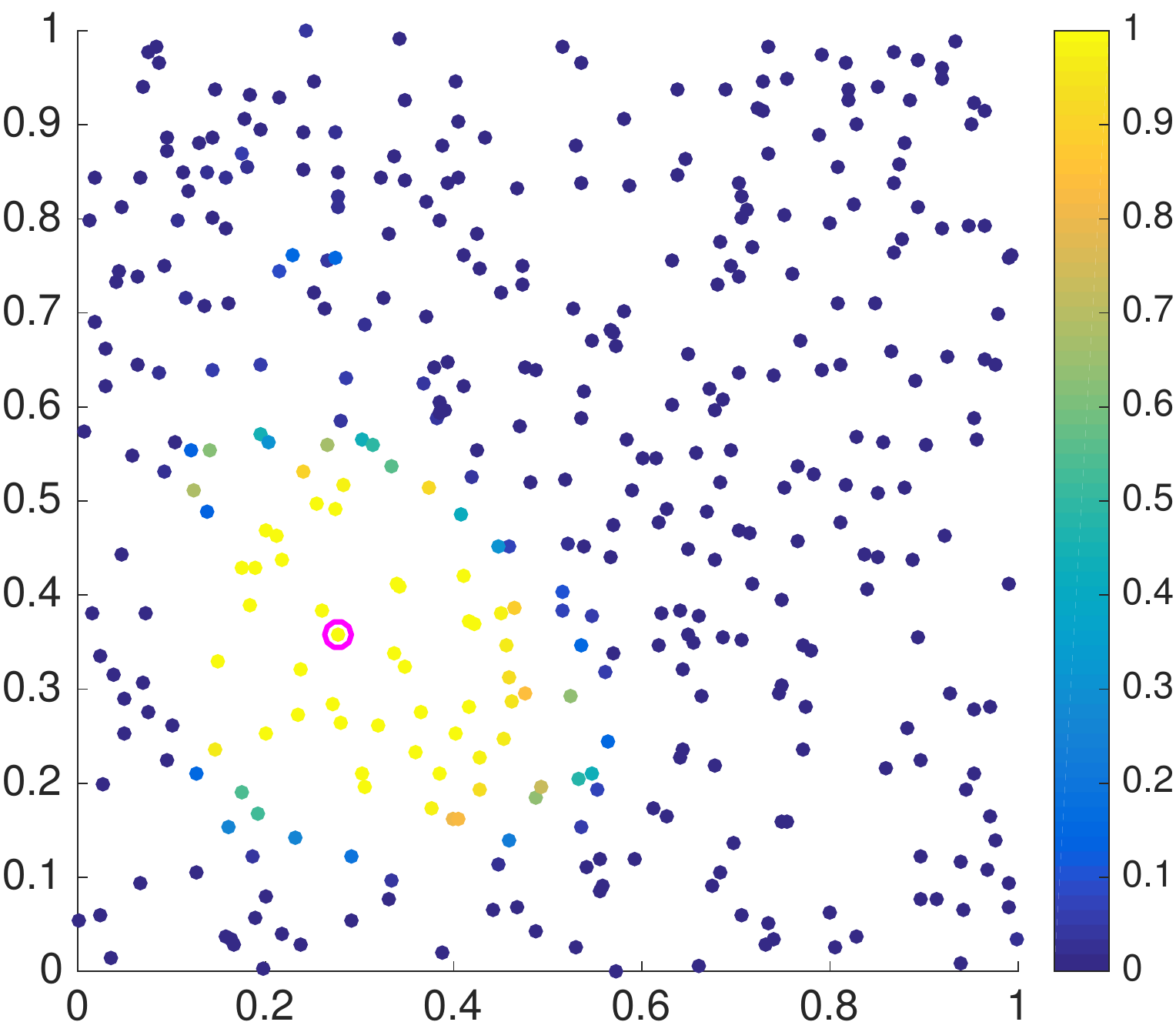}
\hfill
\includegraphics[width = 0.32\textwidth]{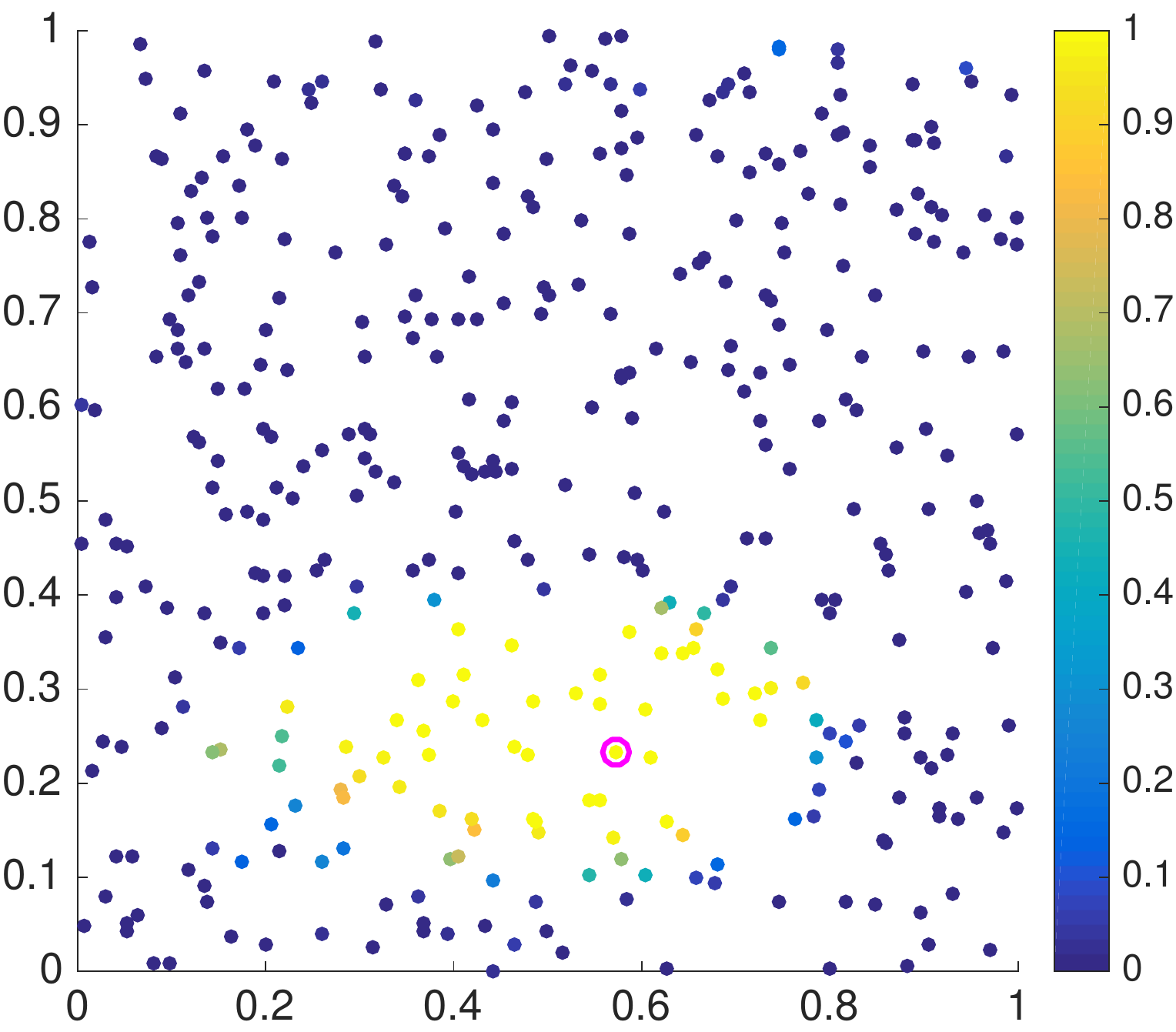}

\includegraphics[width = 0.32\textwidth]{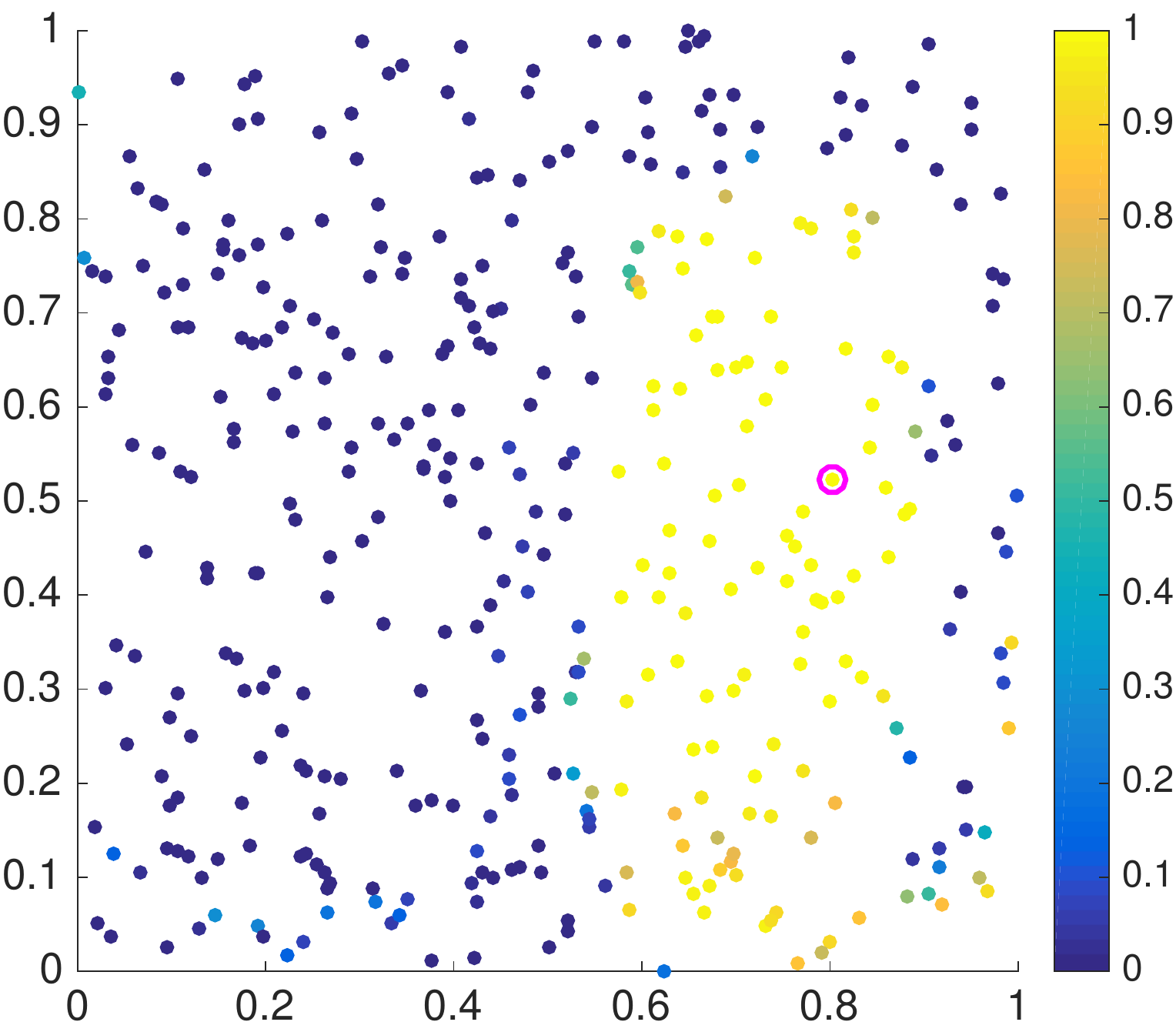}
\hfill
\includegraphics[width = 0.32\textwidth]{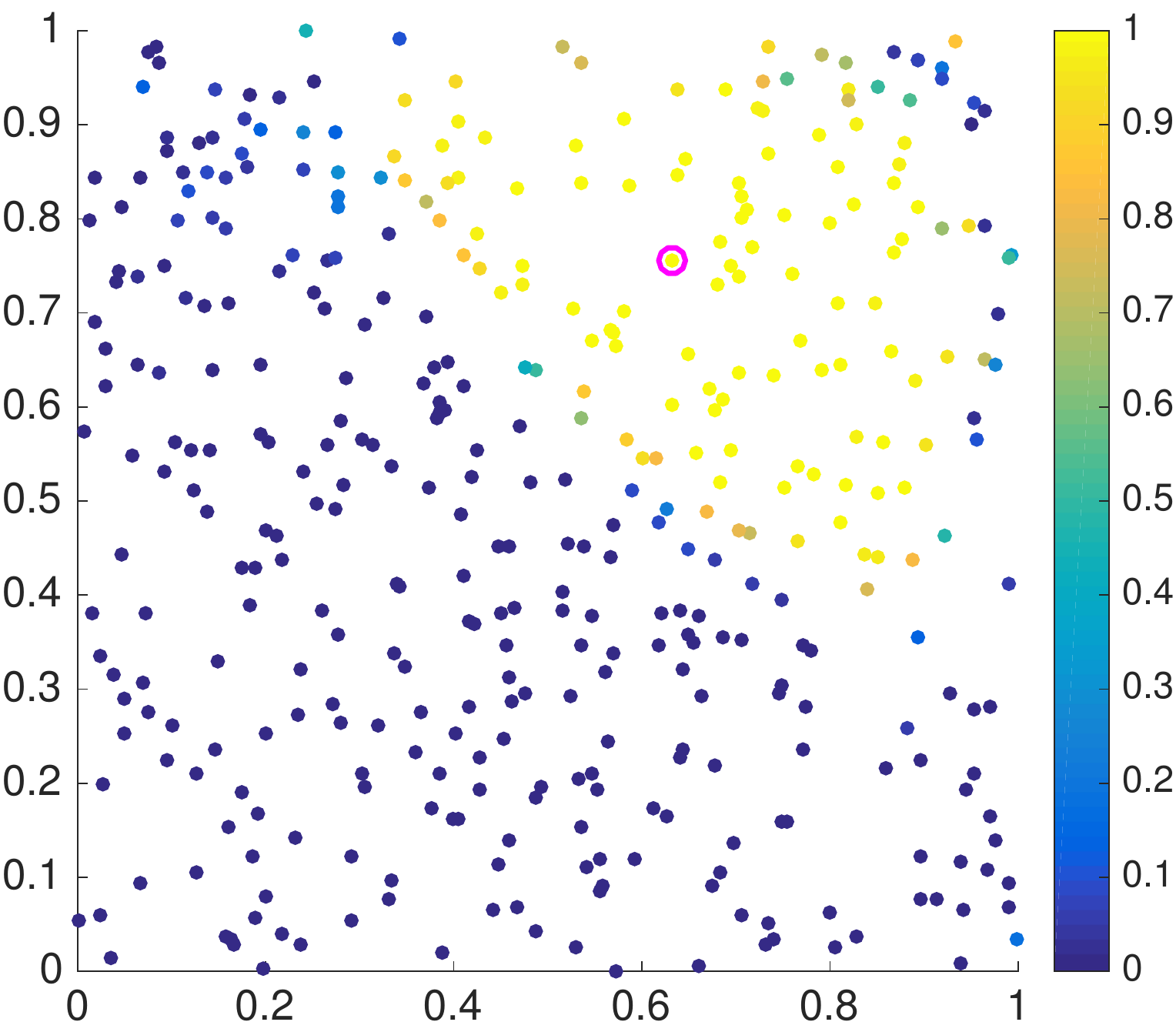}
\hfill
\includegraphics[width = 0.32\textwidth]{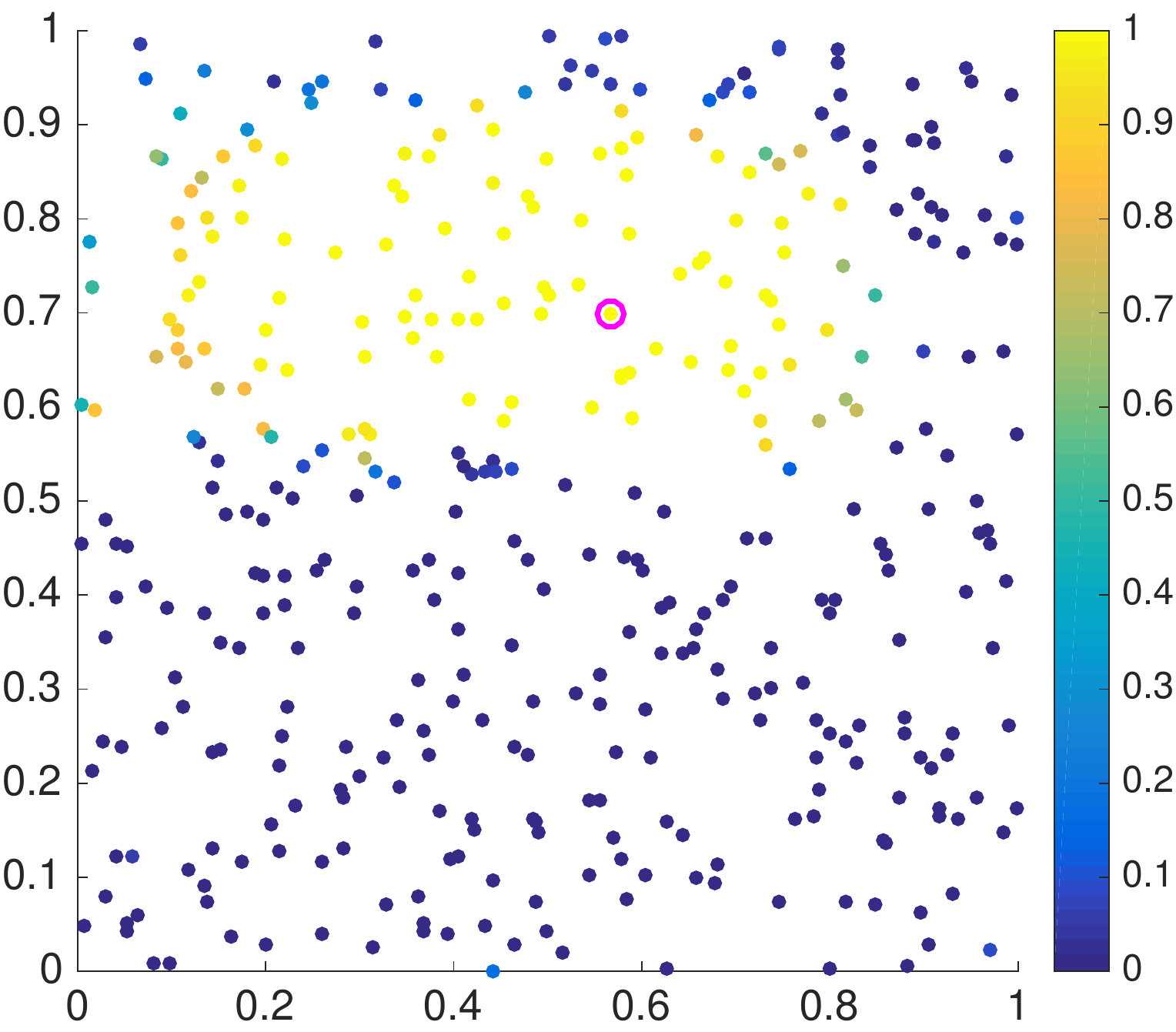}

\caption{The rotating double gyre with randomly chosen~$400$ initial points. The fuzzy affiliations computed with~$m=1.2$ to the cornerstones~$c_1$ (top) and~$c_2$ (bottom), at times $t=0,0.5,1$, from left to right, respectively. The horizontal axis is~$x_1$, the vertical is~$x_2$.}
\label{fig:TransDG_crossdist_affiliations_randpoints}
\end{figure}

Finally, we demonstrate the approach for a scattered set of sparse data points, taking~$400$ initial points randomly distributed in~$X$, and repeating the analysis for their trajectories. We show the resulting fuzzy affiliations in Figure~\ref{fig:TransDG_crossdist_affiliations_randpoints}.

\section{Discussion and outlook}
\label{sec:outlook}

\subsection{The dynamic Laplacian}
\label{ssec:dynLap}

Froyland~\cite{Froyland2015} has introduced the \emph{dynamic Laplacian} as a transport-related tool to find coherent sets. Similarly to our approach, it makes use of a small random perturbation of size~$\ep$, then~$\ep$ is driven to zero.

Numerical methods so far discretize directly the dynamic Laplacian~\cite{FrJu15,BaKo17,FrJu17}. In light of our analysis, which can be used both ways (derive the large-deviation principle in continuous space, then discretize it to finite trajectories, cf.\ Section~\ref{ssec:ldp then discretize}, or discretize the dynamics to finite trajectories, then derive the large-deviation principle on them, cf.\ Section~\ref{subsec:discr traj ldp}), we ask whether there is a discrete dynamic Laplacian that can be derived from a discretization of the perturbed dynamics?

Mimicking the construction in~\cite{Froyland2015} and sketching the idea while skipping details, one should construct a discrete,~$\ep$-dependent transfer operator~$T_{\ep}\in \R^{I\times I}$, that represents transition probabilities of a forward-backward process, then obtain a discrete dynamic Laplace operator~$L_{\text{dyn}} := \frac{d}{d\ep}\big\vert_{\ep=0}T_{\ep}$. A discrete transfer operator~$T_\ep$ that is a consistent approximation of the continuous dynamics can be obtained by a construction as in Section~\ref{subsec:discrete time and space}, by using the transition probabilities~\eqref{eq:onejump}. Technical details aside, we see that the probabilities are linear combinations of terms of the form~$\text{e}^{-\Delta x / \ep}$, where~$\Delta x$ here is a formal distance term that appears in the formulas. Differentiation with respect to~$\ep$ immediately yields that all off-diagonal entries (basically, where~$\Delta x>0$) of~$L_{\text{dyn}}$ are zero, in fact the matrix is the identity.

Thus, this approach of discretizing the dynamics first, and then factoring out the $\ep$-small stochastic perturbation does not give a dynamically meaningful result. In analytic terms the very same problem occurred in a different attempt to introduce a discrete dynamic Laplacian from a discrete transfer operator, see~\cite[Section~IV]{BaKo17}. In general, it would be desirable to understand when and how can the ``first discretize, then factor out~$\ep$'' methods work, such that they can complement the methods that directly discretize the (continuous) dynamic Laplace operator.

\subsection{Other distance measures}

The time-dependent shortest path problem used to compute our semidistances is computationally demanding in our current algorithmic realization, which theoretically limits the number of trajectories that can be handled. Moreover, they do not satisfy the triangle inequality, hence they are not a metric. Although numerical efficiency is not the main focus of this paper, and we demonstrated the usefulness of our semidistances in unraveling the underlying dynamical structure of the example systems, a more cheaply computable metric would enhance the utility and significance of the analysis methods presented here.

Ultimately, one would like to understand the intrinsic, possibly low-dimensional geometric organization of the state space with respect to transport and mixing, as pioneered in~\cite{BaKo17}. Employing proper metrics would allow, e.g., the usage of low-dimensional embedding techniques, such as multidimensional scaling, to represent and better understand this geometric organization. One canonical candidate would be the metric structure related to the dynamic Laplacian, considered in~\cite{KarraschKeller2016}. This will be subject of future studies.

To summarize, although other distance measures could be used to analyze complicated dynamic behavior, we showed that the semidistances we derived in this paper from the physical notion of transport and mixing in the vanishing diffusion setting are natural and effective diagnostic tools.

%Ultimately, one would like to understand the intrinsic, possibly low-dimensional geometric organization of the state space with respect to transport (or mixing), as pioneered in~\cite{BaKo17}. Employing proper metrics would allow, e.g., the usage of low-dimensional embedding techniques, such as multidimensional scaling, to represent and better understand this geometric organization. One canonical candidate would be the metric structure related to the dynamic Laplacian, considered in~\cite{KarraschKeller2016}. This approach seems to require a genuine distance rather than a semidistance.
%
%Indeed, the semidistances introduced in this paper do not satisfy the triangle inequality. Moreover, they have the practical disadvantage of being computationally demanding, which theoretically limits the number of trajectories that can be handled. In this respect a more practical candidate could be the $L^2$-distance, as we mentioned in Remark~\ref{rem:other semidistances}\footnote{\red{Or the dynamic Laplacian?}}. However, numerical efficiency is not the main focus of this paper, and we demonstrated the usefulness of the two transport semidistances in unraveling the underlying dynamical structure of the example systems. The main point is that the two semidistances derived in this paper can be directlt interpretated as large-deviation transport costs.

\section*{Acknowledgments}

This work is supported by the Deutsche Forschungsgemeinschaft (DFG) through the Priority Programme SPP 1881 ``Turbulent Superstructures'', and through the CRC 1114 ``Scaling Cascades in Complex Systems'', projects A01 and C08.

\appendix

\section{Large deviations of the forward-backward conditions}
\label{app:two-time LDP}

In this appendix we explore the conditions~\eqref{eq:forward backward conditions} in the large-deviation regime. The argument is based on the \emph{Laplace Principle}, which states that for any measure $\rho$ and function $f$:
\begin{equation}
  \lim_{\ep\to0}-\ep\log\int\!\e^{-\tfrac{1}{\ep} f(x)}\,\rho(dx) = \inf_{x\in\supp{\rho}} f(x).
\label{eq:Laplace}
\end{equation}
As in \eqref{eq:cont ldp one traj explicit},
\begin{equation}
  -\ep\log\prob\big[ \rv{x}_T\super{\ep}\asymp y \mid \rv{x}_0\super{\ep}=x\big] \xrightarrow[\ep\to0]{} \inf_{x\suber{\cdot}\,: x_0=x,x_T=y}\, \mfrac12\int_0^T\!\lvert \dot x_t-v(t,x_t)\rvert^2\,dt =: \lambda_T(x\smallto y),
\label{eq:cont ldp one traj explicit 2}
\end{equation}
where, contrary to \eqref{eq:cont ldp one traj explicit}, the symbol $y$ now denotes a position at time $T$, that is, $\mu_T(x\smallto \tilde x)=\lambda_T(x\smallto\phi_{0,T}\lbrack \tilde x\rbrack)=\lambda_T(x\smallto y)$.

Fix an $\ep$-independent initial probability measure $\rho_0(dx)=\prob[\rv{x}_0\in dx]$. For the large deviations of the forward condition in~\eqref{eq:forward backward conditions}, it follows from the Laplace principle that
\begin{align*}
  \J^\fw_T(B|A) &:=\lim_{\ep\to0}-\ep \log\prob\big[ \rv{x}_T\super{\ep}\in B \mid \rv{x}_0\in A \big] \\
    &= \lim_{\ep\to0}-\ep\log\int_B\int_A\!\prob\big[ \rv{x}_T\super{\ep}\in dy \mid \rv{x}_0=x\big] \rho_0(dx) \\
    &\stackrel{\eqref{eq:cont ldp one traj explicit 2}}{=}  \lim_{\ep\to0}-\ep\log\int_B\int_A\! \e^{-\tfrac{1}{\ep}\lambda_T(x\smallto y)}\,\rho_0(dx) \\
    &\stackrel{\eqref{eq:Laplace}}{=} \inf_{y\in B}\, \inf_{x\in A\cap\supp\rho_0}\,\lambda_T(x\smallto y).
\end{align*}
Observe that since the initial distribution $\rho_0$ is independent of $\ep$, it only appears in the large deviations through its support $\supp\rho_0$.

The large deviations of the backward conditions in~\eqref{eq:forward backward conditions} can be calculated analogously, but now the conditioning does depend on $\ep$. By Bayes' rule, the rate function of the backward condition is
\begin{align*}
  \J^\bw_T(A|B) & := \lim_{\ep\to0}-\ep \log\prob [\rv{x}_0\in A\mid \rv{x}_T\super{\ep}\in B] \\
    &= \lim_{\ep\to0}-\ep\log\prob[\rv{x}_T\super{\ep}\in B \mid \rv{x}_0\in A] \frac{\prob[\rv{x}_0\in A]}{\prob[\rv{x}_T\super{\ep}\in B]} \\
    &= \lim_{\ep\to0} -\ep\log\int_B\int_A\!\prob\big[ \rv{x}_T\super{\ep}\in dy \mid \rv{x}_0=x\big] \rho_0(dx) \\
    &\hspace{32pt} + \ep\log\int_B\int\!\prob\big[ \rv{x}_T\super{\ep}\in dy \mid \rv{x}_0=x\big] \rho_0(dx) - \ep\log\rho_0(A) \\
    &\stackrel{(\ref{eq:Laplace},\ref{eq:cont ldp one traj explicit 2})}{=} \inf_{y\in B}\, \inf_{x\in A\cap\supp\rho_0}\,\lambda_T(x\smallto y) - \inf_{y\in B}\, \inf_{x\in \supp\rho_0}\,\lambda_T(x\smallto y)\\
    &=\J^\fw_T(B|A) - \J^\fw_T(B|X).
\end{align*}
If we assume that that there is at least one admissible path~$x_{(\cdot)}$ that starts in~$\supp\rho_0$ and ends in~$B$, then in fact~$\J^\fw_T(B|X)=0$, and so~$\J^\fw_T(B|A)=\J^\bw_T(A|B)$. 

These calculations have two important implications. First, observe that while the forward and backward probabilities~$\prob[\rv{x}_T\super{\ep}\in B \mid \rv{x}_0\in A]$ and~$\prob[\rv{x}_0\in A\mid \rv{x}_T\super{\ep}\in B]$ are not equal in general, the forward and backward rate functions are. The same argument even holds if we shrink the sets~$A$ ad~$B$ down to single points~$x$ and~$y$; in that case we obtain for the ``backward rates'' that
\[
  \lambda_T(x{\scriptstyle \leftarrow} y) := \lim_{\ep\to 0} -\ep\log\prob\big[ \rv{x}_0 \asymp x \mid \rv{x}_T\super{\ep} = y\big] = \lambda_T(x\smallto y)\,,
\]
cf.\ Remark \ref{rem:cont time-reversal}.
Apparently, in the large-deviation scaling it does not matter whether we consider the forward or the backward process. Since the forward condition $\J^\fw_T(B|A)\approx0$ in itself does not hold enough information to characterize coherence and the backward condition $\J^\bw_T(A|B)\approx0$ does not add information, these conditions are not helpful to characterize coherence.

Secondly, we see that $\J^\fw_T(B|A)=0$ as soon as $\phi_{0,T}[A\cap\supp\rho_0] \cap B \neq \emptyset$. Naturally, there are many such pairs $A,B$, and the set function $\J^\fw_T$ does not give any quantitative information about which pairs are more coherent than others. Because of this, the large deviations of the forward and backward conditions~\eqref{eq:forward backward conditions} are even less useful to identify coherent sets.

We start to gain useful information about coherence, if there are at least two coherent pairs, say~$A_1,B_1$ and~$A_2,B_2$. Then the rates~$\J^\fw_T(B_2|A_1)$ and~$\J^\fw_T(B_1|A_2)$ are in general large, since coherence of the respective set pairs dictate that it is very unlikely to encounter paths from pair~$\#1$ to pair~$\#2$. Using these rates as measures of farness is one idea this paper exploits.

%Even without this assumption $\J^\fw_T(B|A)\geq\J^\bw(A|B)$, and so the backward condition is not restrictive. 
%Recall from \eqref{eq:forward backward conditions} that the two probabilities should both be large for a coherent set at two different times. This implies that the large-deviation scaling of these probabilities $J^\fw_T(B|A)$ and $J^\bw_T(A|B)$ should both be small. However, by the above calculation, we see that if $J^\fw_T(B|A)$ is small than so is $J^\bw_T(A|B)$. Hence, in the large-deviation scaling, the backward condition does not contain more information that can be used to characterise coherent sets (unless we state how much smaller $J^\bw_T(A|B)$ should be than $J^\fw_T(B|A)$).

\section{Algorithm: shortest path in time-dependent graphs}
\label{app:algo}

Since we could not find an algorithm suited to our purpose\footnote{It should be mentioned here that every time-dependent shortest path problem can be rephrased as a time-independent problem by considering each node at each time point as a distinct node of a large graph, and then it could be solved by standard methods. We do not take this approach here, as it might introduce memory requirement issues.}, we describe in this appendix a solution we came up with to solve the problem of finding shortest paths in a graph with time-dependent non-negative edge weights. Transition is only possible between nodes that are connected by an edge of positive weight.

There are several solutions to the shortest path problem for \emph{time-independent} graphs, such as \emph{Dijkstra's algorithm}~\cite{Dij59} or the \emph{Floyd--Warshall algorithm}~\cite{Flo62}. Each of them use in some sense a ``monotonicity'' argument, namely, that sub-paths of shortest paths are shortest paths themselves. This does not hold for time-dependent graphs, because at every step that we make the environment might change completely, and the number of steps we can make is limited by the number of time instances of the graph.

We propose the following algorithm to compute shortest paths from a specific node~$s$ to all other nodes. Note that we can stay in a node for any time at zero cost. The weight of the transition $i\to j$ at time $t$ is denoted by~$w_t(i\to j)$.

\begin{algorithm}
\caption{Shortest distance in time-dependent graphs}\label{algo1}
\begin{algorithmic}[1]
\State $R_{old} = \{s\}, R_{new} = \emptyset$ (reached states at times $0$ and $1$)
\State $\dist(s)=0$, $\dist(i) = \infty$ for $i\neq s$
\For {$t=1,\ldots,T$}
\While {$R_{old} \neq \emptyset$}
\State $v = \mathrm{arg}\max_{i\in R_{old}} \dist(i)$
\State $R_{old} \gets R_{old}\setminus \{v\}$
\For {$j:\ w_t(v\to j)<\infty$}
\If {$\dist(v) + w_t(v\to j) < \dist(j)$}
\State $\dist(j) = \dist(v) + w_t(v\to j)$
\State $R_{new} \gets R_{new}\cup \{j\}$
\EndIf
\EndFor
\EndWhile
\State $R_{old} = R_{new}$
\EndFor
\end{algorithmic}
\end{algorithm}
It is important to have the $\max$ on line 5, since if we does not start the update procedure at the node which has the maximal distance, then we might erroneously cut off nodes that could still be reached from it.

Algorithm~\ref{algo1} can clearly be extended to keep track of the shortest path as well. The distance of a node~$j$ is updated to a smaller one, whenever there is a path through some other node~$v$ that is shorter than the previous one (line 10). Hence, the new candidate shortest path is the one leading to $v$ and then jumping to~$j$ in the current time step. This is implemented in Algorithm~\ref{algo2} (line 12). Herein,~$\spath(i\smallto j)$ is the shortest path from node~$i$ to node~$j$, such that~$\spath_t(i\smallto j)$ is the node the walker resides in at time~$t=0,1,\ldots,T$ while going through the shortest path, and~$\spath_0(i\smallto j)=i$. If there is no path from~$i$ to~$j$, then~$\spath(i\smallto j)$ is the zero vector. We use~$1:k$ to denote the index set~$1,2,\ldots,k$.
\begin{algorithm}
\caption{Shortest path in time-dependent graphs}\label{algo2}
\begin{algorithmic}[1]
\State $R_{old} = \{s\}, R_{new} = \emptyset$ (reached states at times $0$ and $1$)
\State $\dist(s)=0$, $\dist(i) = \infty$ for $i\neq s$
\State $\spath(s\smallto j)=0\in\R^{T+1}$ for all $j$, $\spath_0(s\smallto s)=s$
\For {$t=1,\ldots,T$}
\While {$R_{old} \neq \emptyset$}
\State $v = \mathrm{arg}\max_{i\in R_{old}} \dist(i)$
\State $R_{old} \gets R_{old}\setminus \{v\}$
\For {$j:\ w_t(v\to j)<\infty$}
\If {$\dist(v) + w_t(v\to j) < \dist(j)$}
\State $\dist(j) = \dist(v) + w_t(v\to j)$
\State $\spath_{1:t-1}(s\smallto j) = \spath_{1:t-1}(s\smallto v)$, $\spath_{t}(s\smallto j) = j$
\State $R_{new} \gets R_{new}\cup \{j\}$
\EndIf
\EndFor
\EndWhile
\State $\spath_t(s\smallto s) = s$
\State $R_{old} = R_{new}$
\EndFor
\end{algorithmic}
\end{algorithm}

\small
\bibliographystyle{abbrv}
\bibliography{library}

\end{document}